\documentclass[oneside,reqno]{amsart}
\usepackage{amsmath,amssymb,amscd,amsthm,amsxtra,amsfonts}
\usepackage{mathtools,mathrsfs,dsfont,xparse}
\usepackage{graphicx,epstopdf}
\usepackage{multirow}
\usepackage{cases}
\usepackage{enumerate}
\usepackage{xcolor}
\usepackage{tikz,pgfplots}
\usetikzlibrary{matrix}
\usepackage{mathtools}
\usepackage[margin=1 in]{geometry}
\mathtoolsset{showonlyrefs}
\usepackage[numbers,sort]{natbib}
\usepackage[all]{xy}
\usepackage{appendix}
\usepackage{comment}
\usepackage{thmtools}

\usepackage{multicol,bbm}
\usepackage[colorlinks=true, pdfstartview=FitV, linkcolor=blue, citecolor=blue, urlcolor=blue]{hyperref}

\allowdisplaybreaks[4]

\newcommand*{\parallelogramm}{%
  \scalebox{0.6}{%
    \rlap{\rotatebox{-30}{\rule[.05ex]{.4pt}{.77em}}}%
    \kern.04em%
    \rlap{\kern.36em\raisebox{0.649519052835em}{\rule{.6em}{.4pt}}}%
    \rule{.6em}{.4pt}\kern-.04em%
    \rotatebox{-30}{\rule[.05ex]{.4pt}{.77em}}%
  }%
}

\theoremstyle{plain}
\newtheorem{theorem}{Theorem}[section]
\newtheorem{proposition}[theorem]{Proposition}
\newtheorem{lemma}[theorem]{Lemma}
\newtheorem{corollary}[theorem]{Corollary}

\theoremstyle{definition}
\newtheorem{definition}[theorem]{Definition}
\newtheorem{remark}[theorem]{Remark}
\numberwithin{equation}{section}
\numberwithin{figure}{section}

\newcommand{\R}{\mathbb{R}}

\newcommand{\Z}{\mathbb{Z}}
\newcommand{\T}{\mathbb{T}}

\newcommand{\abs}[1]{\left| #1\right|}
\newcommand{\norm}[1]{\left\|#1\right\|}
\newcommand{\inner}[1]{\left \langle#1\right \rangle}

\newcommand{\parenthese}[1]{\left(#1\right)}

\begin{document}
\title[Bilinear Strichartz Estimates and Linear Morawetz Estimates]{Bilinear Strichartz Estimates and Linear Morawetz Estimates for Schr\"odinger Operators with Potentials}

\author[Wilson and Yu]{Bobby Wilson  and Xueying Yu}

\address{Bobby Wilson
\newline \indent Department of Mathematics, University of Washington\indent 
\newline \indent  C138 Padelford Hall Box 354350, Seattle, WA 98195,\indent }
\email{blwilson@uw.edu}

\address{Xueying  Yu
\newline \indent Department of Mathematics, Oregon State University\indent 
\newline \indent  Kidder Hall 368, Corvallis, OR 97331 \indent 
}
\email{xueying.yu@oregonstate.edu}

\begin{abstract}
In dimensions three and above, we establish bilinear Strichartz estimates for Schr\"odinger operators with potentials using perturbative arguments and Morawetz estimates. Our approach does not require smallness of the potential. As part of the analysis, we also establish endpoint versions of the bilinear Christ-Kiselev lemma.
\end{abstract}



\maketitle

\setcounter{tocdepth}{1}
\tableofcontents

\parindent = 10pt     
\parskip = 8pt

\section{Introduction}

\subsection{History and Motivation}
Strichartz estimates \cite{Str77,GV85,Yaj87,KT98} are one of the most fundamental tools in the study of dispersive partial differential equations (PDEs). For the linear Schr\"odinger equation on $\R^d$, they provide spacetime bounds of the form
\begin{align}\label{eq Strichartz}
\norm{e^{it\Delta}f}_{L_t^qL_x^r(\R\times\R^d)}
\lesssim \norm{f}_{L^2(\R^d)},
\end{align}
for Schr\"odinger admissible pairs $(q,r)$ (see Proposition~\ref{prop Strichartz} for a precise statement). Combined with Duhamel's formula, these estimates have been used to establish many local and global well-posedness results for nonlinear Schr\"odinger equations.

In more delicate settings (such as at certain critical regularity, on compact or curved manifolds), however, linear Strichartz estimates alone are often insufficient. In particular, they do not capture the cancellation occurring when two wave packets have separated frequency supports. This motivates bilinear Strichartz refinements, which quantify such phenomena well. A classical estimate by Bourgain \cite{Bour98} in $\R^2$ states that
\begin{align}\label{eq Bilinear}
\norm{e^{it \Delta} f_N \, e^{it\Delta} g_L}_{L_{t,x}^2 (\R \times \R^2)} \lesssim \underbrace{\hspace{0.3cm}\frac{L^{\frac{1}{2}}}{N^{\frac{1}{2}}}\hspace{0.3cm}}_{\displaystyle =: C(N,L)} \norm{f_N}_{L^2 (\R^2)} \norm{g_L}_{L^2 (\R^2)}
\end{align}
where $f_N := P_N f$ and $g_L := P_L g$ have Fourier support at frequencies $N$ and $L$ respectively, with $L \ll N$. In general dimensions, the frequency factor $L^{\frac{1}{2}} N^{-\frac{1}{2}}$ is replaced by $C(N,L) = L^{\frac{d-1}{2}} N^{-\frac{1}{2}}$ (see Lemmas~\ref{lem Bourbi} and \ref{lem alldbourbi} for detailed statements).

For comparison, for example in $d=2$, H\"older's inequality and the linear $L_{t,x}^4$ Strichartz estimate only give
\begin{align}
\norm{e^{it \Delta} f_N \, e^{it\Delta} g_L}_{L_{t,x}^2 (\R \times \R^2)} \lesssim  \norm{e^{it \Delta} f_N }_{L_{t,x}^4 (\R \times \R^2)} \norm{e^{it \Delta} g_L }_{L_{t,x}^4 (\R \times \R^2)} \lesssim \norm{f_N}_{L^2 (\R^2)} \norm{g_L}_{L^2 (\R^2)} ,
\end{align}
which does not capture any gain from the frequency separation. The improvement in \eqref{eq Bilinear} reflects the transversality of the Schr\"odinger waves at different frequency scales, whose distinct group velocities limit their spacetime interactions.

The bilinear estimate \eqref{eq Bilinear} can also be viewed as a smoothing effect. 
Indeed, since the product $e^{it\Delta}f_N \, e^{it\Delta}g_L$ is localized at frequency $\sim N$ when $L \ll N$, we have
\begin{align}
\norm{\abs{\nabla}^{\frac{1}{2}} \parenthese{e^{it\Delta}f_N \, e^{it\Delta}g_L}}_{L_{t,x}^2(\R\times\R^2)} \lesssim \norm{f_N}_{L^2(\R^2)} \|\abs{\nabla}^{\frac{1}{2}} g_L \|_{L^2 (\R^2)} .
\end{align}
A naive fractional product rule applying to $(e^{it \Delta} f_N \, e^{it\Delta} g_L )$ would distribute the half-derivative across both factors. However, the bilinear estimate shows that ``half-derivative'' only lands on the lower-frequency solution, which suggests an underlying smoothing effect.

Bilinear refinements have been developed in several directions. They are closely connected with bilinear Fourier restriction estimates for the paraboloid; Tao \cite{Tao03} proved an essentially sharp bilinear restriction theorem for transverse frequency supports, leading to a broader range of bilinear Schr\"odinger estimates. Change-of-variables methods for bilinear estimates appear in related dispersive problems studied by Kenig, Ponce, and Vega \cite{KPV96JAMS,KPV96TAMS}, and in the Schr\"odinger setting in the work of Colliander, Delort, Kenig, and Staffilani \cite{CDKS01}. Shao \cite{Shao09} obtained sharper restriction estimates under cylindrical symmetry, while mixed-norm bilinear estimates have also been useful in critical NLS analysis \cite{Dod22}. 
Although many classical proofs rely on the spacetime Fourier transform, there are also proofs emphasizing physical-space mechanisms. Tao \cite{tao2010physical} obtained a physical-space proof of the classical bilinear Strichartz and local smoothing estimates, while Hernandez \cite{hernandez2018decomposition} used a decomposition of Schr\"odinger solutions into localized tubes. Such methods are particularly relevant to settings in which the exact Fourier representation of the free propagator is unavailable.

The bilinear theory extends beyond Euclidean space. On compact surfaces, Burq, G\'erard, and Tzvetkov \cite{BGT05} proved bilinear eigenfunction estimates, and Hani \cite{Hani12} obtained bilinear refinements of Strichartz estimates on closed manifolds. On tori, bilinear and multilinear bounds have been obtained using lattice-point counting and function spaces adapted to the linear evolution \cite{IPRT3,HTT2}; on irrational tori, Fan, Staffilani, Wang, and Wilson \cite{FSWW} employed a decoupling argument. 
Related bilinear estimates are available for balls and disks \cite{An08,SY22} and for waveguides \cite{Deng23,DFYZZ24,DDFZ24}. These results show that bilinear frequency interactions depend sensitively on the geometry and spectral structure of the underlying space.

Such bilinear estimates are particularly effective in nonlinear dispersive problems.
In Bourgain’s high–low method \cite{Bour98}, they are used to show that the Duhamel term enjoys additional regularity (compared to the initial data) and smallness in the energy space (a smoothing effect), a key step in establishing global well-posedness below the energy space. Bilinear estimates also play an important role in the $I$-method, where they help control the increment of the modified energy \cite{LY24, SY22, CKSTT08, CKSTT02, DPST07, Hani12, Schippa24}. 

In critical problems, such as the mass-critical global well-posedness theory, bilinear estimates are important ingredients for the inverse Strichartz estimate, identifying the frequency and spatial scales at which concentration occurs in the profile decomposition. 
The presence of Galilean invariance in the mass-critical setting makes such bilinear control particularly important. On compact manifolds such as the torus, where dispersion is weaker than in Euclidean space, bilinear (and sometimes trilinear) estimates are also essential in the energy-critical theory \cite{HTT1,HTT2,IPRT3,IPT3,Yue, YY20, SY22, SY20d, KV16, Herr13, HS15, PTW14, Strunk15}.

A direct motivation for the present work comes from \cite{WY23}, where a mass-critical NLS global well-posedness problem on the three-dimensional hyperbolic space was studied. Due to the mass-critical nature of the problem, a bilinear Strichartz estimate is needed in the profile decomposition argument. The authors were able to use a change-of-variables trick to translate the bilinear Strichartz estimate on $\mathbb{H}^3$ to a bilinear estimate on $\mathbb{R}^3$ with a cost of introducing a negligible potential (which can be gauged away) into the linear Schr\"odinger operator. 
Such a simplification, however, is specific to the structure of the problem and is not available for general Schr\"odinger operators.
Recent work of Gavrus, \cite{Gav26},  establishes direct bilinear estimates on hyperbolic space without using the change-of-variables trick.

This naturally raises the question of whether bilinear Strichartz estimates remain valid when the free Schr\"odinger operator is perturbed by a potential.
More generally, one would like to understand to what extent the frequency-separation gains in the free estimates persist under such perturbations. This is the main question addressed in the present work.

\subsection{Goal and Scope of the Current Work}

In this paper, we establish bilinear Strichartz estimates for a class of Schr\"odinger operators of the form
\begin{align}
H_V=-\Delta_{\R^d}+V,
\end{align}
where $d\geq 3$ and $V$ is a real-valued, time-independent potential.

Throughout the paper, we will encounter one of two types of potential: (1) the more general $(\kappa, \delta)$ potentials (see Definition~\ref{defn GoodP}) that may have a pole of order 2 and (2) smooth, bounded potentials of the form $V(x)= \inner{x}^{-\eta}$ for $\eta\geq 2$. Of course, type (2) potentials form a subset of the type (1) potentials. Type (2) potentials will be used when we need to absorb more than 1/2 derivative in the argument, for example in Theorems~\ref{thm bilocdata high modes} and \ref{thm bilocdata lowmodes}. The decay rate of $(\kappa, \delta)$ potentials, $|x|^{-2}$, is critical with respect to the scaling of the Schr\"odinger equation, since the rescaled potential
\begin{align}
V_\lambda(x)=\lambda^2V(\lambda x)
\end{align}
has the same scaling as the Laplacian. 

One of the main difficulties in extending bilinear Strichartz estimates to Schr\"odinger operators with potentials is the loss of commutation between the linear propagator and the usual Littlewood--Paley projections. In the free case, $e^{it\Delta}$ commutes with $P_N$, so it makes no difference whether the frequency localization is imposed before or after the evolution. In the presence of a potential, however, we generally have
\begin{align}
P_Ne^{-itH_V}\neq e^{-itH_V}P_N.
\end{align}
Consequently, frequency-localized initial data need not remain localized at the same frequency under the perturbed evolution.

This leads us to consider two different types of bilinear estimates:
\begin{align}\label{eq subPBilinear}
\text{Type 1:} & \quad \norm{P_Ne^{-itH_V}f \, P_Le^{-itH_V}g}_{L_{t,x}^2(\R\times\R^d)}, \\
\text{Type 2:} & \quad \norm{e^{-itH_V}P_Nf \, e^{-itH_V}P_Lg}_{L_{t,x}^2(\R\times\R^d)}.
\end{align}
The first type measures the interactions between frequency components of solutions to the perturbed equation, whereas the second concerns solutions whose initial data are localized at different frequencies.
Both types are useful in nonlinear applications. In particular, when studying frequency interactions in profile decompositions, one may need to impose frequency projections both before and after the perturbed evolution. This motivates estimates of the schematic form
\begin{align}\label{eq PBilinear}
\begin{aligned}
& \quad \norm{P_Ne^{-itH_V}(P_Mf) \, P_Le^{-itH_V}(P_Kg)}_{L_{t,x}^q(\R\times\R^d)}\\
& \lesssim C(N,L)\phi(N,M)\phi(L,K)\norm{P_Mf}_{L^2(\R^d)}\norm{P_Kg}_{L^2(\R^d)},
\end{aligned}
\end{align}
where $C(N,L)$ represents the frequency-interaction bound and $\phi$ measures the loss associated with the mismatch between the input and output frequencies.
The main goal of this paper is to establish bilinear estimates of both types and to combine them to obtain estimates of the form \eqref{eq PBilinear}.

The difficulty is that classical bilinear Strichartz estimates rely heavily on the Fourier representation of the free Schr\"odinger propagator. In particular, the gain from frequency separation is closely related to the geometry of the paraboloid and the transversality of its frequency supports. When a potential is introduced, the perturbed propagator is no longer a Fourier multiplier, and the standard Fourier analytic arguments cannot be applied directly.

Bilinear estimates for Schr\"odinger evolutions beyond the free equation have been studied using several different approaches. Jao \cite{Jao20} established refined mass-critical Strichartz estimates for Schr\"odinger operators through phase-space methods and bilinear restriction estimates. 
Camps \cite{Cam22} obtained bilinear Strichartz estimates for short-range potentials using distorted Fourier analysis and semiclassical tools. Pusateri and Soffer \cite{PS24} developed estimates for bilinear operators in distorted Fourier variables in the presence of large decaying potentials.
More recently, Schippa and Tataru \cite{ST25} proved sharp bilinear estimates for rough Hamiltonian flows using wave-packet decompositions, while Deng, Song, and Zhang \cite{DSZ26} established bilinear Strichartz estimates for the Schr\"odinger operator with a repulsive inverse-square potential in $\R^3$ through an interaction Morawetz argument.

In the present work, we take a different approach, based on the Duhamel formula and weighted smoothing estimates for the perturbed Schr\"odinger evolution. The starting point is the representation
\begin{align}
e^{-itH_V}f = e^{it\Delta}f -i\int_0^t e^{i(t-s)\Delta}\parenthese{V e^{-isH_V}f}\,ds.
\end{align}
Applying this formula to both factors in a bilinear expression reduces the problem to estimating free--free, free--Duhamel, and Duhamel--Duhamel interactions.
The free--free interaction is controlled by the corresponding bilinear Strichartz estimate for the free Schr\"odinger equation. The main difficulty is to control the terms involving the Duhamel integrals while preserving the gain from frequency separation.

A key ingredient is a family of weighted smoothing estimates associated with $H_V$. For Schr\"odinger operators with potentials of critical decay, Burq, Planchon, Stalker, and Tahvildar-Zadeh \cite{burq2004strichartz} established weighted spacetime estimates using resolvent bounds for the perturbed operator. Related methods have also been applied to dispersive estimates of Schr\"odinger equations on noncompact manifolds by Banica and Duyckaerts \cite{BD07}.
We extend these weighted estimates to a family of sharp Morawetz estimates, in the sense of Ozawa and Rogers \cite{ozawa2013sharp}, for the class of potentials considered here. These estimates provide the smoothing bounds needed to control the Duhamel terms in our bilinear argument. The resolvent framework we use is formulated in dimensions $d\geq 3$, which is the dimension range considered throughout this paper.

The smoothing estimates also require us to work with Sobolev spaces associated with the perturbed operator $H_V$. To relate these spaces to the usual Sobolev spaces, we use quadratic-form comparison, together with the inverse-square Sobolev estimates in \cite{killip2017energy,killip2018sobolev}. We also extend the local smoothing estimate in \cite{killip2017energy} from inverse-square potentials to the more general class of $(\kappa,\delta)$ potentials considered here.

Another important ingredient is an endpoint version of the bilinear Christ--Kiselev lemma. Away from the endpoint, the time truncations in the Duhamel formula can be handled by standard arguments. At the endpoint time exponent $p=2$, however, additional care is needed. We establish endpoint bilinear Christ--Kiselev estimates for both the free--Duhamel and Duhamel--Duhamel interactions, preserving the frequency-separation gain up to an arbitrarily small power loss in the frequency ratio. This allows us to treat these contributions in the bilinear estimates.

Let us emphasize that, although our approach uses a perturbative representation of the Schr\"odinger flow through Duhamel's formula, we do not require the potential $V$ to be small. Instead, our estimates depend on the structural assumptions imposed on $V$ and the corresponding weighted resolvent and smoothing bounds.

Combining these ingredients, we obtain two classes of global-in-time bilinear Strichartz estimates, depending on whether frequency localization is imposed before or after the perturbed evolution. We also control interactions between different input and output frequency scales, yielding a bilinear commutator-type estimate.

\subsection{Main results}\label{ssec Main}

We first establish a family of weighted spacetime estimates for the perturbed Schr\"odinger propagator, which generalize the linear Morawetz estimates for the free Schr\"odinger equation. These estimates will play an important role in the proof of our bilinear estimates.

\begin{restatable}{theorem}{genmorawetz}\label{thm genmorawetz}
Suppose $d \geq 3$, $V$ is a $(\kappa, \delta)$ potential and  $\sigma_1:= \max(\frac{d-2}{2}-\sqrt{\delta}, 0)$.
 Let $u = e^{-itH_V}u_0$ and $-1/2<s<d/2-1-\sigma_1$. Then
\begin{align*}
    \int_{\mathbb{R}} \int_{\mathbb{R}^d}|u(x, t)|^2 \frac{d x d t}{|x|^{2+2 s}} \leq C(d, s, \kappa, \delta)\norm{H_V^{\frac{s}{2}}u_0}_{L^2}^2 .
\end{align*}
\end{restatable}

The main effort of the paper is devoted to establishing bilinear Strichartz estimates for the perturbed Schr\"odinger evolution. As discussed above, the failure of commutation between $P_N$ and $e^{-itH_V}$ leads us to distinguish two types of estimates.

\noindent\underline{Type 1: Bilinear frequency-interaction Strichartz estimates}

Our first class of estimates concerns frequency-localized components of solutions to the perturbed Schr\"odinger equation. We show that, under suitable assumptions on the potential, bilinear estimates for the free Schr\"odinger evolution can be transferred to the perturbed flow.

For time exponents $p>2$, we obtain the same frequency-interaction bound as in the free case, up to a multiplicative constant depending on the potential. At the endpoint $p=2$, a small additional loss appears, which can be further improved under additional decay assumptions on the potential.

\begin{restatable}{theorem}{PBilinear}\label{thm PBilinear}
Let $V$ be a $(\kappa, \delta)$ potential (see Definition~\ref{defn GoodP}) and let $0<\varepsilon < 1$. Suppose that for all $\frac{1}{p}+\frac{1}{q}\leq 1$ and $\varphi, \psi \in L^2 (\R^d)$
\begin{align}
\norm{P_N e^{it\Delta} \varphi \, P_L e^{it\Delta } \psi}_{L_{t}^pL_{x}^q (\R \times \R^d)}  \leq C_{p, q}(N,L) \norm{\varphi}_{L^2 (\R^d)} \norm{\psi}_{L^2 (\R^d)}. 
\end{align}
Let $u_0,v_0 \in L^2 (\R^d)$. If $p>2$, then there exists $c=c(p, \kappa, \delta)$ such that 
\begin{align}
\norm{P_N e^{-itH_V} u_0 \, P_L e^{-itH_V} v_0}_{L_{t}^pL_{x}^q (\R \times \R^d)}  \leq cC_{p, q}(N,L) \norm{u_0}_{L^2 (\R^d)} \norm{v_0}_{L^2 (\R^d)} ,
\end{align}
and if $p=2$, then there exists $c=c(\kappa, \delta, \varepsilon)$ such that 
\begin{align}
    \norm{P_N e^{-itH_V} u_0 \, P_L e^{-itH_V} v_0}_{L_{t}^2L_{x}^q (\R \times \R^d)}  \leq  \left[ C_{2, q}(N, L)+ cC_{p_1,q_1}(N,L)^{1-\varepsilon}\right]  \norm{u_0}_{L^2 (\R^d)} \norm{v_0}_{L^2 (\R^d)} ,
\end{align}
where $p_1,q_1$ are given by
\begin{align*}
    \frac{1}{2}=\frac{d\varepsilon}{d+2}+\frac{1-\varepsilon}{p_1},
    \qquad \frac{1}{q}=\frac{d\varepsilon}{d+2}+\frac{1-\varepsilon}{q_1}.
\end{align*}
\end{restatable}

\begin{restatable}{theorem}{PPBilinear}\label{thm PBilinear1}
Let $V$ be a $(\kappa, \delta)$ potential with $|V(x)| \lesssim \inner{x}^{-2}$ and let $0<\varepsilon < 1/2$.  
There exists $c=c(\kappa, \delta, \varepsilon)$ such that
\begin{align}
& \quad \norm{P_N e^{-itH_V} u_0 \, P_L e^{-itH_V} v_0}_{L_{t,x}^2 (\R \times \R^d)} \\
& \leq  \left[ C(N, L)+ cC(N,L)^{1-2\varepsilon}\min\{N,L\}^{(d-2)\varepsilon} \parenthese{ \frac{1}{\inner{L}^{1/2}}+\frac{1}{\inner{N}^{1/2}}}\right]  \norm{u_0}_{L^2 (\R^d)} \norm{v_0}_{L^2 (\R^d)}
\end{align}
with $u_0, v_0 \in L^2 (\R^d)$.
\end{restatable}

\noindent\underline{Type 2: Bilinear Strichartz estimates with localized data}

Our second class of estimates concerns solutions whose initial data are localized at different frequencies. Since the perturbed propagator does not preserve the usual Fourier localization, the corresponding bilinear estimates require additional control of the frequency interactions generated by the potential.

We treat the low- and high-frequency regimes separately. For low frequencies, we establish bilinear estimates for $(\kappa,\delta)$ potentials satisfying additional spatial decay assumptions.
\begin{restatable}{theorem}{bilocdatahighmodes}\label{thm bilocdata high modes}

        Let $d\geq 3$ and $0< \varepsilon <1/2$. Let $\eta \geq 2$ and $V=\inner{x}^{-\eta}$. Let $f, g \in L^2\left(\mathbb{R}^d\right)$. 
   Then, there exists $C=C(d, \varepsilon)$ such that for $1\leq L\leq N$, the following inequality holds:
        \begin{align*}
            \left\|e^{-itH_V} P_N f \, e^{-itH_V} P_L g\right\|_{L_{t,x}^2 \left(\mathbb{R} \times \mathbb{R}^d\right)} \leq C \parenthese{\frac{L^{\frac{d-1}{2}}}{N^{\frac{1}{2}}}+\frac{L^{\frac{d-1}{2}-1-\varepsilon}}{N^{\frac{1}{2}-\varepsilon}}}\left\|P_Nf\right\|_2\left\|P_Lg\right\|_2 .
        \end{align*}
\end{restatable}

\begin{restatable}{theorem}{bilocdatalowmodes}\label{thm bilocdata lowmodes}
    

 Let $d\geq 3$, $\eta \geq 2$, and $V=\inner{x}^{-\eta}$.
   Assume $0 \leq \varepsilon < \min(2, d/2-1)$ and $\varepsilon \leq \eta-2$.  Let $f, g \in L^2\left(\mathbb{R}^d\right)$. There exists $C=C(d, \varepsilon)$ such that for $ L\leq N \leq 1$, the following inequality holds:
        \begin{align*}
            \left\|e^{-itH_V} P_N f \, e^{-itH_V} P_L g\right\|_{L_{t,x}^2 \left(\mathbb{R} \times \mathbb{R}^d\right)} \leq C\left(\frac{L^{\frac{d-1}{2}}}{N^{\frac{1}{2}}}+ L^{\frac{d-2}{2}}N^{\varepsilon}+N^{\frac{d-2}{2}}L^{\varepsilon} + L^{\varepsilon}+N^{\varepsilon}\right) \norm{P_Nf}_2\norm{P_Lg}_2 .
        \end{align*}
\end{restatable}

\noindent\underline{Commutator Type}

Combining the two types of bilinear estimates, we obtain a further result controlling interactions between solutions with frequency localization imposed both before and after the perturbed evolution. Such estimates quantify the effects of frequency mixing under the perturbed flow and provide a form of bilinear commutator control.
\begin{restatable}{theorem}{commbi}\label{thm commbi}
        Let $\varepsilon>0$, $\eta\geq 2$ and $V=\inner{x}^{-\eta}$. Let $f , g\in L^2(\mathbb{R}^d)$,  and let $N, L, M, K$ be dyadic numbers. If $1\leq L \leq N$, then
        \begin{align*}
            &\norm{ P_N e^{-itH_V}(P_Mf) \, P_Le^{-itH_V}(P_Kg)}_{L_{t,x}^2\left(\mathbb{R} \times \mathbb{R}^d\right)} \\
            &\lesssim_{\varepsilon, \eta} \parenthese{  \frac{L^{\frac{d-1}{2}}} {N^{\frac{1}{2}}}\norm{P_NP_M}_{L^2 \to L^2}\norm{P_LP_K}_{L^2 \to L^2}}\norm{P_Mf}_{L^2 (\R^d)} \norm{P_Kg}_{L^2 (\R^d)}\\
            &\quad+ \frac{L^{\frac{d-1}{2}-\varepsilon}}{N^{\frac{1}{2}-\varepsilon}}\parenthese{\frac{\norm{P_NP_M}_{L^2 \to L^2}}{\inner{L}^{\frac{1}{2}}\inner{K}^{\frac{1}{2}}} +\frac{\norm{P_LP_K}_{L^2 \to L^2}}{\inner{N}^{\frac{1}{2}}\inner{M}^{\frac{1}{2}}} +\frac{1}{\inner{L}^{\frac{1}{2}}\inner{K}^{\frac{1}{2}}\inner{N}^{\frac{1}{2}}\inner{M}^{\frac{1}{2}}}}  \norm{P_Mf}_{L^2 (\R^d)} \norm{P_Kg}_{L^2 (\R^d)}.
        \end{align*}
    \end{restatable}

In the course of proving these bilinear estimates, we also establish an endpoint version of the bilinear Christ–Kiselev lemma.
These results control the time-truncated Duhamel interactions at the endpoint time exponent $p=2$, for both the linear--Duhamel and Duhamel--Duhamel terms.
\begin{restatable}{theorem}{CK}\label{thm CK}
Let $d \geq 3$. Let $N,L$ be dyadic and $0 < \varepsilon \ll 1$.
Suppose that $u_0\in L^2(\R^d)$ and $F,G\in L_t^2L_x^{\frac{2d}{d+2}}(\R\times\R^d)$.
The endpoint versions of bilinear Christ-Kiselev Lemma hold for both 
\begin{enumerate}[(a)]
\item 
Linear-Duhamel interaction
\begin{align}
& \quad \norm{e^{it\Delta} P_N u_0 \parenthese{ \int_{-\infty}^t e^{i(t-s)\Delta} P_L F(s) \, ds }}_{L_{t,x}^2 (\R \times \R^d)} \\
& \leq C_{2,2}(N,L)^{1-\varepsilon} \min\{N,L\}^{\frac{d-2}{2} \varepsilon} \norm{u_0}_{L^2 (\R^d)} \norm{F}_{L_t^{2} L_x^{\frac{2d}{d+2}} (\R \times \R^d)} ;
\end{align}

\item 
Duhamel-Duhamel interaction
\begin{align}
& \quad \norm{\parenthese{ \int_{-\infty}^t e^{i(t-s)\Delta} P_N F (s) \, ds } \parenthese{ \int_{-\infty}^t e^{i(t-s')\Delta} P_L G(s') \, ds' }}_{L_{t,x}^2 (\R \times \R^d)} \\
& \leq C_{2,2}(N,L)^{1-\varepsilon} \min\{N,L\}^{\frac{d-2}{2} \varepsilon} \norm{F}_{L_t^{2} L_x^{\frac{2d}{d+2}}(\R \times \R^d)} \norm{G}_{L_t^{2} L_x^{\frac{2d}{d+2}}(\R \times \R^d)} .
\end{align}
\end{enumerate}
Moreover, both estimates remain valid with the $L_t^2L_x^{\frac{2d}{d+2}}$ norms on the right-hand sides replaced by $L_t^2L_x^{\frac{2d}{d+2},2}$ norms, provided the latter are finite.
\end{restatable}
We remark that the upper bounds in Lemma~\ref{lem CK} are in fact
\begin{align}
C_{2,2}(N,L)^{1-\varepsilon} \min \{N,L\}^{\frac{d-2}{2} \varepsilon} = \min\{N,L\}^{\frac{d-2}{2}} \parenthese{\frac{\min \{N,L\}}{\max \{ N,L\}} }^{\frac{1 - \varepsilon}{2}} ,
\end{align}
where we used the well-known bilinear Strichartz bound
\begin{align}
C_{2,2}(N,L) = \min\{N,L\}^{\frac{d-2}{2}} \parenthese{\frac{\min \{N,L\}}{\max \{ N,L\}} }^{\frac{1}{2}} .
\end{align}
Hence the endpoint Christ--Kiselev lemma only loses an arbitrarily small power of the frequency-separation ratio, $\frac{\min}{\max}$.

\subsection{Organization of the Rest of This Paper}
In Section~\ref{sec Prem}, we briefly review known bilinear estimates on $\R^d$ and collect the harmonic analysis tools, resolvent estimates, and smoothing estimates needed in our analysis.
In Section~\ref{sec sharpmor}, we establish sharp linear Morawetz estimates (Theorem~\ref{thm genmorawetz}).
In Section~\ref{sec CK}, we establish endpoint versions of the bilinear Christ-Kiselev lemma (Theorem~\ref{thm CK}).
In Section~\ref{sec PBilinear}, we prove bilinear frequency-interaction estimates for the perturbed evolution (Theorem~\ref{thm PBilinear} and Theorem~\ref{thm PBilinear1}).
In Sections~\ref{sec BilinearPHigh} and \ref{sec BilinearPLow}, we establish bilinear Strichartz estimates for frequency-localized initial data at high and low frequencies (Theorem~\ref{thm bilocdata high modes} and Theorem~\ref{thm bilocdata lowmodes}), respectively.
Finally, in Section~\ref{sec Commbi}, we combine these estimates to establish the bilinear commutator-type result (Theorem~\ref{thm commbi}).

\subsection{AI Statement} 
ChatGPT Pro was used to assist with proofreading for grammar and mathematical correctness. It also provided assistance in deriving Corollary~\ref{cor BDsmooth}. No AI assistance was used in the rest of the work.

\subsection*{Acknowledgement}
B. W. was supported by NSF grant DMS 1856124, and NSF CAREER Fellowship, DMS 2142064. 
X. Y. was partially supported by NSF DMS-2306429. Part of this work was completed during X.Y.'s participation in the SLMath program ``Recent Trends in Stochastic Partial Differential Equations".


\section{Preliminaries}\label{sec Prem}
In this section, we briefly review known bilinear estimates on $\R^d$ and collect the harmonic analysis tools, resolvent estimates, and smoothing estimates needed in our analysis.


\subsection{Notations}
We use the usual notation that $A \lesssim  B$ or $B \gtrsim A$ to denote an estimate of the form $A \leq C B$, for some constant $0 < C < \infty$ depending only on the {\it a priori} fixed constants of the problem. We also use $a+$ and $a-$ to denote expressions of the form $a + \varepsilon$ and $a - \varepsilon$, for any $0 <  \varepsilon \ll 1$.

Let $\phi:[0, \infty) \rightarrow[0,1]$ be a smooth function such that
\begin{align*}
    \phi(\lambda)=1 \text { for } 0 \leq \lambda \leq 1 \text { and } \phi(\lambda)=0 \text { for } \lambda \geq 2 \text {. }
\end{align*}
For each dyadic number, $N$, we define
\begin{align*}
    \phi_N(\lambda):=\phi(\lambda / N) \text { and } \psi_N(\lambda):=\phi_N(\lambda)-\phi_{N / 2}(\lambda) \text {. }
\end{align*}
Furthermore, define the standard Littlewood-Paley projections as follows:
    \begin{align*}
        P_{\leq N}:=\phi_N\left(\sqrt{-\Delta}\right), \quad P_N:=\psi_N\left(\sqrt{-\Delta}\right), \quad \text { and } P_{>N}:=I-P_{\leq N} .
    \end{align*}


\subsection{A Review on Bilinear Estimates on $\R^d$}

Standard Strichartz estimates appear in what are referred to as linear Strichartz estimates:

\begin{proposition}[Strichartz Estimates \cite{Str77,GV85,Yaj87,KT98}]\label{prop Strichartz}
Let $d\geq 1$. Let $q,r \in [2, \infty]$ satisfy 
    \begin{align*}
        \frac{2}{q}+ \frac{d}{r}= \frac{d}{2} ,
    \end{align*}
    and $(q, r, d) \neq (2, \infty, 2)$. Then there exists $C= C(q, r, d)$ such that 
    \begin{align*}
        \norm{e^{it \Delta}f}_{L^q_tL^r_x(\R \times \R^d)} \leq C \norm{f}_{L^2_x(\R^d)}.
    \end{align*}
\end{proposition}

The following bilinear interaction estimate was proven in \cite{Bour98} to handle self-interaction of waves for the critically-scaling nonlinear Schr\"odinger equation.

\begin{lemma}[Lemma 111 in \cite{Bour98}]\label{lem Bourbi}
Let $f,g \in L^2\left(\mathbb{R}^2\right)$. 
Then, for $L \leq N$, the following inequality holds:
\begin{align}\label{eq bourbi}
    \left\|\left(e^{i t \Delta} P_Nf\right)\left(e^{i t \Delta} P_Lg\right)\right\|_{L_{t,x}^2\left(\mathbb{R} \times \mathbb{R}^2 \right)} \leq C\left(\frac{L}{N}\right)^{1 / 2}\left\|P_Nf\right\|_{L^2} \left\|P_Lg\right\|_{L^2} .
\end{align}
\end{lemma}

One can generalize the argument of Bourgain to all dimensions as seen in \cite{CDKS01} with the main idea initially introduced in \cite{KPV96JAMS, KPV96TAMS}.

\begin{lemma}\label{lem alldbourbi}
Let $f,g \in L^2\left(\mathbb{R}^d\right)$. 
Then, for $L \leq N$, the following inequality holds:
\begin{align*}
    \left\|\left(e^{i t \Delta} P_Nf\right)\left(e^{i t \Delta} P_Lg\right)\right\|_{L_{t,x}^2\left(\mathbb{R} \times \mathbb{R}^d \right)} \leq C\frac{L^{\frac{d-1}{2}}}{N^{\frac{1}{2}}}\left\|P_Nf\right\|_{L^2} \left\|P_Lg\right\|_{L^2} .
\end{align*}
\end{lemma}

Interpolating Lemma~\ref{lem alldbourbi} and the following  bilinear estimates that follow from H\"older and linear Strichartz estimates
\begin{align}
\norm{(e^{it\Delta} P_Nf) (e^{it\Delta} P_Lg) }_{L_{t}^{\infty} L_x^1 (\R \times \R^d)} & \leq C \norm{P_Nf}_{L^2} \norm{P_Lg}_{L^2} \label{eq LinfL1} , \\
\norm{(e^{it\Delta} P_Nf) (e^{it\Delta} P_Lg) }_{L_{t}^{\infty} L_x^{\infty} (\R \times \R^d)} & \leq  (NL)^{\frac{d}{2}}C\norm{P_Nf}_{L^2} \norm{P_Lg}_{L^2} , \\
\norm{(e^{it\Delta} P_Nf) (e^{it\Delta} P_Lg) }_{L_{t}^{1} L_x^{\frac{d}{d-2}} (\R \times \R^d)} & \leq C\norm{P_Nf}_{L^2} \norm{P_Lg}_{L^2} ,
\end{align}
we have 
\begin{lemma}\label{lem alldbourbiallp}
Let $d \geq 3$. Let $f,g \in L^2\left(\mathbb{R}^d\right)$. Let $(p, q)$ be interpolatable from $(\infty, 1)$, $(\infty, \infty)$, $(2, 2)$ and $(1, \frac{d}{d-2})$. Then, for $L \leq N$, there exists $C(N,L)= C_{p, q, d}(N, L)$ such that the following inequality holds:
\begin{align*}
    \left\|\left(e^{i t \Delta} P_Nf\right)\left(e^{i t \Delta} P_Lg\right)\right\|_{L_{t}^pL^q_x\left(\mathbb{R} \times \mathbb{R}^d \right)} \leq C (N, L)\left\|P_Nf\right\|_{L^2} \left\|P_Lg\right\|_{L^2}.
\end{align*}
\end{lemma}

\subsection{Functional and Harmonic Analysis Tools}

We will need the following standard results throughout the proofs of the article. They will also serve as points of reference.

\begin{lemma}[Schur's Test, Lemma 9.2 in \cite{muscalu2013classical}]\label{lem:schur}
Define an integral operator on a measure space $X \times Y$ with the positive product measure $\mu \otimes \nu$ via
    \begin{align*}
        (Tf)(x) = \int_Y K(x, y) f(y) \, \nu(dy),
    \end{align*}
where $K$ is a measurable kernel. One has the following bounds:
    \begin{enumerate}
    \setlength{\itemsep}{0.2cm}
        \item 
        $\|T\|_{1 \to 1} \leq \sup_{y \in Y} \int_X |K(x, y)|\,\mu(dx) = : A$;
        \item 
        $\|T\|_{\infty \to \infty} \leq \sup_{x \in X} \int_Y |K(x, y)|\,\nu(dy) = : B$;
        \item $\|T\|_{p \to p} \leq A^{1/p}B^{1/p'}$, where $1 \leq p \leq \infty$;
        \item $\|T\|_{1 \to \infty} \leq \|K\|_{L^{\infty}(X \times Y)}$.
    \end{enumerate}
\end{lemma}

The $p< q$ condition in the following Christ-Kiselev Lemma will significantly limit the applicability of the Duhamel formula argument that we will employ throughout the rest of the paper.

\begin{lemma}[Christ-Kiselev Lemma, Theorem 1.2 in \cite{christ2001maximal} and Lemma 2.4 in \cite{tao2006nonlinear}]\label{lem:christkiselev}
Let $X, Y$ be Banach spaces, let $I$ be a time interval, and let $K \in C^0(I \times I \to B(X \to Y))$ be a kernel taking values in the space of bounded operators from $X$ to $Y$. Suppose that $1 \leq p < q\leq \infty$ is such that 
    \begin{align*}
        \norm{\int_I K(t, s) f(s)\, ds}_{L^q_t(I \to Y)} \leq A \|f\|_{L^p_t(I \to X)} ,
    \end{align*}
for all $f \in L^p_t(I \to X)$ and some $A>0$. Then one also has
    \begin{align*}
        \norm{\int_{s \in I: s<t} K(t, s) f(s)\, ds}_{L^q_t(I \to Y)} \lesssim_{p, q} A \|f\|_{L^p_t(I \to X)}.
    \end{align*}
\end{lemma}

The frequency/space localization proof of the following theorem that appears in \cite{muscalu2013classical} inspires the proof of Theorem~\ref{thm genmorawetz}.
\begin{proposition}[Hardy's Inequality, Theorem 9.5 in Muscalu-Schlag \cite{muscalu2013classical}]
    For any $0\leq s < d/2$ there is a constant, $C(s, d)$, with the property that
        \begin{align*}
            \norm{|x|^{-s} f}_{L^2_x} \leq C(s, d) \norm{f}_{\dot{H}^s_x} ,
        \end{align*}
    for all $f \in \dot{H}^s(\mathbb{R}^d)$.
\end{proposition}

Given Hardy's inequality, we can establish the following equivalence using interpolation and duality.

\begin{lemma}\label{lem normequiv}
     Let $d\geq 3$,  $f \in C^{\infty}_c(\mathbb{R}^d)$ and let $V$ be a $(\kappa, \delta)$ potential (See Definition~\ref{defn GoodP}). Then for any $s \in [-2, 2]$, there exists $C=C(s, d, \delta, \kappa)$ such that 
        \begin{align*}
            C^{-1}\|(-\Delta+V)^{s/2}f\|_{L^2_x} \leq \|(-\Delta)^{s/2}f\|_{L^2_x} \leq C \|(-\Delta+V)^{s/2}f\|_{L^2_x}
        \end{align*}
\end{lemma}




\subsection{Resolvent Estimates}

For any potential, $V$, and $\mu \in \mathbb{C}\setminus \mathbb{R}$, the resolvent is the operator defined by 
    \begin{align*}
        R_V(\mu):= (-\Delta +V - \mu)^{-1}.
    \end{align*}
As $\mu$ approaches a real number, $\xi$, the resolvent mimics a generalized Fourier projection at frequency $\xi$.  We would like to use the work of Burq, Planchon, Stalker, Tahvildar-Zadeh \cite{burq2004strichartz} on Morawetz and Strichartz estimates following from resolvent estimates. The authors establish general conditions for which their results hold. For the flat case, we have the classical estimate: for all $\lambda\geq 1$, $|\varepsilon|>0$, and $\eta>\frac{1}{2}$, there exists $C_{\eta}>0$ such that 
    \begin{align*}
        \norm{\inner{x}^{-\eta}(-\Delta-\lambda -i\varepsilon )^{-1}\inner{x}^{-\eta}}_{L^2 \to L^2} \leq C_{\eta} \lambda^{-1/2} ,
    \end{align*}
which leads to the following proposition via standard arguments.
\begin{proposition}\label{prop classicalresolvent}
For every $s \geq 0$ and $\eta> \frac{1}{2}$,  there exists $C=C_{s, \eta}>0$ such that
    \begin{align*}
        \sup _{\mu \in \mathbb{C} \backslash \mathbb{R}}\left\|\inner{x}^{-\eta}(-\Delta-\mu)^{-1}\inner{x}^{-\eta}\right\|_{H^{s-1} \rightarrow H^s} \leq C .
    \end{align*}
\end{proposition}

In order to pass to dispersive estimates from resolvent estimates, one can use a time-variable Fourier transform. Consider $G \in L^2_tH^s_x$ for $s\geq 0$. Then, for $\eta>\frac{1}{2}$, let
    \begin{align*}
        u(t):= e^{it\Delta}\int_0^t e^{-is\Delta}\inner{x}^{-\eta} G(s)\, ds ,
    \end{align*}
and observe that 
    \begin{align*}
        (i\partial_t +\Delta)u(t)= i\inner{x}^{-\eta}G(t).
    \end{align*}
Performing a Fourier transform on both sides with respect to $\mathcal{F}_t$, the Fourier transform in time, we get
    \begin{align*}
        (-\xi +\Delta)\mathcal{F}_t(u)(\xi)= i\inner{x}^{-\eta}\mathcal{F}_t(G)(\xi).
    \end{align*}
This observation allows for the following inhomogeneous estimate:
\begin{corollary}\label{cor duhsmooth}
    Let $s\geq 0$ and $\eta>1/2$. If $G \in L^2_tH^{s-1}_x$, then
        \begin{align*}
            \norm{\inner{x}^{-\eta} e^{it\Delta} \int_0^t e^{-is\Delta} \inner{x}^{-\eta}G(s)\,ds   }_{L^2_tH^s_x} \lesssim \norm{G}_{L^2_tH^{s-1}_x} .
        \end{align*}
\end{corollary}

\begin{proof}[Proof of Corollary~\ref{cor duhsmooth}]
    Let $G \in C^{\infty}_c(\mathbb{R} \times \mathbb{R}^d)$. By the previous computation with contour integration (similar to argument for Proposition 2.7 in \cite{burq2004nonlinear}),
        \begin{align*}
            \norm{\inner{x}^{-\eta} e^{it\Delta} \int_0^t e^{-is\Delta} \inner{x}^{-\eta}G(s)\,ds   }_{L^2_tH^s_x}\leq \sup_{\varepsilon>0}\norm{\inner{x}^{-\eta}(-\xi+\varepsilon i +\Delta)\inner{x}^{-\eta}\mathcal{F}_t(G)(\xi) }_{L^2_{\xi}H^s_x} .
        \end{align*}
    By Proposition~\ref{prop classicalresolvent}, we have
       \begin{align*}
            \norm{\inner{x}^{-\eta} e^{it\Delta} \int_0^t e^{-is\Delta} \inner{x}^{-\eta}G(s)\,ds   }_{L^2_tH^s_x}\lesssim \norm{\mathcal{F}_t(G)(\xi) }_{L^2_{\xi}H^{s-1}_x} =\norm{G}_{L^2_tH^{s-1}_x}.
        \end{align*} 
\end{proof}

The $TT^*$ argument of Proposition 2.7 in \cite{burq2004nonlinear} further provides the standard smoothing estimate:
    \begin{corollary}\label{cor locsmoothflat}
        Let $s\geq 0$ and $\eta\geq 1$. If $f \in H^{s}_x$, then
        \begin{align*}
            \norm{\inner{x}^{-\eta} e^{it\Delta} f  }_{L^2_tH^{s+1/2}_x} \lesssim \norm{f}_{H^{s}_x} .
        \end{align*}
    \end{corollary}

The commutator $[\Delta, \inner{x}^{-\eta}]$ plays a crucial role in the arguments above. It is easier to prove the statements for a general class of weights as long as the weights are compactly supported. Nonlocal weights may present difficulties even when they are short range.

The conditions outlined in \cite{burq2004strichartz} provide resolvent estimates for a broad class of potentials are simultaneously critically scaling at infinity and critically singular locally.

\begin{definition}\label{defn GoodP}
Let $d \geq 3$ and $\kappa, \delta> 0$. Then $W \in C^1(\mathbb{R}^d\setminus \{0\})$ is a {\bf $(\kappa, \delta)$ Potential}, if 
        \begin{align*}
        &\quad |W(x)| \leq \frac{\kappa}{|x|^2} , \\ 
          &\left(\frac{d}{2}-1\right)^2+|x|^2 W(x) \geq \delta\quad \forall x \in \mathbb{R}^d \backslash\{0\},\\
        &\left(\frac{d}{2}-1\right)^2-|x|^2 \partial_{|x|}(|x| W(x)) \geq \delta .
        \end{align*}
\end{definition}

\begin{proposition}[Burq, Planchon, Stalker, and Tahvildar-Zadeh \cite{burq2004strichartz}, Special case of Theorem 3]\label{prop resolv}
     Let $W \in C^1\left(\mathbb{R}^d \setminus \{0\}\right)$ be a $(\kappa, \delta)$ Potential.
Then there exists $C_0>0$ such that
    \begin{align*}
        \sup _{\mu \in \mathbb{C} \backslash \mathbb{R}}\left\||x|^{-1}R_W(\mu)|x|^{-1}\right\|_{L^2 \rightarrow L^2} \leq C_0 .
    \end{align*}
\end{proposition}

The authors then apply a classical result of Kato (Theorem 4 in \cite{burq2004strichartz}) to derive the following Morawetz estimate:
\begin{corollary}[Section 3.1.1 in \cite{burq2004strichartz}]\label{cor BPSTZ resolvent}
        Let $V$ be a $(\kappa, \delta)$ potential let $u = e^{-itH_V}u_0$. Then 
    \begin{align}\label{eq linmor}
        \||x|^{-1}u\|_{L_{t,x}^2 (\R \times \R^d)}  \leq \frac{1}{\delta^2 \sqrt{2\pi}}\|u_0\|_{L^2_x (\R^d)} , 
    \end{align}
    and thus
        \begin{align*}
            \int_{\mathbb{R}} \int_{|x-z| \leq R} | u(t, x)|^2 \,dx \,dt \lesssim R^2\norm{ u_0}^2_{L_x^2 (\R^d)},
        \end{align*}
        uniformly for $z \in \mathbb{R}^d$  and $R>0$.
\end{corollary}


\subsection{Smoothing Estimates}

We also have the following smoothing estimate of Banica-Duyckaerts
\begin{proposition}[Banica, Duyckaerts \cite{BD07}, Special case of Theorem 2]\label{thm BDsmooth}
    Let  $V \in C^1\left(\mathbb{R}^d \setminus \{0\}\right)$ be a $(\kappa, \delta)$ potential, and let $H=-\Delta+V$. Then there exists $C_0>0$ such that
    \begin{align*}
        \norm{\langle x\rangle^{-1}e^{-itH} f }_{L^2_t H^{\frac{1}{2}}_x} \leq C_0 \norm{f}_{L^2_x} .
    \end{align*}
\end{proposition}

In the case that the potential has better regularity than the $(\kappa, \delta)$ condition, we can combine the preceding theorem with Corollary~\ref{cor duhsmooth} and Corollary~\ref{cor locsmoothflat} to get smoothing estimates across all regularities. For simplicity, we consider the case that $V(x) = \inner{x}^{-\eta}$ for $\eta\geq 2$. In which case $V$ is a $(\kappa, \delta)$ potential.

\begin{corollary}\label{cor BDsmooth}
    Let  $\eta\geq 2$, $s \in \mathbb{R}$, and let $H=-\Delta+\inner{x}^{-\eta}$. Then there exists $C_s>0$ such that
    \begin{align}\label{eq SE}
        \norm{\langle x\rangle^{-1}e^{-itH} f }_{L^2_t H^{s}_x} \leq C_s \norm{f}_{H^{s-\frac{1}{2}}_x} .
    \end{align}
\end{corollary}

\begin{proof}[Proof of Corollary~\ref{cor BDsmooth}]

We first show that \eqref{eq SE} holds true when $s\geq 1$, which serves as the base step in the induction argument. Note that 
    \begin{align*}
        \norm{\langle x\rangle^{-1}e^{-itH} f }_{L^2_t H^{1}_x} \leq \norm{\langle x\rangle^{-1}e^{it\Delta} f }_{L^2_t H^{1}_x} +\norm{\inner{x}^{-1} e^{it\Delta} \int_0^t e^{-is\Delta} \inner{x}^{-\eta}e^{-isH} f\,ds   }_{L^2_tH^{1}_x} . 
    \end{align*}
    By Corollary~\ref{cor locsmoothflat} and estimate \eqref{eq linmor}
    \begin{align*}
        &\norm{\langle x\rangle^{-1}e^{it\Delta} f }_{L^2_t H^{1}_x} +\norm{\inner{x}^{-1} e^{it\Delta} \int_0^t e^{-is\Delta} \inner{x}^{-\eta}e^{-isH} f\,ds   }_{L^2_tH^{1}_x}\\
        &\lesssim \norm{ f }_{ H^{1/2}_x}+\norm{ \inner{x}^{-\eta+1}e^{-itH} f  }_{L^2_tL^2_x}\\
        &\lesssim \norm{ f }_{ H^{1/2}_x}+\norm{ |x|^{-1}e^{-itH} f  }_{L^2_tL^2_x}\\
        &\lesssim \norm{ f }_{ H^{1/2}_x} .
    \end{align*} 

    Next, we prove that \eqref{eq SE} holds true when $s=m$ with $m\geq3$ being an odd integer. 
    
   For $m\geq3$ and odd, we have
    \begin{align*}
        &\norm{\langle x\rangle^{-1}e^{-itH} f }_{L^2_t H^{m}_x}\\
        &\leq\norm{\langle x\rangle^{-1}e^{it\Delta} f }_{L^2_t H^{m}_x} +\norm{\inner{x}^{-1} e^{it\Delta} \int_0^t e^{-is\Delta} \inner{x}^{-\eta}e^{-isH} f\,ds   }_{L^2_tH^{m}_x}\\
        &\lesssim\norm{ f }_{ H^{m-1/2}_x} +\norm{ \inner{x}^{-\eta+1}e^{-itH} f }_{L^2_tH^{m-1}_x}\\
        &\lesssim\norm{ f }_{ H^{m-1/2}_x} +\norm{\langle x\rangle^{-\eta+1}e^{it\Delta} f }_{L^2_t H^{m-1}_x} +\norm{\inner{x}^{-\eta+1} e^{it\Delta} \int_0^t e^{-is\Delta} \inner{x}^{-\eta}e^{-isH} f\,ds   }_{L^2_tH^{m-1}_x} .
    \end{align*} 
    Since $\eta-1\geq1$, Corollaries~\ref{cor duhsmooth} and \ref{cor locsmoothflat} imply
    \begin{align*}
        &\norm{\langle x\rangle^{-\eta+1}e^{it\Delta} f }_{L^2_t H^{m-1}_x} +\norm{\inner{x}^{-\eta+1} e^{it\Delta} \int_0^t e^{-is\Delta} \inner{x}^{-\eta}e^{-isH} f\,ds   }_{L^2_tH^{m-1}_x}\\
        &\lesssim \norm{f}_{H^{m-\frac{3}{2}}_x} +\norm{ \inner{x}^{-1}e^{-itH} f }_{L^2_tH^{m-2}_x} . 
    \end{align*}
    With the base step, $m=1$, induction provides the desired estimate for $m \in 2\mathbb{Z}_+-1$.

    We then consider the negative odd integer case in \eqref{eq SE}. Define
    \begin{align}
        \mathcal{P}= 1 - \Delta +V ,   
    \end{align}
    and conjugated operator
    \begin{align}
    \mathcal{P}_{\inner{x}} = \inner{x}^{-1} \mathcal{P} \inner{x} = \inner{x}^{-1} (1 - \Delta +V) \inner{x}.
    \end{align}
    For a smooth function $g$, direct differentiation gives
    \begin{align}
        \mathcal{P}_{\inner{x}} g= [\inner{x}^{-1} (1 - \Delta +V) \inner{x}] g & = (1-\Delta+V)g -2\inner{x}^{-1}\nabla(\inner{x})\cdot\nabla g -\inner{x}^{-1}\Delta(\inner{x})g\\
        & = (1-\Delta+V)g -\frac{2x}{\inner{x}^2}\cdot\nabla g -\frac{d}{\inner{x}^2}g +\frac{|x|^2}{\inner{x}^4}g \\
        & = : \mathcal{P} g + b(x) \cdot \nabla g  + c(x) g \label{eq16}
    \end{align}
    where $b,c$ are smooth functions with all derivatives bounded.
    Thus $\mathcal{P}_{\inner{x}}$ has smooth bounded coefficients and bounded derivatives. 

    Multiplication by a smooth function with bounded derivatives of all orders defines a bounded operator on $H^s_x$ for every $s\in\mathbb{R}$. Indeed, this follows from the product rule for nonnegative integer $s$, duality for negative integers, and interpolation for all remaining indices. Expanding the conjugated operator by the product rule and using \eqref{eq16}, we obtain, for every $s\in\mathbb{R}$,
\begin{align}\label{eq:Bdirect}
    \norm{\mathcal{P}_{\inner{x}} g}_{H^s_x}
    = \norm{\inner{x}^{-1}(1-\Delta+V)\inner{x}g}_{H^s_x}
    \lesssim_s \norm{g}_{H^{s+2}_x}.
\end{align}
    We now claim that
    \begin{align}\label{eq:Breverse}
        \norm{g}_{H^{s+2}_x}
        \lesssim_s\norm{\mathcal{P}_{\inner{x}} g}_{H^s_x}+\norm{g}_{H^s_x}.
    \end{align}
    Note that by Plancherel, $\|g\|_{H^{s+2}} \sim \|(\Delta-1) g\|_{H^s}$. Since $\mathcal{P}_{\inner{x}} -  (1-\Delta)$ is composed of first order and zeroth order terms, we can observe that
    \begin{align*}
        \norm{g}_{H^{s+2}_x} \sim \|(\Delta-1) g\|_{H^s}\leq\norm{\mathcal{P}_{\inner{x}} g}_{H^s_x}+C_s\norm{g}_{H^{s+1}_x}.
    \end{align*}
    Using
    \begin{align*}
        \norm{g}_{H^{s+1}_x}
        \leq\varepsilon\norm{g}_{H^{s+2}_x}
        +\frac{1}{\varepsilon}\norm{g}_{H^s_x} ,
    \end{align*}
    we write
    \begin{align}
        \norm{g}_{H^{s+2}_x} \leq\norm{\mathcal{P}_{\inner{x}} g}_{H^s_x}+C_s \varepsilon \norm{g}_{H^{s+2}_x} + \frac{C_s}{\varepsilon} \norm{g}_{H^s_x}. 
    \end{align}
    Choosing $\varepsilon$ so that $C_s\varepsilon<\frac{1}{2}$ and absorbing the $\norm{g}_{H^{s+2}_x}$ term into the left-hand side, we prove the claimed estimate \eqref{eq:Breverse}.

    We also need the elliptic inverse estimate
    \begin{align}\label{eq:ellipticinverse}
        \norm{\mathcal{P}^{-1}g}_{H^{s+2}_x}\lesssim_s\norm{g}_{H^s_x},\qquad s\in\mathbb{R}.
    \end{align}
    Since $V\geq0$ (hence the spectrum of the operator $\sigma(\mathcal{P}) \subset [1,\infty)$), the spectral theorem gives 
    \begin{align}
        \norm{\mathcal{P}^{-1}g}_{L^2_x}\leq\norm{g}_{L^2_x} .
    \end{align}
    Recall again that $\|g\|_{H^{s}} \sim \|(\Delta-1)^{-1} g\|_{H^{s+2}}$. The identity $(1-\Delta) \mathcal{P}^{-1}g=g-V\mathcal{P}^{-1}g$ implies $\mathcal{P}^{-1}g=(1-\Delta)^{-1}g-(1-\Delta)^{-1}V\mathcal{P}^{-1}g$. Therefore,
    \begin{align*}
        \norm{\mathcal{P}^{-1}g}_{H^{j+2}_x}\lesssim_j\norm{g}_{H^j_x}+\norm{V\mathcal{P}^{-1}g}_{H^j_x},\qquad j\geq0.
    \end{align*}
     Since $V= \inner{x}^{-\eta}$, we can remove the potential
     \begin{align*}
        \norm{\mathcal{P}^{-1}g}_{H^{j+2}_x}\lesssim_j\norm{g}_{H^j_x}+\norm{\mathcal{P}^{-1}g}_{H^j_x},\qquad j\geq0.
    \end{align*}
    Induction proves \eqref{eq:ellipticinverse} for nonnegative integer $s$. Self-adjointness and duality prove it for integer $s\leq -2$; interpolation supplies the remaining indices.

    Suppose that the desired smoothing estimate is known at
    some index $m \in \mathbb{Z}$:
    \begin{align}\label{eq:BDknownindex}
        \norm{\inner{x}^{-1} e^{-itH}f}_{L^2_tH^m_x}
        \lesssim_m\norm{f}_{H^{m-\frac{1}{2}}_x}.
    \end{align}

    We lower the index by two. By \eqref{eq:Bdirect},
    \eqref{eq:BDknownindex}, and
    \eqref{eq:ellipticinverse},
    \begin{align*}
        \norm{\inner{x}^{-1} e^{-itH}f}_{L^2_tH^{m-2}_x}
        &=
        \norm{\mathcal{P}_{\inner{x}}\parenthese{\inner{x}^{-1} e^{-itH}\mathcal{P}^{-1}f}}_{L^2_tH^{m-2}_x}\\
        &\lesssim_m
        \norm{\inner{x}^{-1} e^{-itH}\mathcal{P}^{-1}f}_{L^2_tH^m_x}\\
        &\lesssim_m\norm{\mathcal{P}^{-1}f}_{H^{m-\frac{1}{2}}_x}\\
        &\lesssim_m\norm{f}_{H^{m-\frac52}_x}.
    \end{align*}
    This is the desired estimate at $m-2$. With base step, $m=1$, induction provides the desired estimate for $2\mathbb{Z}_-+3$. Therefore, the smoothing estimate holds for $m \in 2\mathbb{Z}-1.$ Interpolation provides the estimate \eqref{eq SE} for $s \in \mathbb{R}$.
\end{proof}




    

We may also obtain smoothing estimates in a more direct fashion that is given in the CLAY Math notes of Killip and Visan:

\begin{proposition}[Proposition 4.14 in \cite{KVClay}]\label{prop KVlocsmoothing}
    Fix $\varphi \in C^{\infty}_c(\mathbb{R}^d)$. Then for all $f \in L^2(\mathbb{R}^d)$ and $R>0$,
        \begin{align}
            \int_{\mathbb{R}}\int_{\mathbb{R}^d} \left| (|\nabla|^{\frac{1}{2}} e^{it\Delta}f)(x)\right|^2 \varphi(x/R) \,dxdt \lesssim_{\varphi}R\|f\|^2_{L^2(\mathbb{R}^d)} ,
        \end{align}
        and so for $\varepsilon>0$,
        \begin{align}
            \int_{\mathbb{R}} \int_{\mathbb{R}^d}  \left| (|\nabla|^{\frac{1}{2}} e^{it\Delta}f)(x)\right|^2 \langle x\rangle^{-1-\varepsilon}\,dxdt\lesssim_{\varepsilon}\|f\|^2_{L^2(\mathbb{R}^d)} .
        \end{align}
\end{proposition}

This proposition is essentially an endpoint for the Morawetz estimate, Proposition~\ref{thm ozawarogers}, that we will discuss in Section~\ref{sec sharpmor}.

Here is an implication of the above estimate:

\begin{corollary}\label{cor:homsmooth}
    Let $F: \mathbb{R} \to L^2(\mathbb{R}^d)$. Let $\eta\geq 1$ and $s\geq 0$, then 
    \begin{align*}
        \norm{ \int_{\mathbb{R}}e^{-it\Delta}\langle x\rangle^{-\eta}F(t)\,dt }_{ H^{s+1/2}(\mathbb{R}^d)} \lesssim_{\eta, s} \|F\|_{L^2(\mathbb{R}; H^{s}(\mathbb{R}^d))}.
    \end{align*}
\end{corollary}

\begin{proof}[Proof of Corollary~\ref{cor:homsmooth}]
First, Corollary~\ref{cor BPSTZ resolvent} with $V=0$ and $\inner{x}^{-\eta}\leq|x|^{-1}$  gives, by duality,
\begin{align*}
    \norm{\int_{\mathbb{R}}e^{-it\Delta}\inner{x}^{-\eta}F(t)\,dt}_{L^2_x}\lesssim\norm{F}_{L_{t,x}^2}.
\end{align*}
This shows the low-frequency part of the inhomogeneous Sobolev norm.

We observe that for $m \in \mathbb{Z}_{\geq 0}$,
     \begin{align*}
        \norm{ \int_{\mathbb{R}}e^{-it\Delta}\langle x\rangle^{-\eta}F(t)\,dt }_{ \dot{H}^{m+1/2}(\mathbb{R}^d)} \lesssim \sum_{k_1+k_2=m}\norm{ \int_{\mathbb{R}}e^{-it\Delta}\nabla^{k_1}\langle x\rangle^{-\eta}\nabla^{k_2}F(t)\,dt }_{ \dot{H}^{1/2}(\mathbb{R}^d)}.
    \end{align*}
Now note that, by Proposition~\ref{prop KVlocsmoothing}, 
    \begin{align}
        \int_{\mathbb{R}} \int_{\mathbb{R}^d}  \left| (\nabla^{k_1}\inner{x}^{-\eta})(|\nabla|^{\frac{1}{2}} e^{it\Delta}f(x))\right|^2\,dxdt&\lesssim_{k_1} \int_{\mathbb{R}} \int_{\mathbb{R}^d}  \left| (\inner{x}^{-\eta}|\nabla|^{\frac{1}{2}} e^{it\Delta}f(x))\right|^2\,dxdt\\
        &\lesssim_{\eta, k_1} \norm{f}^2_{L^2(\mathbb{R}^d)} . 
    \end{align}
Duality implies 
    \begin{align*}
        \norm{ \int_{\mathbb{R}}e^{-it\Delta}\nabla^{k_1}\langle x\rangle^{-\eta}\nabla^{k_2}F(t)\,dt }_{ \dot{H}^{1/2}(\mathbb{R}^d)} \lesssim_{\eta, k_1}\norm{\nabla^{k_2}F}_{L^2(\mathbb{R}; L^2(\mathbb{R}^d))}.
    \end{align*}
Therefore,
    \begin{align*}
        \norm{ \int_{\mathbb{R}}e^{-it\Delta}\langle x\rangle^{-\eta}F(t)\,dt }_{ \dot{H}^{m+1/2}(\mathbb{R}^d)} \lesssim_{\eta, m} \norm{F}_{L^2(\mathbb{R}; H^m(\mathbb{R}^d))}.
    \end{align*}
    Combining this with the preceding $L^2_x$ estimate gives the $H^{m+1/2}_x$ bound. Interpolation completes the argument.
\end{proof}



    \subsection{The Smoothing Estimate of Killip-Miao-Visan-Zhang-Zheng \cite{killip2017energy}}

First, for $a \in \mathbb{R}$, define the operator
    \begin{align*}
        \mathcal{L}_a := -\Delta + \frac{a}{|x|^2} .
    \end{align*}

    The following lemma was proven:
    \begin{lemma}[Local Smoothing, Lemma 2.9 in \cite{killip2017energy}]\label{lem KMVZZ}
        Fix $a> - \frac{1}{4}$ and let $u = e^{-it\mathcal{L}_a}u_0$. Then 
        \begin{align*}
            \int_{\mathbb{R}} \int_{\mathbb{R}^3} \frac{|\nabla u(t, x)|^2}{R\left\langle R^{-1} x\right\rangle^3}+\frac{|u(t, x)|^2}{R|x|^2} \,dx \,dt \lesssim \norm{u_0}_{L^2_x}\norm{\nabla u_0}_{L^2_x}+R^{-1}\norm{u_0}^2_{L^2_x} , \\
            \int_{\mathbb{R}} \int_{|x-z| \leq R} \frac{1}{R}|\nabla u(t, x)|^2 \,dx \,dt \lesssim\norm{u_0}_{L_x^2}\norm{\nabla u_0}_{L_x^2}+R^{-1}\norm{u_0 }_{L_x^2}^2,
        \end{align*}
        uniformly for $z \in \mathbb{R}^3$  and $R>0$.
    \end{lemma}

The lemma holds in general dimensions, but is not presented in such a way in \cite{killip2017energy}. The generalization is a formality that we sequester to Appendix~\ref{app KMVZZsmoothing}.
    \begin{lemma}[Local Smoothing] \label{lem KMVZZ smoothing}
       Let $V$ be a $(\kappa, \delta)$ potential and let $u = e^{-it(-\Delta+ V)}u_0$. Then 
        \begin{align*}
            \int_{\mathbb{R}} \int_{\mathbb{R}^d} \frac{|\nabla u(t, x)|^2}{R\left\langle R^{-1} x\right\rangle^3}+\frac{|u(t, x)|^2}{R|x|^2} \,dx \,dt \lesssim \norm{u_0}_{L^2_x}\norm{\nabla u_0}_{L^2_x}+R^{-1}\norm{u_0}^2_{L^2_x} , \\
            \int_{\mathbb{R}} \int_{|x-z| \leq R} \frac{1}{R}|\nabla u(t, x)|^2 \,dx \,dt \lesssim\norm{u_0}_{L_x^2}\norm{\nabla u_0}_{L_x^2}+R^{-1}\norm{u_0 }_{L_x^2}^2,
        \end{align*}
        uniformly for $z \in \mathbb{R}^d$  and $R>0$.
    \end{lemma}

    From this lemma, we have the local smoothing estimate with sharp decay of the weight, resembling the second estimate of Proposition~\ref{prop KVlocsmoothing}:
\begin{corollary}\label{cor KMVZZ smoothing}
        Let $V$ be a $(\kappa, \delta)$ potential, $\varepsilon, R>0$, and let $u = e^{-it(-\Delta+ V)}u_0$. Then 
        \begin{align*}
            \int_{\mathbb{R}} \int_{\mathbb{R}^d} \frac{|\nabla u(t, x)|^2}{R\left\langle R^{-1} x\right\rangle^{1+\varepsilon}}\,dx \,dt \lesssim_{\varepsilon} \norm{u_0}_{L_x^2}\norm{\nabla u_0}_{L_x^2}+R^{-1}\norm{u_0 }_{L_x^2}^2.
        \end{align*}
\end{corollary}

\begin{proof}[Proof of Corollary~\ref{cor KMVZZ smoothing}]
    \begin{align*}
            \int_{\mathbb{R}} \int_{\mathbb{R}^d} \frac{|\nabla u(t, x)|^2}{R\left\langle R^{-1} x\right\rangle^{1+\varepsilon}}\,dx \,dt  &\lesssim \sum_{k=0}^{\infty}\int_{\mathbb{R}} \int_{|x|\leq 2^kR} \frac{|\nabla u(t, x)|^2}{R2^{(1+\varepsilon)k}}\,dx \,dt \\
            &\lesssim \sum_{k=0}^{\infty}2^{-\varepsilon k}\int_{\mathbb{R}} \int_{|x|\leq 2^kR} \frac{|\nabla u(t, x)|^2}{R2^{k}}\,dx \,dt\\
            &\lesssim \sum_{k=0}^{\infty}2^{-\varepsilon k}\left(\norm{u_0}_{L_x^2}\norm{\nabla u_0}_{L_x^2}+(R2^{k})^{-1}\norm{u_0 }_{L_x^2}^2\right) .
        \end{align*}
        The last inequality is due to the second estimate in Lemma~\ref{lem KMVZZ smoothing}.
\end{proof}

We note that it is clear from the proof of Lemma~\ref{lem KMVZZ smoothing} that for $V\equiv 0$, we have essentially the same estimate as Proposition~\ref{prop KVlocsmoothing}:
        \begin{align*}
            \int_{\mathbb{R}} \int_{\mathbb{R}^d} \frac{|\nabla e^{it\Delta}u_0|^2}{R\left\langle R^{-1} x\right\rangle^{1+\varepsilon}}\,dx \,dt \lesssim_{\varepsilon} \norm{u_0}_{L_x^2}\norm{\nabla u_0}_{L_x^2}.
        \end{align*}



    \subsection{Lorentz Spaces and Endpoint Strichartz Estimates}

        Recall the definition of Lorentz spaces: For  $ p \in (0, \infty)$ and $q \in (0, \infty)$. If $\mu$  is Lebesgue measure on $\mathbb{R}^d$
        \begin{align*}
            \norm{f}_{L^{p, q}(\mathbb{R}^d)} := p^{1/q} \left( \int_0^{\infty}  \parenthese{t \mu \{ |f|\geq t\}^{\frac{1}{p}}}^q \,\frac{dt}{t} \right)^{1/q} . 
        \end{align*}
        For $q=\infty$,
        \begin{align*}
            \norm{f}_{L^{p, \infty}(\mathbb{R}^d)}:= \sup_{t > 0} t \mu \{ |f|\geq t\}^{\frac{1}{p}} .
        \end{align*}
        In similar fashion as Banica-Duyckaerts \cite{BD07}, we want to describe the process by which one combines O'Neil's inequality, the endpoint Strichartz estimates of Keel and Tao \cite{KT98}, and local smoothing to obtain suitable  Strichartz estimates for Schr\"odinger with potential. First, O'Neil's inequality states that 
            \begin{align*}
                \|fg\|_{L^{\frac{2d}{d+2}, 2}(\mathbb{R}^d)} \lesssim \|f\|_{L^2 (\R^d)}\|g\|_{L^{d, \infty} (\R^d)} .
            \end{align*}
            The relevant Strichartz estimate of Keel and Tao is as follows:
        \begin{proposition}[Specific case of Theorem 1.2 in \cite{KT98}]\label{prop lorentzstrichartz}
        Let $d \geq 3$.
            There is a constant $C>0$ such that if $u:\R \times \R^d \to \mathbb{R}$ is a solution to 
            \begin{align*}
                \left\{\begin{array}{c} i\partial_t u + \Delta u = F \\ u(0)=u_0, \end{array}  \right.
            \end{align*}
            then 
            \begin{align*}
                \|u\|_{L^2( \mathbb{R}, L^{\frac{2d}{d-2}, 2})} \leq C \left( \|u_0\|_{L^2} + \|F\|_{L^2(\mathbb{R}, L^{\frac{2d}{d+2}, 2})} \right).
            \end{align*}
        \end{proposition}
    This exact statement does not appear in \cite{KT98}, but it follows from the last paragraph of Section 5 of the paper.

\subsection{Lions-Peetre Interpolation Theorem}\label{sec LP}
The definitions and results in this subsection are adapted from Sections 2.1-2.4 in the lecture notes on interpolation theory of Tice \cite{TiceInterpolation}.
\begin{definition}
Suppose that $X_0, X_1$ are Banach spaces over the field $\mathbb{F}$. We say $X_0$ and $X_1$ are compatible if there exists a Hausdorff topological vector space $Z$ over $\mathbb{F}$ such that $X_i \hookrightarrow Z$ is a subspace for $i=0,1$ and the inclusions are continuous.  
\end{definition}
For example, $L^p$ and $L^q$ are compatible, as they are both continuously embedded into the space of measurable functions.

\begin{definition}
Let $X_0$ and $X_1$ be compatible Banach spaces. We define the map $K : (X_0 + X_1) \times \R_+ \to \R$ by
\begin{align}
K(x,t) = \inf \{ \norm{x_0}_{X_0} + t \norm{x_1}_{X_1} : x = x_0 + x_1 \text{ for } x_i \in X_i, i = 0,1\} .
\end{align}

For $1 \leq p \leq \infty$ and $\theta \in (0,1)$, we define the space
\begin{align}
(X_0 ,X_1)_{\theta , p} = \{ x \in X_0 + X_1 : \norm{x}_{\theta ,p} < \infty\} , 
\end{align}
where $\norm{\cdot}_{\theta , p} : X_0 + X_1 \to [0, \infty]$ is defined by
\begin{align}
\norm{x}_{\theta ,p} = 
\begin{cases}
\displaystyle \parenthese{\int_{\R_+} \parenthese{\frac{K(x,t)}{t^{\theta}}}^p \frac{dt}{t}}^{\frac{1}{p}} & p < \infty \\
\displaystyle \sup_{t >0} \frac{K(x,t)}{t^{\theta}} & p = \infty
\end{cases}
\end{align}
\end{definition}

\begin{proposition}
Let $X_0$ and $X_1$ be compatible Banach spaces, $1 \leq p \leq \infty$, and $\theta \in (0,1)$. Then there exists a constant $C \in \R_+$ such that if $x \in X_0 \cap X_1$, then
\begin{align}
\norm{x}_{\theta ,p} \leq C \norm{x}_{X_0}^{1-\theta} \norm{x}_{X_1}^{\theta}
\end{align}
for all $x \in X_0 \cap X_1$.
\end{proposition}

\begin{proposition}
Let $X_0$ and $X_1$ be compatible Banach spaces and $Y_0$ and $Y_1$ be compatible Banach spaces over the same field. Suppose that $T: X_0 + X_1 \to Y_0 + Y_1$ is a linear map such that $T(X_i) \subset Y_i$ and $T \in \mathcal{L} (X_i, Y_i)$ for $i =0,1$. Then for every $1 \leq p \leq \infty$ and $\theta \in (0,1)$, we have that $T \in \mathcal{L} ( (X_0,X_1)_{\theta ,p} , (Y_0 , Y_1)_{\theta ,p})$ and
\begin{align}
\norm{T}_{\mathcal{L} ((X_0 , X_1)_{\theta ,p} , (Y_0, Y_1)_{\theta , p})} \leq \norm{T}_{\mathcal{L} (X_0 , Y_0)}^{1-\theta} \norm{T}_{\mathcal{L} (X_1 , Y_1)}^{\theta} .
\end{align}
\end{proposition}

\begin{proposition}

Let $1< p_0,p_1 < \infty$, $p_0 \neq p_1$, and $1 \leq q_0, q_1 \leq \infty$. Let $1 \leq q \leq \infty$, $\theta \in (0,1)$, and define $1 < p_{\theta} < \infty$ by
\begin{align}
\frac{1}{p_{\theta}} = \frac{1 - \theta}{p_0} + \frac{\theta}{p_1} .
\end{align}
Then
\begin{align}
(L^{p_0, q_0} (X) , L^{p_1 ,q_1} (X))_{\theta ,q} = L^{p_{\theta} ,q} (X) ,
\end{align}
where the equality holds with equivalence of norms.
In particular, when $q_0 = p_0$ and $q_1 = p_1$
\begin{align}
(L^{p_0} (X) , L^{p_1} (X))_{\theta , q} = L^{p_{\theta} ,q} (X) .
\end{align} 
\end{proposition}

\section{Linear Morawetz Estimates} \label{sec sharpmor}

The primary focus of this section will be generalizing the following estimate which is referred to as ``Sharp Morawetz estimates" in the paper of Ozawa and Rogers \cite{ozawa2013sharp}. The Morawetz estimate for the flat Schr\"odinger equation that Ozawa and Rogers reprove is originally due to Kato and Yajima \cite{kato1989some}. Ozawa and Rogers provide an argument that achieves the best possible constant for the inequality.

\begin{proposition}[Theorem 1.5 from \cite{ozawa2013sharp}] \label{thm ozawarogers}
 Suppose that $i \partial_t u+\Delta u=0$, on $\mathbb{R}^{d+1}$, with $d \geq 2$, and let $-1 / 2<s<d / 2-1$. Then
\begin{align*}
    \int_{\mathbb{R}} \int_{\mathbb{R}^d}|u(x, t)|^2 \frac{d x d t}{|x|^{2+2 s}} \leq C(d, s)\norm{u(0)}_{\dot{H}^s}^2 .
\end{align*}
Equality is attained if and only if $u(0)$ is radially symmetric.
\end{proposition}

A recent generalization of this theorem for higher order Laplacians with potential appears in a paper of Mizutani and Yao (\cite{mizutani2020global}, Theorem 1.1) and are referred to as ``Global Kato smoothing estimate". The focus of the result is on the higher order of the Laplacian, so the decay of the potential is not allowed to be critical. We can generalize Proposition~\ref{thm ozawarogers} for the critical inverse square potential case. The result relies heavily on the Riesz kernel estimates of \cite{killip2018sobolev} (Lemma 2.2) which, in turn, rely heavily on the Heat kernel estimates of Liskevich and Sobol \cite{liskevich2003estimates} (Theorem 1.2 and the discussion directly preceding Section 2) as well as Milman and Semenov \cite{milman2004global}. Generalizing these estimates for all $(\kappa, \delta)$ potentials requires standard heat kernel estimates.


First we state some tools that we will need to establish the main linear Morawetz estimate. Riesz kernel estimates that appear in \cite{killip2018sobolev} will be crucial. These estimates are important because they are for critically scaling potentials. 
In order to get from Riesz kernel estimates to Hardy estimates, we can use the standard Hardy-Littlewood-Sobolev Inequality.
For $\alpha \in (0, d)$, define $k_{\alpha}(x):= |x|^{-\alpha}$ and  $T_{\alpha} f(x) := (k_{\alpha} * f)(x)$. 

\begin{proposition}[Hardy-Littlewood-Sobolev Inequality]
    Let $1<p<q<\infty$ and 
    \begin{align*}
        \frac{d-\alpha}{d}= \frac{1}{p}-\frac{1}{q}.
    \end{align*}
    Then
    \begin{align*}
        \norm{T_{\alpha}f}_{L^q(\mathbb{R}^d)} \lesssim \norm{f}_{L^p(\mathbb{R}^d)}.
    \end{align*}
\end{proposition}

 Let
    \begin{align*}
        \sigma:= \frac{d-2}{2}-\frac{1}{2}\sqrt{(d-2)^2+4a}.
    \end{align*}
    then we have a lemma and proposition from \cite{killip2018sobolev} generalizing Hardy's inequality:

    \begin{lemma}[Lemma 2.2 from \cite{killip2018sobolev}]\label{lem KMVZZ Riesz}
    Let $d \geq 3$ and suppose $0<s<d$ and $d-s-2 \sigma>0$. Let $\mathcal{L}_a = -\Delta + \frac{a}{\abs{x}^2}$ with $a \geq - (\frac{d-2}{2})^2$.
    Then the Riesz potentials
    \begin{align*}
        \mathcal{L}_a^{-\frac{s}{2}}(x, y):=\frac{1}{\Gamma(s / 2)} \int_0^{\infty} e^{-t \mathcal{L}_a}(x, y) t^{\frac{s}{2}} \frac{d t}{t}
    \end{align*}
    satisfy
    \begin{align*}
        \mathcal{L}_a^{-\frac{s}{2}}(x, y) \sim |x-y|^{s-d}\left(\min(\tfrac{|x|}{|x-y|}, \tfrac{|y|}{|x-y|}, 1)\right)^{-\sigma} .
    \end{align*}
    \end{lemma}

    By comparison principles, we can extend these estimates to all admissible potentials.
\begin{definition}\label{def HV}
    Let $V$ be a $(\kappa, \delta)$ potential. We define the Schr\"odinger operator with potential 
    \begin{align}
        H_V := -\Delta + V .
    \end{align}
\end{definition}
    
    \begin{lemma}\label{lem HVbdd}
         Let $V$ be a $(\kappa, \delta)$ potential. Then there exist $\sigma_1=\sigma_1(\delta)$ and $\sigma_2=\sigma_2(\kappa)$, such that  
        \begin{align*}
           |x-y|^{s-d}\left(\min(\tfrac{|x|}{|x-y|}, \tfrac{|y|}{|x-y|}, 1)\right)^{-\sigma_2} \lesssim H^{-\frac{s}{2}}_{V}(x, y) \lesssim |x-y|^{s-d}\left(\min(\tfrac{|x|}{|x-y|}, \tfrac{|y|}{|x-y|}, 1)\right)^{-\sigma_1} .  
        \end{align*}
    \end{lemma}

    \begin{proof}[Proof of Lemma~\ref{lem HVbdd}]
       By definition, we can estimate $V$ by 
        \begin{align*}
            \left(- \left(\frac{d}{2}-1\right)^2+\delta \right)|x|^{-2} \leq V(x) \leq \kappa|x|^{-2} .
        \end{align*} 
        Let $a_1:=- \left(\frac{d}{2}-1\right)^2+\delta$ and $a_2:= \kappa$. 
        Then for $n \in \mathbb{Z}_+$ and $\psi \in C^{\infty}_c(\mathbb{R}^d)$ with $\psi \geq 0$,
        \begin{align*}
            e^{-\frac{1}{n} a_2 |x|^{-2}}\psi(x) \leq e^{-\frac{1}{n} V(x)}\psi(x)\leq e^{-\frac{1}{n} a_1 |x|^{-2}}\psi(x)
        \end{align*}
        and thus 
        \begin{align*}
            \left(e^{-\frac{1}{n} (-\Delta)}e^{-\frac{1}{n} a_2 |x|^{-2}} \right)^n\psi(x) \leq \left(e^{-\frac{1}{n} (-\Delta)}e^{-\frac{1}{n} V(x)} \right)^n\psi(x)\leq \left(e^{-\frac{1}{n} (-\Delta)}e^{-\frac{1}{n} a_1 |x|^{-2}} \right)^n\psi(x).
        \end{align*}
        The Trotter product formula (Corollary 2.7 in Chapter IV of \cite{engel2006short})  implies that 
            \begin{align*}
                e^{-t \mathcal{L}_{a_2}}(x, y) \leq e^{-tH_{V}}(x, y) \leq e^{-t \mathcal{L}_{a_1}}(x, y).
            \end{align*}
        Using these estimates with the identity
            \begin{align*}
                H_{V}^{-\frac{s}{2}}(x, y):=\frac{1}{\Gamma(s / 2)} \int_0^{\infty} e^{-t H_{V}}(x, y) t^{\frac{s}{2}} \frac{d t}{t}
            \end{align*}
       and Lemma~\ref{lem KMVZZ Riesz} provides the desired estimate.
    \end{proof}

We can observe that if $V\geq 0$, then we can take $\sigma_1=0$ and 
    \begin{align*}
        H^{-\frac{s}{2}}_{V}(x, y) \lesssim |x-y|^{s-d}.
    \end{align*}
This inequality with the Hardy-Littlewood-Sobolev inequality implies the following corollary.
    \begin{corollary}
     Suppose $d \geq 3$, $1<p<\infty$ and $0< s< \frac{d}{p}$. Let $V\geq 0$ be a $(\kappa, \delta)$ potential. Then
    \begin{align*}
        \left\||x|^{-s} f(x)\right\|_{L^p\left(\mathbb{R}^d\right)} \lesssim\left\|H_{V}^{\frac{s}{2}} f\right\|_{L^p\left(\mathbb{R}^d\right)}.
    \end{align*}
\end{corollary}

As shown in \cite{killip2018sobolev}, for potentials with more complicated Riesz kernel estimates, the relationship between the dimension and regularity is affected.

\begin{proposition}[Proposition 3.2 in \cite{killip2018sobolev}] \label{prop KMVZZ Hardy}
     Suppose $d \geq 3,0<s<d, d-s-2 \sigma>0$, and $1<p<\infty$. Then
    \begin{align*}
        \left\||x|^{-s} f(x)\right\|_{L^p\left(\mathbb{R}^d\right)} \lesssim\left\|\mathcal{L}_a^{\frac{s}{2}} f\right\|_{L^p\left(\mathbb{R}^d\right)}
    \end{align*}
holds if and only if
    \begin{align*}
        s+\sigma<\frac{d}{p}<d-\sigma .
    \end{align*}
\end{proposition}

These estimates can be extended to estimates for $(\kappa, \delta)$ potentials

\begin{proposition} \label{prop gen KMVZZ Hardy}
     Suppose $d \geq 3$, $V$ is a $(\kappa, \delta)$ potential and $\sigma_1:= \max(\frac{d-2}{2}-\sqrt{\delta}, 0)$. Let $0<s<d, d-s-2 \sigma_1>0$, and $1<p<\infty$. Then
    \begin{align*}
        \left\||x|^{-s} f(x)\right\|_{L^p\left(\mathbb{R}^d\right)} \lesssim\left\|H_V^{\frac{s}{2}} f\right\|_{L^p\left(\mathbb{R}^d\right)}
    \end{align*}
holds if and only if
    \begin{align*}
        s+\sigma_1<\frac{d}{p}<d-\sigma_1 .
    \end{align*}
\end{proposition}

Since kernels of the form $a|x|^{-2}$ for $a>-(d-2)^2/4$ fall under the umbrella of $(\kappa, \delta)$ potentials, we can proceed in our discussion considering only $(\kappa, \delta)$ potentials.

Our first goal is to prove a localized version of Hardy's inequality. Of course Hardy's inequality concerns nonlocal operators so complete localization is impossible. However, quantifying the extent to which the estimates are nonlocal will prove to be useful.
\begin{corollary}[Corollary of Proposition~\ref{prop gen KMVZZ Hardy}] \label{cor locHardy}
        Suppose $d \geq 3$, $V$ is a $(\kappa, \delta)$ potential and  $\sigma_1:= \max(\frac{d-2}{2}-\sqrt{\delta}, 0)$. Let $0<s<d$, $d-s-2 \sigma_1>0$, $1<p<\infty$, and $0<\varepsilon<d-p(s+\sigma_1)$. Let $B_R \subset \mathbb{R}^d$ be the ball of radius $R>0$ centered at the origin. Then there exists a $C=C(d, s, p, \varepsilon)$ such that 
    \begin{align*}
        \left\||x|^{-s} f\right\|_{L^p(B_R)} \leq C \left\|H_V^{\frac{s}{2}} f\right\|_{L^p(w_{R, \varepsilon})}
    \end{align*}
holds if and only if
    \begin{align*}
        s+\sigma_1<\frac{d}{p}<d-\sigma_1 .
    \end{align*}
    Where $w_{R}$ is a weight adapted to $B_R$ satisfying $w_{R, \varepsilon}(x) \sim \min( 1, (R^{-1}\mbox{dist}(B_R, x))^{-\varepsilon})$.
\end{corollary}

\begin{proof}[Proof of Corollary~\ref{cor locHardy}]
    It suffices to prove that 
    \begin{align*}
        \left\||x|^{-s}H_V^{-\frac{s}{2}} f(x)\right\|^p_{L^p(B_R)} \lesssim\left\| f\right\|^p_{L^p(w_{R})}.
    \end{align*}
    Lemma~\ref{lem HVbdd} implies that if $H_V^{-\frac{s}{2}}f(x)= \int K(x, y) f(y)\,dy$, then 
    \begin{align*}
        \left\||x|^{-s}H_V^{-\frac{s}{2}} f\right\|^p_{L^p(B_R)}&= \int_{B_R}\left||x|^{-s}\int K(x, y)f(y) \,dy \right|^p\,dx\\
        &\lesssim  \int_{B_R}\left||x|^{-s}\int K(x, y) \chi_{B_{10R}}(y)f(y) \,dy \right|^p\,dx+\int_{B_R}\left||x|^{-s}\int K(x, y) \chi_{B_{10R}^c}(y)f(y) \,dy \right|^p\,dx\\
        &\lesssim  \left\||x|^{-s}H_V^{-\frac{s}{2}} \chi_{B_{10R}}f\right\|^p_{L^p(\mathbb{R}^d)}+\int_{B_R}\left||x|^{-s}\int K(x, y) \chi_{B_{10R}^c}(y)f(y) \,dy \right|^p\,dx.
    \end{align*}
    Proposition~\ref{prop gen KMVZZ Hardy} implies 
    \begin{align*}
        \left\||x|^{-s}H_V^{-\frac{s}{2}} \chi_{B_{10R}}f\right\|^p_{L^p(\mathbb{R}^d)} \lesssim \norm{\chi_{B_{10R}}f}^p_{L^p(\mathbb{R}^d)}.
    \end{align*}
    Now by the estimate of Lemma~\ref{lem KMVZZ Riesz} we have
    \begin{align*}
        \int_{B_R}\left||x|^{-s}\int K(x, y) \chi_{B_{10R}^c}(y)f(y) \,dy \right|^p&\,dx\lesssim \int_{B_R}\left(\int_{B_{10R}^c}|x|^{-s}|x-y|^{s-d}\frac{|x|^{-\sigma_1}}{|x-y|^{-\sigma_1}} f(y) \,dy \right)^p\,dx\\
        &\lesssim \int_{B_R}|x|^{-p(s+\sigma_1)}\,dx \sup_{|x|\leq R}\left(\int_{B_{10R}^c}|x-y|^{(\sigma_1+s)-d} f(y) \,dy \right)^p.
    \end{align*}
    Now, since $s+\sigma_1 <d/p$, $\int_{B_R}|x|^{-p(s+\sigma_1)}\,dx \lesssim R^{\varepsilon_0}$ where $\varepsilon_0=d-p(s+\sigma_1)$. Therefore, for $\beta \in (0, \varepsilon_0)$
    \begin{align*}
        \int_{B_R}\left||x|^{-s}\int K(x, y) \chi_{B_{10R}^c}(y)f(y) \,dy \right|^p \, dx
        &\lesssim  R^{\varepsilon_0} \sup_{|x|\leq R}\left(\int_{B_{10R}^c}|x-y|^{(\sigma_1+s)-d} f(y) \,dy \right)^p\\
        &= R^{\varepsilon_0} \sup_{|x|\leq R}\left(\int_{B_{10R}^c}|x-y|^{(\sigma_1+s)-d/p-d/p'} f(y) \,dy \right)^p\\
        &= R^{\varepsilon_0} \sup_{|x|\leq R}\left(\int_{B_{10R}^c}|x-y|^{(\sigma_1+s)-d/p+\frac{1}{p}\beta-d/p'-\frac{1}{p}\beta} f(y) \,dy \right)^p.
    \end{align*}
    Now H\"older's inequality implies
    \begin{align*}
       & \quad R^{\varepsilon_0} \sup_{|x|\leq R}\left(\int_{B_{10R}^c}|x-y|^{(\sigma_1+s)-d} f(y) \,dy \right)^p
        \\
        &\lesssim R^{\varepsilon_0} \sup_{|x|\leq R}\left(\int_{B_{10R}^c}|x-y|^{-d-\frac{p'}{p}\beta} \,dy \right)^{p/p'}\left(\int_{B_{10R}^c}|x-y|^{p(\sigma_1+s)-d+\beta} |f(y)|^p \,dy \right)\\
        &\lesssim R^{\varepsilon_0-\beta} \sup_{|x|\leq R}\left(\int_{B_{10R}^c}|x-y|^{-\varepsilon_0 +\beta} |f(y)|^p \,dy \right) .
    \end{align*}
    Choose $\varepsilon = \varepsilon_0 - \beta$, then
    \begin{align}
        R^{\varepsilon_0-\beta} \sup_{|x|\leq R}\left(\int_{B_{10R}^c}|x-y|^{-\varepsilon_0 +\beta} |f(y)|^p \,dy \right) \lesssim \sup_{|x|\leq R}\left(\int_{B_{10R}^c} \parenthese{\frac{|x-y|}{R}}^{-\varepsilon} |f(y)|^p \,dy \right) \lesssim \norm{f}_{L^p (w_{R, \varepsilon})}^p .
    \end{align}
\end{proof}

At this point, we should re-establish some notation for Littlewood-Paley projectors. Given a $(\kappa, \delta)$ potential, $V$, and dyadic $N>0$, 
we will develop Littlewood-Paley theory in the same fashion as in \cite{killip2018sobolev}.

Let $\phi:[0, \infty) \rightarrow[0,1]$ be a smooth function such that
\begin{align*}
    \phi(\lambda)=1 \text { for } 0 \leq \lambda \leq 1 \text { and } \phi(\lambda)=0 \text { for } \lambda \geq 2 \text {. }
\end{align*}
For each dyadic number, $N$, we define
\begin{align*}
    \phi_N(\lambda):=\phi(\lambda / N) \text { and } \psi_N(\lambda):=\phi_N(\lambda)-\phi_{N / 2}(\lambda) \text {. }
\end{align*}
Furthermore, define the Littlewood-Paley projections as follows:
    \begin{align*}
        P_{\leq N}^V:=\phi_N\left(\sqrt{H_V}\right), \quad P_N^V:=\psi_N\left(\sqrt{H_V}\right), \quad \text { and } P_{>N}^V:=I-P_{\leq N}^V .
    \end{align*}

Define
    \begin{align*}
        \widetilde{P}^V_N f:=\left(e^{-H_V/N^2} - e^{-4H_V/N^2} \right)f .
    \end{align*}
If $P^a_N$ and $\widetilde{P}^a_N$ represent $P^V_N$ and $\widetilde{P}^V_N$, respectively, for $V=a|x|^{-2}$, then the following statement was proven. 
\begin{proposition}[Theorem 5.3 from \cite{killip2018sobolev}]\label{thm sobequiv}
     Fix $s \geq 0$. Let $1<p<\infty$ when $a \geq 0$ and let $r_0<p<r_0^{\prime}=\frac{d}{\sigma}$ when $-\left(\frac{d-2}{2}\right)^2 \leq a<0$. Then for any $f \in C_c^{\infty}\left(\mathbb{R}^d\right)$,
\begin{align*}
    \left\|\left(\sum_{N \in 2^{\mathbb{Z}}} N^{2 s}\left|P_N^a f\right|^2\right)^{\frac{1}{2}}\right\|_{L^p\left(\mathbb{R}^d\right)} \sim\left\|\mathcal{L}_a^{\frac{s}{2}} f\right\|_{L^p\left(\mathbb{R}^d\right)} \sim\left\|\left(\sum_{N \in 2^{\mathbb{Z}}} N^{2 s}\left|\left(\widetilde{P}_N^a\right)^k f\right|^2\right)^{\frac{1}{2}}\right\|_{L^p\left(\mathbb{R}^d\right)},
\end{align*}
provided the integer $k \geq 1$ satisfies $2 k>s$.
\end{proposition}

We will only need the equivalent statement for general potentials in the $p=2$ case. Thus, it suffices to appeal to the spectral theorem.

\begin{lemma}\label{lem sobequiv}
     Fix $s \in \mathbb{R}$. Then for any $f \in C_c^{\infty}\left(\mathbb{R}^d\right)$,
\begin{align*}
    \left\|\left(\sum_{N \in 2^{\mathbb{Z}}} N^{2 s}\left|P_N^V f\right|^2\right)^{\frac{1}{2}}\right\|_{L^2\left(\mathbb{R}^d\right)} \sim\left\|H_V^{\frac{s}{2}} f\right\|_{L^2\left(\mathbb{R}^d\right)}.
\end{align*}

\end{lemma}

On the way to proving the Morawetz estimate, we would like to establish basic frequency localization estimates so that they have some weak spatial localization. We start with a local Bernstein inequality.

\begin{lemma}\label{lem locbernclassic}
     Suppose $d \geq 3$, $N$ dyadic, $R>0$, and $1<p<\infty$.  If $R\geq 1/N$, then
    \begin{align*}
        \norm{NP_N f}_{L^p(B_R)} \lesssim \norm{\nabla f}_{L^p(w_{R})}.
    \end{align*}
    Where $w_{R}$ is a weight adapted to $B_R$ satisfying $w_{R}(x) \sim \inner{ x/R}^{-d}$.

\end{lemma}

\begin{proof}[Proof of Lemma~\ref{lem locbernclassic}]
    
    We observe that $NP_N|\nabla|^{-1}$ is a classic Mikhlin multiplier, as it is a Fourier multiplier of the form
    \begin{align*}
        \mathcal{F}\left(NP_N|\nabla|^{-1}f\right)(\xi) = N\chi_{\{|\xi| \sim N\}}\left(\sum_{j=1}^d|\xi_j|\right)^{-1}\widehat{f}(\xi).
    \end{align*} So (\cite{muscalu2013classical}, Proof of Theorem 8.2) there is a kernel $K_N$ such that 
    \begin{align*}
        |K_N(x)| &\lesssim |x|^{-d}
    \end{align*}
    for $|x|\geq R$ and 
    \begin{align*}
        NP_N|\nabla|^{-1}f(x)= (K_N * f)(x) . 
    \end{align*}
    Thus, an argument similar to what appears in Corollary~\ref{cor locHardy} gives us
    \begin{align*}
        \left\|NP_N f\right\|_{L^p(B_R)} \lesssim \left\||\nabla| f\right\|_{L^p(w_{R})}.
    \end{align*}
    Similarly, for each $j \in \{ 1,..., d\}$, $|\xi_j|/(i\xi_j)$ is a Mikhlin multiplier. Thus,
    \begin{align*}
        \left\|NP_N f\right\|_{L^p(B_R)} \lesssim \left\||\nabla| f\right\|_{L^p(w_{R})}\lesssim \sum_j\left\||\partial_{x_j}| f\right\|_{L^p(w_{R})}\lesssim \sum_j\left\|\partial_{x_j} f\right\|_{L^p(w_{R})}\lesssim \left\|\nabla f\right\|_{L^p(w_{R})} .
    \end{align*}
\end{proof}

With local Bernstein, we can now take our local smoothing estimate, Corollary~\ref{cor KMVZZ smoothing}, and prove a discretized version.

\begin{lemma}[Local Smoothing] \label{lem projlocsmooth}
     Suppose $d \geq 3$ and $V$ is a $(\kappa, \delta)$ potential. Let $R>0$, $N>0$, $u_0 \in L^2(\mathbb{R}^d)$, and $u = e^{-itH_V}u_0$.  If $P^V_N$ is the Littlewood-Paley projector for $H_V$, then
     \begin{align*}
         \norm{R^{-1/2}\chi_{B_R} NP^V_{N}u}^2_{L_{t,x}^2} \lesssim  (N+R^{-1})\norm{P^V_{N}u_0}^2_{L^2_x}.
     \end{align*}
\end{lemma}

\begin{proof}[Proof of Lemma~\ref{lem projlocsmooth}]
If $R < N^{-1}$, then the second estimate of Corollary~\ref{cor BPSTZ resolvent}, implies
    \begin{align*}
         \norm{R^{-1/2}\chi_{B_R} NP^V_{N}u}^2_{L_{t,x}^2} \lesssim \norm{R^{-1/2}\chi_{B_{N^{-1}}} NP^V_{N}u}^2_{L_{t,x}^2} \lesssim  R^{-1}\norm{P^V_{N}u_0}^2_{L^2_x}\lesssim  (N+R^{-1})\norm{P^V_{N}u_0}^2_{L^2_x}.
     \end{align*}

 When $R \geq N^{-1}$, first note that 
\begin{align*}
    \int_{B_R} \left|  NP^V_N u\right|^2 \lesssim \int_{B_R} \left|  NP_N P^V_N u\right|^2+\int_{B_R} \left|  N(P^V_N -P_N)P^V_Nu\right|^2 .
\end{align*}
Let $w_R$ be a weight satisfying $w_R(x) \sim \langle R^{-1}x \rangle^{-d}$. Since $R \geq N^{-1}$, Local Bernstein, Lemma~\ref{lem locbernclassic}, implies  
    \begin{align*}
        \int_{B_R} \left|  NP_N P^V_N u\right|^2\lesssim   \int \left|  \nabla P^V_N u\right|^2 w_R.
    \end{align*}
    By the estimate of the kernel of $P^V_N-P_N$  (Proposition 6.2 in \cite{killip2018sobolev}) and Hardy's inequality applied to $P^V_Nu \cdot w^{1/2}_R$ 
    \begin{align*}
        \int_{B_R} \left|  N(P^V_N-P_N) P^V_N u\right|^2\lesssim  \int \left|  |x|^{-1} P^V_N u\right|^2 w_R\lesssim  \int \left|  \nabla P^V_N u\right|^2 w_R +R^{-2}\int \left|  P^V_N u\right|^2 w_R . 
    \end{align*}
Corollary~\ref{cor KMVZZ smoothing}, Corollary~\ref{cor BPSTZ resolvent} along with Lemma~\ref{lem normequiv} implies that
    \begin{align*}
        \frac{1}{R}\int \int \left|  \nabla P^V_N u\right|^2 w_R +R^{-2}\int \left|  P^V_N u\right|^2 w_R\,dt&\lesssim \norm{P^V_Nu_0}_{L_x^2}\norm{\nabla P^V_Nu_0}_{L_x^2}+R^{-1}\norm{P^V_Nu_0 }_{L_x^2}^2 \\
        &\lesssim \norm{P^V_Nu_0}_{L_x^2}\norm{(-\Delta)^{1/2} P^V_Nu_0}_{L_x^2}+R^{-1}\norm{P^V_Nu_0 }_{L_x^2}^2 \\
        &\lesssim \norm{P^V_Nu_0}_{L_x^2}\norm{H_V^{1/2} P^V_Nu_0}_{L_x^2}+R^{-1}\norm{P^V_Nu_0 }_{L_x^2}^2\\
        &\lesssim N\norm{P^V_Nu_0}^2_{L_x^2}+R^{-1}\norm{P^V_Nu_0 }_{L_x^2}^2 .
    \end{align*}
\end{proof}

It should be noted that the proof of Lemma~\ref{lem projlocsmooth} should be attainable without introducing the classical Littlewood-Paley projector in the argument. However, one would probably need to develop two estimates: (1) A generalization of the decay estimate for the kernel of $NP^V_NH^{-1}_V$ which is mostly done in \cite{killip2018sobolev}, and  (2) a generalization of the local smoothing estimate (Corollary~\ref{cor KMVZZ smoothing}) with a nonlocal operator $H_V^{1/2}$ replacing $\nabla$. Unlike in Ozawa-Rogers \cite{ozawa2013sharp}, we do not attempt to determine a sharp constant for the estimate and we do not treat the $d=2$ case. The proof is a discretization argument that ends with an application of Schur's test.

 

\genmorawetz*

\begin{proof}[Proof of Theorem~\ref{thm genmorawetz}]

Suppose that $u:[0, \infty) \to H^s(\mathbb{R}^d)$ satisfies $i \partial_t u+\Delta u-Vu=0$ and let $u_0:= u(0)$. 
Let $\{ \chi_{\ell}\}_{\ell \in \mathbb{Z}}$ be a smooth partition of unity of $\mathbb{R}^d$ such that 
    \begin{align*}
        \mbox{supp}(\chi_{\ell}) \subset \{ x \in \mathbb{R}^d ~:~ 2^{\ell-1} \leq |x| \leq 2^{\ell+1}\}.
    \end{align*}

    
Now let $s$ satisfy 
    \begin{align*}
        -1/2<s<d/2-1-\sigma_1
    \end{align*}
    and let $s_0$ satisfy 
    \begin{align*}
        s<s_0<d/2-1-\sigma_1.
    \end{align*}
    
    \begin{align*}
        \norm{|x|^{-(1+s)} u}^2_{L_{t,x}^2} &\leq 10 \sum_{\ell \in \mathbb{Z}}\norm{\chi_{\ell}|x|^{-(1+s)} u}^2_{L_{t,x}^2}\leq 100 \sum_{\ell \in \mathbb{Z}}2^{-2(1+s)\ell}\norm{\chi_{\ell} u}^2_{L_{t,x}^2}\\
        &\leq 1000 \sum_{\ell \in \mathbb{Z}}2^{-2(1+s)\ell}\left(\norm{\chi_{\ell} P^V_{\leq 2^{-\ell}}u}^2_{L_{t,x}^2} +\norm{\chi_{\ell} P^V_{\geq 2^{-\ell}}u}^2_{L_{t,x}^2} \right) .
    \end{align*}
We handle the two terms in the last inequality separately.  First, the case in which the restriction in space is at a smaller scale than the smoothing of the Littlewood-Paley projector is considered. We note that $\mbox{supp}\chi_{\ell} \subset B(0, 2^{\ell+1})$ and when $k \leq - \ell$,  $\mbox{supp}(\chi_{\ell}) \subset B(0, 2^{-k+1}).$ Therefore, 
    \begin{align*}
        \sum_{\ell \in \mathbb{Z}}2^{-2(1+s)\ell}\norm{\chi_{\ell} P^V_{\leq 2^{-\ell}}u}^2_{L_{t,x}^2} &\lesssim \sum_{\ell \in \mathbb{Z}}2^{-2(1+s)\ell}\left( \sum_{k \leq -\ell}\norm{\chi_{\ell} P^V_{2^k}u}_{L_{t,x}^2} \right)^2\\
        &\lesssim \sum_{\ell \in \mathbb{Z}}2^{-2(1+s)\ell}\left( \sum_{k \leq -\ell}2^{(s_0+1)\ell}\norm{2^{-(s_0+1)\ell}\chi_{\ell} P^V_{2^k}u}_{L_{t,x}^2(B(0, 2^{-k+1}))} \right)^2 .
    \end{align*}

Now the local Hardy inequality, Corollary~\ref{cor locHardy}, dictates that
    \begin{align*}
        \norm{2^{-(s_0+1)\ell}\chi_{\ell} P^V_{2^k}u}_{L^{2}_x(B(0, 2^{-k+1}))} &\lesssim \norm{|x|^{-(s_0+1)} P^V_{2^k}u}_{L^{2}_x(B(0, 2^{-k+1}))} \\
        &\lesssim \norm{H_V^{(s_0+1)} P^V_{2^k}u}_{L^{2}_x(w_{2^{-k+1}})}\lesssim 2^{k(s_0+1)}2^{-k}\norm{ 2^kP^V_{k}u}_{L^2_x(w_{2^{-k+1}})}\\
        & \lesssim 2^{k(s_0+1)}2^{-k}\norm{ |x|^{-1}P^V_{2^k}u}_{L^2_x(\mathbb{R}^d)}
    \end{align*}
We can now use the Morawetz estimate of \cite{burq2004strichartz} (Corollary~\ref{cor BPSTZ resolvent}) to get 
    \begin{align*}
        \sum_{\ell \in \mathbb{Z}}2^{-2(1+s)\ell}\norm{\chi_{\ell} P^V_{\leq 2^{-\ell}}u}^2_{L_{t,x}^2} &\lesssim\sum_{\ell \in \mathbb{Z}}2^{-2(1+s)\ell}\left( \sum_{k \leq -\ell}2^{(s_0+1)\ell}2^{k(s_0+1)}2^{-k}\norm{ |x|^{-1}P^V_{2^k}u}_{L_{t,x}^2} \right)^2\\
        &=\sum_{\ell \in \mathbb{Z}}2^{-2(1+s)\ell}\left( \sum_{k \leq -\ell}2^{(s_0+1)\ell}2^{k(s_0+1)}2^{-k}\norm{ |x|^{-1}e^{-itH_V}(P^V_{2^k}u_0)}_{L_{t,x}^2} \right)^2\\
        &\lesssim \sum_{\ell \in \mathbb{Z}}2^{-2(1+s)\ell}\left( \sum_{k \leq -\ell}2^{\ell (s_0+1)}2^{k(s_0+1)}2^{-k}\norm{P^V_{2^k}u_0}_{L^2_x} \right)^2
    \end{align*}
By simple multiplication by one and rearrangement, we now have
    \begin{align*}
        \sum_{\ell \in \mathbb{Z}}2^{-2(1+s)\ell}\norm{\chi_{\ell} P^V_{\leq 2^{-\ell}}u}^2_{L_{t,x}^2} &\lesssim \sum_{\ell \in \mathbb{Z}}2^{-2s\ell}\left( \sum_{k \leq -\ell}2^{-\ell}2^{(s_0+1)\ell}2^{k(s_0+1)}2^{-k}\norm{P^V_{2^k}u_0}_{L^2_x} \right)^2\\
        &= \sum_{\ell \in \mathbb{Z}}2^{-2s\ell}\left( \sum_{k \leq -\ell}2^{(k+\ell)s_0}\norm{P^V_{2^k}u_0}_{L^2_x} \right)^2\\
        &= \sum_{\ell \in \mathbb{Z}}\left( \sum_{k \leq -\ell}2^{(k+\ell)(s_0-s)}2^{sk}\norm{P^V_{2^k}u_0}_{L^2_x} \right)^2
    \end{align*}
    Schur's Test (Lemma~\ref{lem:schur}) with $f(k) := 2^{sk}\norm{P^V_{2^k}u_0}_{L^2_x}$ and kernel $K(k, \ell):= 2^{(k+\ell)(s_0-s)}$. Since $s_0-s>0$, we have
    \begin{align*}
        \sum_{\ell \in \mathbb{Z}}2^{-2(1+s)\ell}\norm{\chi_{\ell} P^V_{\leq 2^{-\ell}}u}^2_{L_{t,x}^2} &\lesssim  \sum_{\ell \in \mathbb{Z}}\left( \sum_{k \leq -\ell}2^{(k+\ell)(s_0-s)}2^{sk}\norm{P^V_{2^k}u_0}_{L^2_x} \right)^2\\
        &\lesssim  \sum_{k \in \mathbb{Z}}2^{2sk}\norm{P^V_{2^k}u_0}_{L^2_x}^2 . 
    \end{align*}
    
For the second term,  we have 
    \begin{align*}
        \sum_{\ell \in \mathbb{Z}}2^{-2(1+s)\ell}\norm{\chi_{\ell} P^V_{\geq 2^{-\ell}}u}^2_{L_{t,x}^2} &\lesssim \sum_{\ell \in \mathbb{Z}} 2^{-2(s+1/2)\ell} \sum_{k \geq -\ell} 2^{-2sk} 2^{2sk}\norm{2^{-\ell/2}\chi_{\ell} P^V_{2^k}u}^2_{L_{t,x}^2}\\
        &\lesssim \sum_{\ell \in \mathbb{Z}}  \sum_{k \geq -\ell} 2^{-2(s+\frac{1}{2})(k+\ell)} 2^{2(s+\frac{1}{2})k}\norm{2^{-\ell/2}\chi_{\ell} P^V_{2^k}u}^2_{L_{t,x}^2} .
    \end{align*}
Now Local Smoothing, Lemma~\ref{lem projlocsmooth}, implies
        \begin{align*}
        \sum_{\ell \in \mathbb{Z}}2^{-2(1+s)\ell}\norm{\chi_{\ell} P^V_{\geq 2^{-\ell}}u}^2_{L_{t,x}^2} &\lesssim\sum_{\ell \in \mathbb{Z}}  \sum_{k \geq -\ell} 2^{-2(s+\frac{1}{2})(k+\ell)} 2^{2(s+\frac{1}{2})k}\norm{2^{-\ell/2}\chi_{\ell} P^V_{2^k}u}^2_{L_{t,x}^2}\\
        &\lesssim \sum_{\ell \in \mathbb{Z}}  \sum_{k \geq -\ell} 2^{-2(s+\frac{1}{2})(k+\ell)} 2^{2(s+\frac{1}{2})k}(2^{k}+2^{-\ell})2^{-2k}\norm{P^V_{2^k}u_0}^2_{L^2_x}\\
        &\lesssim \sum_{\ell \in \mathbb{Z}}  \sum_{k \geq -\ell} 2^{-2(s+\frac{1}{2})(k+\ell)} 2^{2sk}\norm{P^V_{2^k}u_0}^2_{L^2_x}\\
        &= \sum_{k \in \mathbb{Z}}   2^{2sk}\norm{P^V_{2^k}u_0}^2_{L^2_x}\left(\sum_{\ell \geq -k} 2^{-2(s+\frac{1}{2})(k+\ell)} \right) .
    \end{align*}
Since $s>-\frac{1}{2}$, we have both terms bounded in a satisfactory way. Now using Lemma~\ref{lem sobequiv}, we have 
    \begin{align*}
        \norm{|x|^{-(1+s)} u}^2_{L_{t,x}^2} \lesssim \sum_{k \in \mathbb{Z}}   2^{2sk}\norm{P^V_{2^k}u_0}^2_{L^2_x}  \sim \norm{H_V^{\frac{s}{2}}u_0}^2_{L^2_x} . 
    \end{align*}

    We then finish the proof of Theorem~\ref{thm genmorawetz}.
\end{proof}

Using Lemma~\ref{lem normequiv} we get the equivalent statement with $H_V$ replaced with $-\Delta$ for a smaller range of $s$:
\begin{corollary}\label{cor flatgenmorawetz}
Suppose $d \geq 3$, $V$ is a $(\kappa, \delta)$ potential and  $\sigma_1:= \max(\frac{d-2}{2}-\sqrt{\delta}, 0)$.
 
Let $u = e^{-itH_V}u_0$ and $-1/2<s<\min(2, d/2-1-\sigma_1)$. Then
\begin{align*}
    \int_{\mathbb{R}} \int_{\mathbb{R}^d}|u(x, t)|^2 \frac{d x d t}{|x|^{2+2 s}} \leq C(d, s)\norm{(-\Delta)^{\frac{s}{2}}u_0}_{L^2}^2 .
\end{align*}
\end{corollary}


\section{Endpoint of Bilinear Christ-Kiselev Lemmas}
\label{sec CK}
In this section, we present two bilinear Christ-Kiselev type Lemmas at the endpoint, which will be used in the next sections.

In general, Christ–Kiselev lemma misses the double $L_t^2$-type estimates (see Lemma~\ref{lem:christkiselev}), however, when handling Schr\"odinger operators, it is sometimes anticipated to obtain the missing double endpoint. In particular, we refer to a generalized double endpoint Christ–Kiselev type estimate for inhomogeneous Strichartz estimates in \cite{GLNY18}. 
\begin{proposition}[Theorem 1.3 in \cite{GLNY18} with $a=2$]\label{prop GLNY18}
Let $d \geq 3$. Assume $p,r > \frac{4d-2}{2d-3}$. Then 
\begin{align}
\norm{\int_0^t e^{i (t-s) \Delta} P_0 g(s) \, ds}_{L_t^2 \mathcal{L}_{\rho}^p L_{\omega}^2} \lesssim \norm{g}_{L_t^2 \mathcal{L}_{\rho}^{r'} L_{\omega}^2} ,
\end{align}
where $L^p (\R^d)$ denotes the usual Lebesgue space, and $\mathcal{L}^p (\R^+) = L^p (\R^+ :r^{d-1} dr)$. 
\begin{align}
\norm{f}_{\mathcal{L}_r^p L_{\omega}^q} = \norm{ \norm{f(r\omega)}_{L_{\omega}^q}}_{\mathcal{L}_r^p} , 
\end{align}
with $x = r \omega$, $\omega \in \mathbb{S}^{d-1}$. Here $P_0$ is the projection to frequencies $\abs{\xi} \sim 1$.
\end{proposition}

Inspired by Proposition~\ref{prop GLNY18}, we will show the following bilinear versions of the double-endpoint Christ–Kiselev lemmas. 
\CK*

We split the proof of Theorem~\ref{thm CK} into Lemma~\ref{lem CK} treating the linear-Duhamel interaction and Lemma~\ref{lem CK'} treating the Duhamel-Duhamel interaction.

\begin{lemma}[Endpoint of bilinear Christ-Kiselev Lemma, linear-Duhamel]\label{lem CK}
Let $d \geq 3$. For every $0 < \varepsilon <1$, the following bounds hold
\begin{align}\label{eq CK1}
\norm{e^{it\Delta} P_N u_0 \parenthese{ \int_{-\infty}^t e^{i(t-s)\Delta} P_L F(s) \, ds }}_{L_{t,x}^2} \lesssim C(N,L)^{1-\varepsilon} \min \{ N,L\}^{\frac{d-2}{2} \varepsilon}  \norm{u_0}_{L^2} \norm{F}_{L_t^{2} L_x^{\frac{2d} , {d+2},2}} 
\end{align}
and
\begin{align}\label{eq CK2}
\norm{e^{it\Delta} P_N u_0 \parenthese{ \int_{-\infty}^t e^{i(t-s)\Delta} P_L F(s) \, ds }}_{L_{t,x}^2} \lesssim C(N,L)^{1-\varepsilon} \min \{N,L\}^{\frac{d-2}{2} \varepsilon} \norm{u_0}_{L^2} \norm{F}_{L_t^{2} L_x^{\frac{2d}{d+2}}} .
\end{align}
\end{lemma}

\begin{remark}
Note that the coefficient $C(N,L)$ in the $L_{t,x}^2$ bilinear Strichartz estimate setting takes the form
\begin{align}
C(N,L) = \min\{N,L\}^{\frac{d-2}{2}} \parenthese{\frac{\min \{N,L\}}{\max \{ N,L\}} }^{\frac{1}{2}} . 
\end{align}
The upper bounds in Lemma~\ref{lem CK} are in fact
\begin{align}
C(N,L)^{1-\varepsilon} \min \{N,L\}^{\frac{d-2}{2} \varepsilon} = \min\{N,L\}^{\frac{d-2}{2}} \parenthese{\frac{\min \{N,L\}}{\max \{ N,L\}} }^{\frac{1 - \varepsilon}{2}} , 
\end{align}
hence we only lose an arbitrary small power of the $\frac{\min}{\max}$ ratio.
\end{remark}

\begin{corollary}
Let $V$ be a $(\kappa , \delta)$ potential and $v = e^{-itH_V} v_0$.
After we obtain the endpoint of Christ-Kiselev Lemma (Lemma~\ref{lem CK}), by taking $F = Vv$, H\"older inequality for Lorentz norms, and Corollary~\ref{cor BPSTZ resolvent}, we have
\begin{align}
\norm{Vv}_{L_t^2 L_x^{\frac{2d}{d+2},2} } \lesssim \norm{|x| V}_{L_t^{\infty} L_x^{d,\infty}} \norm{|x|^{-1}v}_{L_{t}^2 L_x^{2,2}}  = \norm{|x| V}_{L_t^{\infty} L_x^{d,\infty}} \norm{|x|^{-1}v}_{L_{t,x}^2} \lesssim \norm{v_0}_{L^2} .
\end{align}
Hence
\begin{align}
\norm{e^{it\Delta} P_N u_0 \parenthese{ \int_{-\infty}^t e^{i(t-s)\Delta} P_L (V v)(s) \, ds }}_{L_{t,x}^2} & \lesssim C(N,L)^{1-\varepsilon} \min \{ N,L\}^{\frac{d-2}{2} \varepsilon} \norm{u_0}_{L^2} \norm{Vv}_{L_t^{2} L_x^{\frac{2d}{d+2},2}} \\
& \lesssim C(N,L)^{1-\varepsilon} \min \{ N,L\}^{\frac{d-2}{2} \varepsilon} \norm{u_0}_{L^2} \norm{v_0}_{L^2} . 
\end{align}  
\end{corollary}

\begin{proof}[Proof of Lemma~\ref{lem CK}]
We only show the first bound \eqref{eq CK1} and the second bound follows from the embedding $L^{\frac{2d}{d+2}} \subset L^{\frac{2d}{d+2},2}$.

The proof is inspired by the proof of Lemma 2.2 in \cite{Vil07} and that of Theorem 1.3 in \cite{GLNY18}. The main difference in our argument is that we use a dyadic decomposition of the time interval adapted to the smaller frequency.

Let $m = \min\{N,L\}$ denote the smaller frequency. We decompose the interval $(-\infty,t)$ into dyadic subintervals
\begin{align}
(-\infty,t) = \bigcup_{j=-\infty}^{\infty} I_j,
\qquad
I_j = [t-2\delta_j,t-\delta_j],
\qquad
\delta_j = m^{-2}2^j.
\end{align}
We note that the length of each subinterval, $|I_j|= m^{-2}2^j$, are adapted to the smaller frequency $m$. 
\begin{align}
& \quad \norm{e^{it\Delta} P_N u_0 \parenthese{ \int_{-\infty}^t e^{i(t-s)\Delta} P_L F (s) \, ds }}_{L_{t,x}^2} = \norm{e^{it\Delta} P_N u_0  \, e^{it\Delta}P_L \parenthese{\sum_j \int_{I_j} e^{-is\Delta}  F(s) \, ds }}_{L_{t,x}^2} \\
& = \norm{e^{it\Delta} P_N u_0  \, \sum_j  e^{it\Delta}P_L \parenthese{ \int_{I_j} e^{-is\Delta} F(s) \, ds }}_{L_{t,x}^2} \leq \sum_j \underbrace{\norm{e^{it\Delta} P_N u_0  \, e^{it\Delta}P_L \parenthese{\int_{I_j} e^{-is\Delta}  F(s) \, ds }}_{L_{t,x}^2} }_{T_j} .
\end{align}
We will separate the discussion into two cases: $j < 0$ and $j \geq 0$. The plan is as follows.
\begin{enumerate}[(A)]
\item 
When $j \geq 0$, we obtain the following two bounds for $T_j$
\begin{align}
T_j & \lesssim C(N,L) 2^{j \gamma}  m^{-2 \gamma} \norm{u_0}_{L^2} \norm{F}_{L_t^{2} L_x^{\frac{2d}{d+2}+}}  \label{eq BDD1} ,\\
T_j & \lesssim  m^{d-2} 2^{-j\frac{d-2}{2}} \norm{u_0}_{L^2}  \norm{F}_{L_t^2 L_x^1} , \label{eq BDD2}
\end{align}
where $\gamma = \frac{(d-2)\varepsilon}{4(1-\varepsilon)}$. Then an interpolation between these bounds gives
\begin{align}
T_j \lesssim C(N,L)^{1-\varepsilon} m^{\frac{d-2}{2} \varepsilon} 2^{-j \frac{d-2}{4} \varepsilon}   \norm{u_0}_{L^2} \norm{F}_{L_t^{2} L_x^{\frac{2d}{d+2},2}} . \label{eq1}
\end{align}
The interpolation we used here is
\begin{align}
(L_t^2 L_x^{r} , L_t^2 L_x^{1})_{\varepsilon ,2} = L_t^2 L_x^{\frac{2d}{d+2},2} \quad \text{where} \quad \frac{1-\varepsilon}{r} + \varepsilon = \frac{d+2}{2d} ,
\end{align}
from Lions-Peetre interpolation theorem (see Section~\ref{sec LP}). Interpolation with weight $\varepsilon$ on \eqref{eq BDD1} gives \eqref{eq1}.

\item
When $j < 0$, we show the following bound holds
\begin{align}
T_j & \lesssim C(N,L)^{1-\varepsilon} \min \{N,L\}^{\frac{d-2}{2} \varepsilon} 2^{\frac{j}{2} \varepsilon}   \norm{u_0}_{L^2} \norm{F}_{L_t^{2} L_x^{\frac{2d}{d+2},2}} , \label{eq BDD3} 
\end{align}

\item
Both $j < 0$ and $j \geq 0$ bounds allow us to sum $T_j$ over $j$ and get
\begin{align}
\norm{e^{it\Delta} P_N u_0 \parenthese{ \int_{-\infty}^t e^{i(t-s)\Delta} P_L F(s) \, ds }}_{L_{t,x}^2} \lesssim C(N,L)^{1-\varepsilon} \min \{N,L\}^{\frac{d-2}{2} \varepsilon}  \norm{u_0}_{L^2} \norm{F}_{L_t^{2} L_x^{\frac{2d}{d+2},2}} , 
\end{align}
hence concludes the desired bound.
\end{enumerate}

Note that in the plan above and in the proof of Bound \eqref{eq BDD1}, we choose 
\begin{align}
\gamma = \frac{(d-2)\varepsilon}{4(1-\varepsilon)} , \quad \frac{1}{q} = \frac{1}{2} + \gamma, \quad \frac{1}{r} = \frac{d+2}{2d} - \frac{2\gamma}{d} , 
\end{align}
for small enough $\varepsilon >0$. 
Here $(q,r)$ is a dual Strichartz pair and we write instead $q= 2-$ and $r = \frac{2d}{d+2}+$ in the calculation for simplicity. We will treat both $\gamma$ and $\varepsilon$ as very small parameters.

\noindent{\bf The first bound, \eqref{eq BDD1}:}

For each $T_j$, using duality, we write
\begin{align}
T_j & = \norm{e^{it\Delta} P_N u_0  \, e^{it\Delta}P_L \parenthese{\int_{I_j} e^{-is\Delta}  F(s) \, ds }}_{L_{t,x}^2} \\
& = \sup_h \abs{\iint e^{it\Delta} P_N u_0  \, e^{it\Delta}P_L \parenthese{\int_{I_j} e^{-is\Delta}  F(s) \, ds } \overline{h}(t,x) \, dxdt} \\
& = \sup_h \abs{ \sum_{a \in \delta_j \Z} \int_{[a,a+ \delta_j]} \int e^{it\Delta} P_N u_0  \, e^{it\Delta} P_L \parenthese{\int_{I_j} e^{-is\Delta}  F(s) \, ds } \overline{h}(t,x) \, dxdt} , 
\end{align}
where the supremum is taken over $h$ with $\norm{h}_{L_{t,x}^2} = 1$. 
Note that in the last line, we break the interval $(-\infty , \infty)$ into $\cup [a, a+\delta_j]$.

Now on each $[a, a+\delta_j]$, by H\"older inequality, Christ-Kiselev Lemma (Lemma~\ref{lem:christkiselev})
\begin{align}
& \quad \abs{\int_{[a,a+\delta_j]} \int e^{it\Delta} P_N u_0  \, e^{it\Delta}P_L \parenthese{\int_{I_j} e^{-is\Delta}  F(s) \, ds } \overline{h}(t,x) \, dxdt } \\
& \lesssim \norm{e^{it\Delta} P_N u_0  \, e^{it\Delta}P_L \parenthese{\int_{I_j} e^{-is\Delta} \chi_{[a-2\delta_j, a]} (s) F(s) \, ds }}_{L_t^{2} L_x^2} \norm{\chi_{[a,a+\delta_j]} (t) h(t,x)}_{L_t^{2}L_x^2}\\
& \lesssim C(N,L) \norm{u_0}_{L^2} \norm{\chi_{[a-2\delta_j ,a]}(s) F}_{L_s^{2-} L_x^{\frac{2d}{d+2}+}} \norm{\chi_{[a,a+\delta_j]} (t) h}_{L_{t,x}^2} \\
& \lesssim C(N,L) \delta_j^{\gamma} \norm{u_0}_{L^2} \norm{\chi_{[a-2\delta_j ,a]}(s) F}_{L_s^{2} L_x^{\frac{2d}{d+2}+}} \norm{\chi_{[a,a+\delta_j]} (t) h}_{L_{t,x}^2} .
\end{align}
Note that we can write
\begin{align}
\int_{I_j} = \int_{t-2\delta_j}^{t-\delta_j} =\int_{s<t-\delta_j}-\int_{s<t-2\delta_j} .
\end{align}
then apply Christ-Kiselev for both integrals. Note that for any spatial Banach space $Y$, we are able to sum
\begin{align}
\sum_a \norm{\chi_{[a-2\delta_j, a]} F}_{L_t^2 Y}^2 = 2 \norm{F}_{L_t^2 Y}^2 . 
\end{align}

Then when summing in $a \in \delta_j \Z$, by we have
\begin{align}
& \quad \sum_a \abs{ \int_{[a,a+\delta_j]} \int e^{it\Delta} P_N u_0  \, e^{it\Delta}P_L \parenthese{\int_{I_j} e^{-is\Delta}  F(s) \, ds } \overline{h}(t,x) \, dxdt} \\
& \lesssim C(N,L) \delta_j^{\gamma}  \norm{u_0}_{L^2} \norm{ \norm{\chi_{[a-2\delta_j ,a]}(s) F}_{L_s^{2} L_x^{\frac{2d}{d+2}+}} }_{\ell_a^2} \norm{ \norm{\chi_{[a,a+\delta_j]} (t) h}_{L_{t,x}^2} }_{\ell_a^2} \\
& \lesssim C(N,L) \delta_j^{\gamma}  \norm{u_0}_{L^2} \norm{F}_{L_t^{2} L_x^{\frac{2d}{d+2}+}} \norm{\norm{\chi_{[a,a+\delta_j]} (t) h}_{\ell_a^1 L_{x}^2} }_{L_t^2} \\
& \lesssim C(N,L) \delta_j^{\gamma}  \norm{u_0}_{L^2} \norm{F}_{L_t^{2} L_x^{\frac{2d}{d+2}+}}  \norm{h}_{L_{t,x}^2} .
\end{align}

Note that for the unrestricted (in time) estimate, we have
\begin{align}
\norm{e^{it\Delta} P_N u_0 \parenthese{ \int_{\mathbb{R}} e^{i(t-s)\Delta} P_L F(s) \, ds }}_{L_{t,x}^2} & \lesssim C(N,L) \norm{u_0}_{L^2} \norm{\int_{\mathbb{R}} e^{-is\Delta}  F(s) \, ds  }_{L^2} \\
& \lesssim C(N,L) \norm{u_0}_{L^2} \norm{F}_{L_t^{2-} L_x^{\frac{2d}{d+2}+}} , 
\end{align}
then by Christ-Kiselev Lemma
\begin{align}
\norm{e^{it\Delta} P_N u_0  \, e^{it\Delta}P_L \parenthese{\int_{I_j} e^{-is\Delta}  F(s) \, ds }}_{L_{t,x}^2} \lesssim C(N,L) \norm{u_0}_{L^2} \norm{F}_{L_t^{2-} L_x^{\frac{2d}{d+2}+}} . 
\end{align}

\noindent{\bf The second bound, \eqref{eq BDD2}:}

By H\"older inequality and dispersive estimate, we write
\begin{align}
T_j & = \norm{\int_{I_j} e^{it\Delta} P_N u_0  \, e^{i(t-s)\Delta} P_L F(s) \, ds }_{L_{t,x}^2} \\
& \leq \norm{\int_{I_j} \norm{e^{it\Delta} P_N u_0  \, e^{i(t-s)\Delta} P_L F(s)}_{L_x^2} \, ds }_{L_{t}^2}  \\
& \leq \norm{\int_{I_j} \norm{e^{it\Delta} P_N u_0}_{L_x^2}  \norm{e^{i(t-s)\Delta} P_L F(s)}_{L_x^{\infty}} \, ds }_{L_{t}^2} \\
& \lesssim \norm{\int_{I_j} \norm{P_N u_0}_{L_x^2}  \frac{1}{\abs{t-s}^{\frac{d}{2}}} \norm{P_L F(s)}_{L_x^{1}} \, ds }_{L_{t}^2} \\
& \lesssim \norm{u_0}_{L_x^2} \norm{\int_{I_j}   \frac{1}{\abs{t-s}^{\frac{d}{2}}} \norm{P_L F(s)}_{L_x^{1}} \, ds }_{L_{t}^2} .
\end{align}

For the second factor above, we use duality again
\begin{align}\label{eq4}
& \norm{\int_{I_j} \frac{1}{\abs{t-s}^{\frac{d}{2}}} \norm{P_L F(s)}_{L_x^{1}} \, ds }_{L_{t}^2} = \sup_h \abs{\int_{\R} \int_{I_j} \frac{1}{\abs{t-s}^{\frac{d}{2}}} \norm{P_L F(s)}_{L_x^{1}}  h(t) \, dsdt } , 
\end{align}
where the supremum is taken over $h$ with $\norm{h}_{L_{t}^2} = 1$ and $h \geq 0$.

Next, we subdivide $\R = \cup_{a \in \delta_j \Z} [a, a+\delta_j]$. 
Then by Fubini and H\"older, we have
\begin{align}
& \quad \int_{[a,a+\delta_j]} \int_{I_j} \frac{1}{\abs{t-s}^{\frac{d}{2}}} \norm{P_L F(s)}_{L_x^{1}} h(t) \, ds dt \\
& = \iint \chi_{\parallelogramm} h(t) \frac{1}{\abs{t-s}^{\frac{d}{2}}} \norm{P_L F(s)}_{L_x^{1}} \, dtds \\
& = \int h(t) \parenthese{\int \chi_{\parallelogramm}  \frac{1}{\abs{t-s}^{\frac{d}{2}}} \norm{P_L F(s)}_{L_x^{1}} \, ds} \, dt \\
& \lesssim \delta_j^{-d/2} \int \chi_{[a, a+\delta_j]} h(t) \parenthese{\int \chi_{\parallelogramm}   \norm{P_L F(s)}_{L_x^{1}} \, ds} \, dt \\
& \lesssim \delta_j^{-d/2} \int  \chi_{[a, a+\delta_j]} (t) h(t) \norm{\chi_{\parallelogramm} }_{L_s^2}  \norm{\chi_{[a -2\delta_j ,a ]} (s) P_L F(s)}_{L_s^2L_x^1} \, dt \\
& \lesssim \delta_j^{-d/2} \delta_j^{1/2}\int \chi_{[a, a+\delta_j]} (t) h(t)  \norm{\chi_{[a -2\delta_j ,a ]} (s) P_L F(s)}_{L_s^2L_x^1} \, dt \\
& \lesssim \delta_j^{-d/2} \delta_j^{1/2} \norm{\chi_{[a -2\delta_j ,a ]} (s) P_L F(s)}_{L_s^2L_x^1} \int \chi_{[a, a+\delta_j]} (t) h(t) \, dt \\
& \lesssim \delta_j^{-d/2} \delta_j^{1/2} \norm{\chi_{[a -2\delta_j ,a ]} (s) P_L F(s)}_{L_s^2L_x^1} \norm{\chi_{[a, a+\delta_j]} (t) h(t) }_{L_t^2}  \norm{\chi_{[a, a+\delta_j]} (t)}_{L_t^2}\\
& \lesssim \delta_j^{-d/2} \delta_j \norm{\chi_{[a -2\delta_j ,a ]} (s) P_L F(s)}_{L_s^2L_x^1} \norm{\chi_{[a, a+\delta_j]} (t) h(t) }_{L_t^2} , 
\end{align}
where $\chi_{\parallelogramm} = \chi_{\{a \leq t \leq a + \delta_j, \, t-2\delta_j \leq s \leq t - \delta_j \}}$ (see Figure~\ref{figure} for the parallelogram decomposition).

Now we need to sum over $a$, then by Cauchy-Schwarz, 
\begin{align}
& \quad \int_{\R} \int_{I_j} \frac{1}{\abs{t-s}^{\frac{d}{2}}} \norm{P_L F(s)}_{L_x^{1}} h(t) \, dsdt \\
& = \sum_a \int_{[a,a+\delta_j]} \int_{I_j} \frac{1}{\abs{t-s}^{\frac{d}{2}}} \norm{P_L F(s)}_{L_x^{1}} h(t) \, ds dt \\
& \lesssim \sum_a \delta_j^{-d/2} \delta_j \norm{\chi_{[a -2\delta_j ,a ]} (s) P_L F(s)}_{L_s^2L_x^1} \norm{\chi_{[a, a+\delta_j]} (t) h(t) }_{L_t^2}  \\
& \lesssim \delta_j^{-d/2} \delta_j \norm{\norm{\chi_{[a -2\delta_j ,a ]} (s) P_L F(s)}_{L_s^2L_x^1} }_{\ell_a^2} \norm{\norm{\chi_{[a, a+\delta_j]} (t) h(t) }_{L_t^2}}_{\ell_a^2} \\
& = \delta_j^{-d/2} \delta_j \norm{\chi_{[a -2\delta_j ,a ]} (s) P_L F(s)}_{L_s^2 \ell_a^2 L_x^1} \norm{\chi_{[a, a+\delta_j]} (t) h(t) }_{L_t^2 \ell_a^2} \\
& \lesssim \delta_j^{-d/2} \delta_j \norm{\chi_{[a -2\delta_j ,a ]} (s) P_L F(s)}_{L_s^2 \ell_a^1 L_x^1} \norm{\chi_{[a, a+\delta_j]} (t) h(t) }_{L_t^2 \ell_a^1} \\
& \lesssim \delta_j^{-d/2} \delta_j \norm{P_L F(s)}_{L_s^2 L_x^1} \norm{h(t)}_{L_t^2} \\
& \lesssim \delta_j^{-d/2} \delta_j \norm{P_L F(s)}_{L_s^2 L_x^1} \\
& = m^{d-2} 2^{-j \frac{d-2}{2}} \norm{P_L F}_{L_t^2 L_x^1} . \label{eq5}
\end{align}

\noindent{\bf The $j < 0$ bound, \eqref{eq BDD3}:}

When $j < 0$, recall we plan to show the following bound
\begin{align}
T_j & \lesssim C(N,L)^{1-\varepsilon} m^{\frac{d-2}{2} \varepsilon} 2^{\frac{j}{2}\varepsilon}  \norm{u_0}_{L^2} \norm{F}_{L_t^{2} L_x^{\frac{2d}{d+2},2}} . 
\end{align}

Note that for the unrestricted (in time) estimate, we have
\begin{align}
\norm{e^{it\Delta} P_N u_0 \parenthese{ \int_{\mathbb{R}} e^{i(t-s)\Delta} P_L F(s) \, ds }}_{L_{t,x}^2} & \lesssim C(N,L) \norm{u_0}_{L^2} \norm{\int_{\mathbb{R}} e^{-is\Delta}  F(s) \, ds  }_{L^2} \\
& \lesssim C(N,L) \norm{u_0}_{L^2} \norm{F}_{L_t^{2} L_x^{\frac{2d}{d+2},2}} , \label{eq2}
\end{align}
and by H\"older and Bernstein
\begin{align}
\norm{e^{it\Delta} P_N u_0 \parenthese{ \int_{\mathbb{R}} e^{i(t-s)\Delta} P_L F(s) \, ds }}_{L_t^{\infty} L_x^2} & \lesssim \min \{N,L\}^{\frac{d}{2}}  \norm{u_0}_{L^2} \norm{\int_{\mathbb{R}} e^{-is\Delta}  F(s) \, ds  }_{L^2} \\
& \lesssim \min \{N,L\}^{\frac{d}{2}} \norm{u_0}_{L^2} \norm{F}_{L_t^{2} L_x^{\frac{2d}{d+2},2}} .
\end{align}

Then interpolating the two unrestricted bounds and applying Christ-Kiselev gives
\begin{align}
\norm{e^{it\Delta} P_N u_0  \, e^{it\Delta}P_L \parenthese{\int_{I_j} e^{-is\Delta}  F(s) \, ds }}_{L_t^{\frac{2}{1-\varepsilon}} L_x^2} \lesssim C(N,L)^{1-\varepsilon} \min \{N,L\}^{\frac{d}{2} \varepsilon}  \norm{u_0}_{L^2} \norm{F}_{L_t^{2} L_x^{\frac{2d}{d+2},2}} .
\end{align}

Now
\begin{align}
& \quad \abs{\int_{[a,a+\delta_j]} \int e^{it\Delta} P_N u_0  \, e^{it\Delta}P_L \parenthese{\int_{I_j} e^{-is\Delta}  F(s) \, ds } \overline{h}(t,x) \, dxdt } \\
& \lesssim \norm{e^{it\Delta} P_N u_0  \, e^{it\Delta}P_L \parenthese{\int_{I_j} e^{-is\Delta} \chi_{[a-2\delta_j, a]} (s) F(s) \, ds }}_{L_t^{\frac{2}{1-\varepsilon}} L_x^2} \norm{\chi_{[a,a+\delta_j]} (t) h(t,x)}_{L_t^{\frac{2}{1+\varepsilon}}L_x^2}\\
& \lesssim C(N,L)^{1-\varepsilon} \min \{N,L\}^{\frac{d-2}{2} \varepsilon} 2^{\frac{j}{2} \varepsilon} \norm{u_0}_{L^2} \norm{\chi_{[a-2\delta_j ,a]}(s) F}_{L_s^{2} L_x^{\frac{2d}{d+2},2}} \norm{\chi_{[a,a+\delta_j]} (t) h}_{L_{t,x}^2} ,
\end{align}
which allows us to sum in $a \in \delta_j \Z$, 
\begin{align}
& \quad \sum_a \abs{ \int_{[a,a+\delta_j]} \int e^{it\Delta} P_N u_0  \, e^{it\Delta}P_L \parenthese{\int_{I_j} e^{-is\Delta}  F(s) \, ds } \overline{h}(t,x) \, dxdt} \\
& \lesssim C(N,L)^{1-\varepsilon} \min \{N,L\}^{\frac{d-2}{2} \varepsilon} 2^{\frac{j}{2}\varepsilon} \norm{u_0}_{L^2} \norm{ \norm{\chi_{[a-2\delta_j ,a]}(s) F}_{L_s^{2} L_x^{\frac{2d}{d+2},2}} }_{\ell_a^2} \norm{ \norm{\chi_{[a,a+\delta_j]} (t) h}_{L_{t,x}^2} }_{\ell_a^2} \\
& \lesssim C(N,L)^{1-\varepsilon} \min \{N,L\}^{\frac{d-2}{2} \varepsilon} 2^{\frac{j}{2} \varepsilon}  \norm{u_0}_{L^2} \norm{F}_{L_t^{2} L_x^{\frac{2d}{d+2},2}}  \norm{h}_{L_{t,x}^2} ,
\end{align}
hence \eqref{eq BDD3} follows.

Now we finish the proof of Lemma~\ref{lem CK}. 

\begin{figure}[!htbp]
\centering

\begin{tikzpicture}[scale=0.9, x=1cm, y=1cm]
\draw[->] (-5.5,-1) -- (5.5,-1) node[right] {$t$};
\draw[->] (-2.5,-8.5) -- (-2.5,5.5) node[above] {$s$};

\def\a{0}

\def\twojmm{0.5}   
\def\twojm{1}      
\def\twoj{2}       
\def\twojp{4}      
\def\twojpp{8}     


\coordinate (A) at (\a,\a-\twojp);
\coordinate (B) at (\a,\a-\twoj);
\coordinate (C) at (\a+\twoj,\a);
\coordinate (D) at (\a+\twoj,\a-\twoj);

\filldraw[
    fill=blue!15,
    draw=blue,
    thick
]
(A) -- (B) -- (C) -- (D) -- cycle;


\coordinate (L1) at (\a-\twoj,\a-3*\twoj);
\coordinate (L2) at (\a-\twoj,\a-2*\twoj);
\coordinate (L3) at (\a,\a-\twoj);
\coordinate (L4) at (\a,\a-2*\twoj);

\draw[dashed, thick]
(L1) -- (L2) -- (L3) -- (L4) -- cycle;


\coordinate (R1) at (\a+\twoj,\a-\twoj);
\coordinate (R2) at (\a+\twoj,\a);
\coordinate (R3) at (\a+2*\twoj,\a+\twoj);
\coordinate (R4) at (\a+2*\twoj,\a);

\draw[dashed, thick]
(R1) -- (R2) -- (R3) -- (R4) -- cycle;


\coordinate (A1) at (\a-2*\twojm,\a-4*\twojm);
\coordinate (A2) at (\a-2*\twojm,\a-3*\twojm);
\coordinate (A3) at (\a-\twojm,\a-2*\twojm);
\coordinate (A4) at (\a-\twojm,\a-3*\twojm);

\draw[dashed]
(A1) -- (A2) -- (A3) -- (A4) -- cycle;

\coordinate (B1) at (\a-\twojm,\a-3*\twojm);
\coordinate (B2) at (\a-\twojm,\a-2*\twojm);
\coordinate (B3) at (\a,\a-\twojm);
\coordinate (B4) at (\a,\a-2*\twojm);

\draw[dashed]
(B1) -- (B2) -- (B3) -- (B4) -- cycle;

\coordinate (C1) at (\a,\a-2*\twojm);
\coordinate (C2) at (\a,\a-\twojm);
\coordinate (C3) at (\a+\twojm,\a);
\coordinate (C4) at (\a+\twojm,\a-\twojm);

\draw[dashed]
(C1) -- (C2) -- (C3) -- (C4) -- cycle;

\coordinate (D1) at (\a+\twojm,\a-\twojm);
\coordinate (D2) at (\a+\twojm,\a);
\coordinate (D3) at (\a+2*\twojm,\a+\twojm);
\coordinate (D4) at (\a+2*\twojm,\a);

\draw[dashed]
(D1) -- (D2) -- (D3) -- (D4) -- cycle;

\coordinate (E1) at (\a+2*\twojm,\a);
\coordinate (E2) at (\a+2*\twojm,\a+\twojm);
\coordinate (E3) at (\a+3*\twojm,\a+2*\twojm);
\coordinate (E4) at (\a+3*\twojm,\a+\twojm);

\draw[dashed]
(E1) -- (E2) -- (E3) -- (E4) -- cycle;

\coordinate (F1) at (\a+3*\twojm,\a+\twojm);
\coordinate (F2) at (\a+3*\twojm,\a+2*\twojm);
\coordinate (F3) at (\a+4*\twojm,\a+3*\twojm);
\coordinate (F4) at (\a+4*\twojm,\a+2*\twojm);

\draw[dashed]
(F1) -- (F2) -- (F3) -- (F4) -- cycle;


\coordinate (AA1) at (\a-4*\twojmm,\a-6*\twojmm);
\coordinate (AA2) at (\a-4*\twojmm,\a-5*\twojmm);
\coordinate (AA3) at (\a-3*\twojmm,\a-4*\twojmm);
\coordinate (AA4) at (\a-3*\twojmm,\a-5*\twojmm);

\draw[dashed]
(AA1) -- (AA2) -- (AA3) -- (AA4) -- cycle;

\coordinate (BB1) at (\a-3*\twojmm,\a-5*\twojmm);
\coordinate (BB2) at (\a-3*\twojmm,\a-4*\twojmm);
\coordinate (BB3) at (\a-2*\twojmm,\a-3*\twojmm);
\coordinate (BB4) at (\a-2*\twojmm,\a-4*\twojmm);

\draw[dashed]
(BB1) -- (BB2) -- (BB3) -- (BB4) -- cycle;

\coordinate (CC1) at (\a-2*\twojmm,\a-4*\twojmm);
\coordinate (CC2) at (\a-2*\twojmm,\a-3*\twojmm);
\coordinate (CC3) at (\a-\twojmm,\a-2*\twojmm);
\coordinate (CC4) at (\a-\twojmm,\a-3*\twojmm);

\draw[dashed]
(CC1) -- (CC2) -- (CC3) -- (CC4) -- cycle;

\coordinate (DD1) at (\a-\twojmm,\a-3*\twojmm);
\coordinate (DD2) at (\a-\twojmm,\a-2*\twojmm);
\coordinate (DD3) at (\a,\a-\twojmm);
\coordinate (DD4) at (\a,\a-2*\twojmm);

\draw[dashed]
(DD1) -- (DD2) -- (DD3) -- (DD4) -- cycle;

\coordinate (EE1) at (\a,\a-2*\twojmm);
\coordinate (EE2) at (\a,\a-\twojmm);
\coordinate (EE3) at (\a+\twojmm,\a);
\coordinate (EE4) at (\a+\twojmm,\a-\twojmm);

\draw[dashed]
(EE1) -- (EE2) -- (EE3) -- (EE4) -- cycle;

\coordinate (FF1) at (\a+\twojmm,\a-\twojmm);
\coordinate (FF2) at (\a+\twojmm,\a);
\coordinate (FF3) at (\a+2*\twojmm,\a+\twojmm);
\coordinate (FF4) at (\a+2*\twojmm,\a);

\draw[dashed]
(FF1) -- (FF2) -- (FF3) -- (FF4) -- cycle;

\coordinate (GG1) at (\a+2*\twojmm,\a);
\coordinate (GG2) at (\a+2*\twojmm,\a+\twojmm);
\coordinate (GG3) at (\a+3*\twojmm,\a+2*\twojmm);
\coordinate (GG4) at (\a+3*\twojmm,\a+\twojmm);

\draw[dashed]
(GG1) -- (GG2) -- (GG3) -- (GG4) -- cycle;

\coordinate (HH1) at (\a+3*\twojmm,\a+\twojmm);
\coordinate (HH2) at (\a+3*\twojmm,\a+2*\twojmm);
\coordinate (HH3) at (\a+4*\twojmm,\a+3*\twojmm);
\coordinate (HH4) at (\a+4*\twojmm,\a+2*\twojmm);

\draw[dashed]
(HH1) -- (HH2) -- (HH3) -- (HH4) -- cycle;

\coordinate (II1) at (\a+4*\twojmm,\a+2*\twojmm);
\coordinate (II2) at (\a+4*\twojmm,\a+3*\twojmm);
\coordinate (II3) at (\a+5*\twojmm,\a+4*\twojmm);
\coordinate (II4) at (\a+5*\twojmm,\a+3*\twojmm);

\draw[dashed]
(II1) -- (II2) -- (II3) -- (II4) -- cycle;

\coordinate (JJ1) at (\a+5*\twojmm,\a+3*\twojmm);
\coordinate (JJ2) at (\a+5*\twojmm,\a+4*\twojmm);
\coordinate (JJ3) at (\a+6*\twojmm,\a+5*\twojmm);
\coordinate (JJ4) at (\a+6*\twojmm,\a+4*\twojmm);

\draw[dashed]
(JJ1) -- (JJ2) -- (JJ3) -- (JJ4) -- cycle;

\coordinate (KK1) at (\a+6*\twojmm,\a+4*\twojmm);
\coordinate (KK2) at (\a+6*\twojmm,\a+5*\twojmm);
\coordinate (KK3) at (\a+7*\twojmm,\a+6*\twojmm);
\coordinate (KK4) at (\a+7*\twojmm,\a+5*\twojmm);

\draw[dashed]
(KK1) -- (KK2) -- (KK3) -- (KK4) -- cycle;

\coordinate (LL1) at (\a+7*\twojmm,\a+5*\twojmm);
\coordinate (LL2) at (\a+7*\twojmm,\a+6*\twojmm);
\coordinate (LL3) at (\a+8*\twojmm,\a+7*\twojmm);
\coordinate (LL4) at (\a+8*\twojmm,\a+6*\twojmm);

\draw[dashed]
(LL1) -- (LL2) -- (LL3) -- (LL4) -- cycle;

\draw[dashed]
(-5.5,-6) -- (5.5,5)
node[above right] {$t=s+{2^{-2}\delta_j}$};

\draw[dashed]
(-5.5,-6.5) -- (5.5,4.5)
node[above right] {$t=s+{2^{-1}\delta_j}$};

\draw[dashed]
(-5.5,-7.5) -- (5.5,3.5)
node[above right] {$t=s+{\delta_j}$};

\draw[dashed]
(-5.5,-9.5) -- (5.5,1.5)
node[above right] {$t=s+{2\delta_j}$};



\node[below] at (A)
{$(a,\; a-{2\delta_j})$};

\node[below] at (B)
{$(a,\; a-{\delta_j})$};

\node[right] at (C)
{$(a+{\delta_j},\; a)$};

\node[right] at (D)
{$(a+{\delta_j},\; a-{\delta_j})$};

\node at (3.5,-3)
{\small
$\chi_{\{
a \le t \le a+{\delta_j},\;
t-{2\delta_j}\le s\le t-{\delta_j}
\}}$
};

\end{tikzpicture}

\caption{Parallelogram decomposition {with $\delta_j=m^{-2}2^j$ and $m=\min\{N,L\}$}.}

\label{figure}
\end{figure}
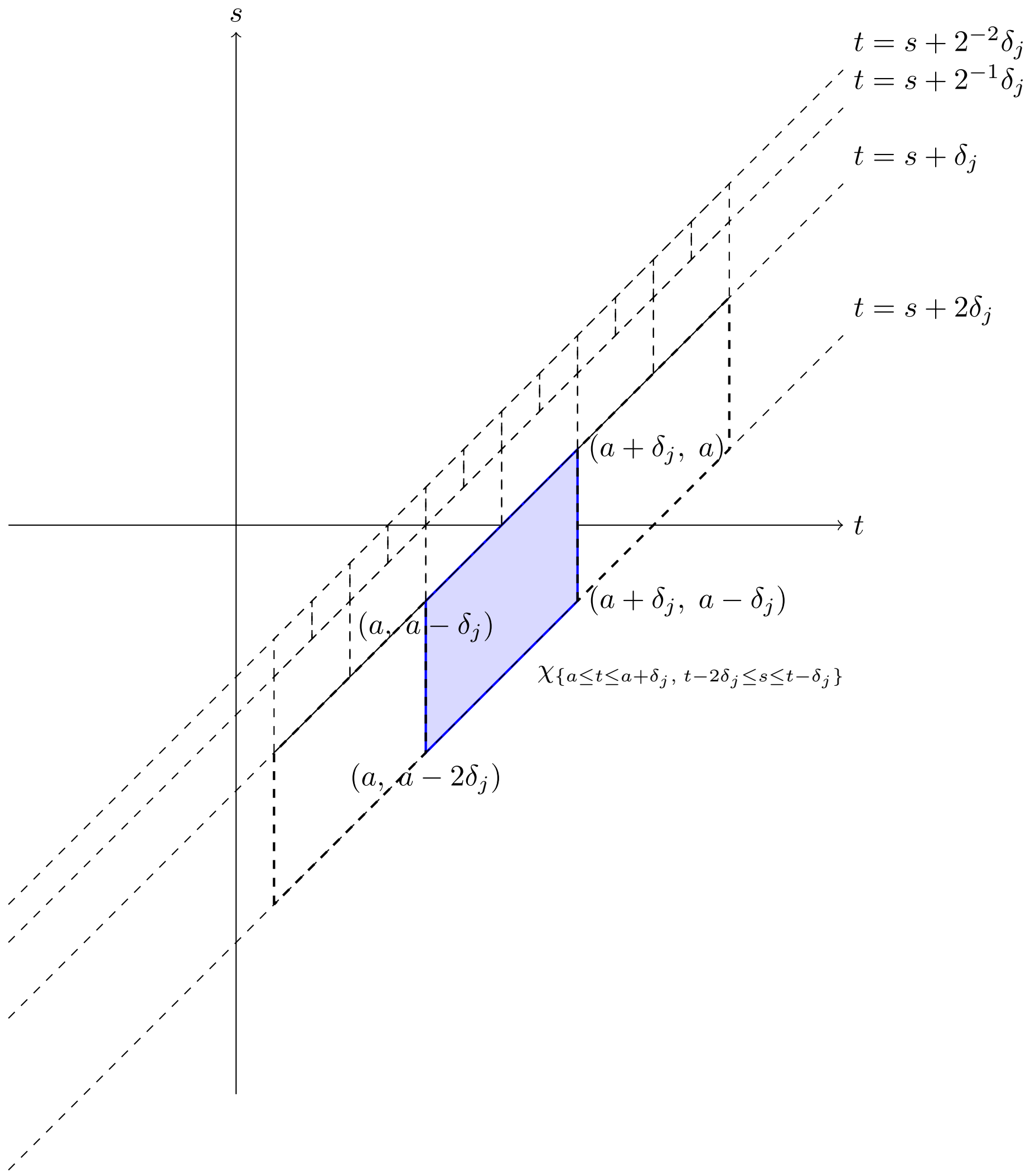

\end{proof}

\begin{lemma}[Endpoint of bilinear Christ-Kiselev Lemma, Duhamel-Duhamel]\label{lem CK'}
Let $d\geq 3$. 
Assume the bilinear estimate in Lemma~\ref{lem CK}. Then we have
\begin{align}\label{eq CK3}
& \quad \norm{\parenthese{ \int_{-\infty}^t e^{i(t-s)\Delta} P_N F (s) \, ds } \parenthese{ \int_{-\infty}^t e^{i(t-s')\Delta} P_L G(s') \, ds' }}_{L_{t,x}^2} \\
& \lesssim C(N,L)^{1-\varepsilon} \min \{N,L\}^{\frac{d-2}{2} \varepsilon} \norm{F}_{L_t^{2} L_x^{\frac{2d}{d+2},2}} \norm{G}_{L_t^{2} L_x^{\frac{2d}{d+2},2}} ,
\end{align}
and
\begin{align}\label{eq CK4}
& \quad \norm{\parenthese{ \int_{-\infty}^t e^{i(t-s)\Delta} P_N F (s) \, ds } \parenthese{ \int_{-\infty}^t e^{i(t-s')\Delta} P_L G(s') \, ds' }}_{L_{t,x}^2}  \\
& \lesssim C(N,L)^{1-\varepsilon} \min \{ N,L\}^{\frac{d-2}{2} \varepsilon} \norm{F}_{L_t^{2} L_x^{\frac{2d}{d+2}}} \norm{G}_{L_t^{2} L_x^{\frac{2d}{d+2}}} .
\end{align}
\end{lemma}

\begin{remark}
Again, we remark that the upper bounds in Lemma~\ref{lem CK'} are also given in the form of
\begin{align}
C(N,L)^{1-\varepsilon} \min \{N,L\}^{\frac{d-2}{2} \varepsilon} = \min\{N,L\}^{\frac{d-2}{2}} \parenthese{\frac{\min \{N,L\}}{\max \{ N,L\}} }^{\frac{1 - \varepsilon}{2}} .
\end{align}   
\end{remark}

\begin{proof}[Proof of Lemma~\ref{lem CK'}]
We only show the bound, \eqref{eq CK3}. The second bound, \eqref{eq CK4}, follows similarly from the same proof. 
In the rest of the proof, we choose 
\begin{align}
\rho = \frac{\varepsilon}{2} , \quad \theta = \frac{\varepsilon}{2-\varepsilon} \quad, \gamma = \frac{(d-2)\theta}{4(1-\theta)}, \quad \frac{1}{q} = \frac{1}{2} + \gamma, \quad \frac{1}{r} = \frac{d+2}{2d} - \frac{2\gamma}{d}
\end{align}
for small enough $\varepsilon >0$. 
Here $\rho + (1-\rho) \theta = \varepsilon$ and $(q,r)$ is a dual Strichartz pair. We will write instead $q =2-$ and $r = \frac{2d}{d+2}+$ for simplicity. Also we treat all $\rho , \theta, \gamma , \varepsilon$ as very small parameters.

By Lemma~\ref{lem CK} and Strichartz estimates, we see that
\begin{align}
& \quad \norm{\parenthese{ \int_{-\infty}^{\infty} e^{i(t-s)\Delta} P_N F (s) \, ds } \parenthese{ \int_{-\infty}^t e^{i(t-s')\Delta} P_L G(s') \, ds' }}_{L_{t,x}^2} \\
& =\norm{ e^{it\Delta} P_N \underbrace{\parenthese{\int_{-\infty}^{\infty} e^{-is\Delta} F (s) \, ds }}_{u_0 \text{ in Lemma~\ref{lem CK}}} \parenthese{ \int_{-\infty}^t e^{i(t-s')\Delta} P_L G(s') \, ds' }}_{L_{t,x}^2} \\
& \lesssim C(N,L)^{1-\rho} \min \{N,L\}^{\frac{d-2}{2} \rho} \norm{\int e^{-is\Delta} F(s) \, ds}_{L_x^{2}} \norm{G}_{L_t^{2} L_x^{\frac{2d}{d+2},2}} \\
& \lesssim C(N,L)^{1-\rho} \min \{N,L\}^{\frac{d-2}{2} \rho} \norm{F}_{L_t^{2-} L_x^{\frac{2d}{d+2}+}} \norm{G}_{L_t^{2} L_x^{\frac{2d}{d+2},2}} . \label{eq6}
\end{align}

Then by Christ-Kiselev Lemma, we have 
\begin{align}
& \quad \norm{\parenthese{ \int_{-\infty}^{t} e^{i(t-s)\Delta} P_N F (s) \, ds } \parenthese{ \int_{-\infty}^t e^{i(t-s')\Delta} P_L G(s') \, ds' }}_{L_{t,x}^2} \\
& \lesssim C(N,L)^{1-\rho} \min \{N,L\}^{\frac{d-2}{2} \rho} \norm{F}_{L_t^{2-} L_x^{\frac{2d}{d+2}+}}  \norm{G}_{L_t^{2} L_x^{\frac{2d}{d+2},2}} . \label{eq3}
\end{align}

Like we did in Lemma~\ref{lem CK}, we denote $m = \min \{N,L\}$. Then we break the interval $(-\infty , t)$ into $\cup_{j=-\infty}^{\infty} I_j$ with $\delta_j = m^{-2}2^j$ and $I_j = [t- 2\delta_j , t - \delta_j]$, and write
\begin{align}
\eqref{eq CK3} & = \norm{\parenthese{ \int_{-\infty}^t e^{i(t-s)\Delta} P_N F (s) \, ds } \parenthese{ \int_{-\infty}^t e^{i(t-s')\Delta} P_L G(s') \, ds' }}_{L_{t,x}^2} \\
& \leq \sum_{j} \underbrace{\norm{ \parenthese{\int_{I_j} e^{i(t-s)\Delta} P_N F(s) \,ds} \,  \parenthese{\int_{-\infty}^t e^{i(t-s')\Delta} P_L G(s') \, ds' }}_{L_{t,x}^2} }_{T_{j}} . 
\end{align}

Similarly, we estimate $T_j$ in different regimes: $j < 0$ and $j \geq 0$. 

We observe that for fixed $j$ and $t \in [a, a+ \delta_j]$,
\begin{align}
\int_{I_j} e^{i(t-s)\Delta} P_N F(s) \, ds = \int_{I_j} e^{i(t-s)\Delta} P_N \chi_{[a-2\delta_j ,a]} (s) F(s) \, ds ,
\end{align}
and we are able to sum over $\{ [a-2\delta_j ,a]\}_a$ in $\ell_a^2$ sense.

When $j \geq 0$, by \eqref{eq3} and H\"older in time, we have
\begin{align}
& \quad \norm{\chi_{[a,a+ \delta_j]} (t) \parenthese{\int_{I_j} e^{i(t-s)\Delta} P_N F (s) \, ds}. \parenthese{\int_{-\infty}^t e^{i(t-s')\Delta} P_L G (s') \, ds'} }_{L_{t,x}^2} \\
& \lesssim C(N,L)^{1-\rho} m^{\frac{d-2}{2} \rho} \norm{\chi_{[a-2\delta_j ,a]}F}_{L_t^{2-} L_x^{\frac{2d}{d+2}+}} \norm{G}_{L_t^2 L_x^{\frac{2d}{d+2},2}} \\
& \lesssim C(N,L)^{1-\rho} m^{\frac{d-2}{2} \rho} \delta_j^{\gamma} \norm{\chi_{[a-2\delta_j ,a]}F}_{L_t^{2} L_x^{\frac{2d}{d+2}+}} \norm{G}_{L_t^2 L_x^{\frac{2d}{d+2},2}} .
\end{align}
Again, we used
\begin{align}
\int_{I_j} = \int_{s < t-\delta_j} - \int_{s< t-2\delta_j}
\end{align}
and then applied Christ-Kiselev.

Then summing over $a \in \delta_j \Z$ gives
\begin{align}
T_j \lesssim C(N,L)^{1-\rho} m^{\frac{d-2}{2} \rho} \delta_j^{\gamma} \norm{F}_{L_t^{2} L_x^{\frac{2d}{d+2}+}}  \norm{G}_{L_t^{2} L_x^{\frac{2d}{d+2},2}} . \label{eq7}
\end{align}
On the other hand, by H\"older, the dispersive estimate and similar computation from \eqref{eq4} to \eqref{eq5} (where we handled the bound \eqref{eq BDD2} in Lemma~\ref{lem CK}), we have
\begin{align}
T_j & \leq \norm{\int_{I_j} \norm{e^{i(t-s)\Delta} P_N F (s)}_{L_x^{\infty}}  \norm{\int_{-\infty}^t e^{i(t-s')\Delta} P_L G(s') \, ds'}_{L_x^{2}} \, ds }_{L_{t}^2} \\
& \lesssim \norm{\int_{I_j}   \frac{1}{\abs{t-s}^{\frac{d}{2}}} \norm{P_N F(s)}_{L_x^{1}} \, ds }_{L_{t}^2} \norm{\int_{-\infty}^t e^{i(t-s')\Delta} P_L G(s') \, ds' }_{L_t^{\infty} L_x^2} \\
& \lesssim  \delta_j^{-d/2} \delta_j \norm{F(s)}_{L_s^2 L_x^1} \norm{G}_{L_t^{2} L_x^{\frac{2d}{d+2},2}} \\
& = m^{d-2} 2^{-j \frac{d-2}{2}} \norm{F}_{L_t^2 L_x^1} \norm{G}_{L_t^2 L_x^{\frac{2d}{d+2} ,2}} . \label{eq8}
\end{align}
Now, interpolation between \eqref{eq7} and \eqref{eq8} (with weight $\theta$ on \eqref{eq8}), we have 
\begin{align}
T_j \lesssim C(N,L)^{1-\varepsilon} m^{\frac{d-2}{2} \varepsilon} 2^{-j \frac{d-2}{4}\theta} \norm{F}_{L_t^{2} L_x^{\frac{2d}{d+2},2}}  \norm{G}_{L_t^{2} L_x^{\frac{2d}{d+2},2}} . 
\end{align}

When $j < 0$, we first need similar unrestricted estimates (as what we did in the proof of \eqref{eq BDD3} in Lemma~\ref{lem CK}) before applying Christ-Kiselev. In fact, by Lemma~\ref{lem CK} and Strichartz, we have
\begin{align}
& \quad \norm{\parenthese{\int_{-\infty}^{\infty} e^{i(t-s)\Delta} P_N F(s) \, ds} \parenthese{\int_{-\infty}^t e^{i(t-s')\Delta} P_L G(s') \, ds'}}_{L_{t,x}^2} \\
& \lesssim C(N,L)^{1-\rho} m^{\frac{d-2}{2} \rho} \norm{F}_{L_t^2 L_x^{\frac{2d}{d+2},2}} \norm{G}_{L_t^2 L_x^{\frac{2d}{d+2},2}} ,
\end{align}
and
\begin{align}
& \quad \norm{\parenthese{\int_{-\infty}^{\infty} e^{i(t-s)\Delta} P_N F(s) \, ds} \parenthese{\int_{-\infty}^t e^{i(t-s')\Delta} P_L G(s') \, ds'}}_{L_t^{\infty} L_x^2} \\
& \lesssim m^{\frac{d}{2}} \norm{\int_{-\infty}^{\infty} e^{-is\Delta} F(s) \, ds}_{L_t^{\infty} L_x^2} \norm{\int_{-\infty}^t e^{i(t-s')\Delta} P_L G(s') \, ds}_{L_t^{\infty} L_x^2} \\
& \lesssim m^{\frac{d}{2}} \norm{F}_{L_t^2 L_x^{\frac{2d}{d+2},2}} \norm{G}_{L_t^2 L_x^{\frac{2d}{d+2},2}} .
\end{align}
Then interpolating the two unrestricted estimates above, we have
\begin{align}
& \quad \norm{\parenthese{\int_{-\infty}^{\infty} e^{i(t-s)\Delta} P_N F(s) \, ds} \parenthese{\int_{-\infty}^t e^{i(t-s')\Delta} P_L G(s') \, ds'}}_{L_t^{\frac{2}{1-\theta}} L_x^2} \\
& \lesssim C(N,L)^{1-\varepsilon} m^{\frac{d-2}{2} \varepsilon + \theta} \norm{F}_{L_t^2 L_x^{\frac{2d}{d+2},2}} \norm{G}_{L_t^2 L_x^{\frac{2d}{d+2},2}} .
\end{align}

Similar treatment as what we did in the bound \eqref{eq BDD3} in Lemma~\ref{lem CK}, we have for $h$ with $\norm{h}_{L_{t,x}^2}=1$ 
\begin{align}
T_j & \leq \sup_h \sum_{a \in \delta_j \Z} \abs{ \int_{[a, a+ \delta_j]} \int_{\R^d} \parenthese{\int_{I_j} e^{i(t-s)\Delta} P_N F(s) \, ds } \parenthese{\int_{-\infty}^t e^{i(t-s')\Delta} P_L G(s') \, ds'} \overline{h}(t,x) \, dxdt} \\
& \lesssim \sup_h \sum_a \norm{\parenthese{\int_{I_j} e^{i(t-s)\Delta}  P_N \chi_{[a-2\delta_j, a]} (s) F(s) \, ds } \parenthese{\int_{-\infty}^t e^{i(t-s')\Delta} P_L G(s') \, ds'}}_{L_t^{\frac{2}{1-\theta}} L_x^2} \\
& \hspace{1cm} \times \norm{\chi_{[a,a+\delta_j]} (t) h(t,x)}_{L_t^{\frac{2}{1+\theta}}L_x^2}\\
& \lesssim C(N,L)^{1-\varepsilon} m^{\frac{d-2}{2} \varepsilon + \theta} \delta_j^{\theta/2} \sup_h \sum_a \norm{ \chi_{[a - 2\delta_j ,a]} F}_{L_t^2 L_x^{\frac{2d}{d+2},2}} \norm{G}_{L_t^{2} L_x^{\frac{2d}{d+2},2}} \norm{\chi_{[a,a+\delta_j]} h}_{L_{t,x}^2} \\
& \lesssim C(N,L)^{1-\varepsilon} m^{\frac{d-2}{2} \varepsilon} 2^{j \theta/2} \norm{G}_{L_t^2 L_x^{\frac{2d}{d+2},2}} \sup_h \norm{\norm{\chi_{[a-2\delta_j,a]} F}_{L_t^2 L_x^{\frac{2d}{d+2},2}}}_{\ell_{a}^2} \norm{\norm{\chi_{[a,a+\delta_j]} h}_{L_{t,x}^2}}_{\ell_{a}^2}\\
& \lesssim C(N,L)^{1-\varepsilon} m^{\frac{d-2}{2} \varepsilon}2^{j \theta/2} \norm{F}_{L_t^2 L_x^{\frac{2d}{d+2},2}} \norm{G}_{L_t^{2} L_x^{\frac{2d}{d+2},2}}  \norm{h}_{L_{t,x}^2} .
\end{align}

Now both $j < 0$ and $j \geq 0$ bounds allow us to sum $T_j$ and get the desired bound \eqref{eq CK3} in this lemma.

We then finish the proof of Lemma~\ref{lem CK'}. 

\end{proof}

\begin{theorem}\label{thm weightedCK}
Let $0<\varepsilon< \frac{1}{2}$. Then
\begin{align}
\norm{e^{it\Delta}P_N f\, P_L\int_{-\infty}^{t}e^{i(t-s)\Delta}F(s)\,ds}_{L^2_{t,x}} & \lesssim \frac{\min\{N,L\}^{\frac{d-1}{2}-\varepsilon}} {\max\{N,L\}^{\frac{1}{2}-\varepsilon}\inner{L}^{\frac{1}{2}}} \norm{P_Nf}_{L^2_x}\norm{\inner{x}F}_{L^2_{t,x}} . \label{eq weightedLD}
\end{align}
Moreover, 
\begin{align}
& \quad \norm{P_N\int_{-\infty}^{t}e^{i(t-s)\Delta}F(s)\,ds\, P_N\int_{-\infty}^{t}e^{i(t-s')\Delta}G(s')\,ds'}_{L^2_{t,x}} \\
& \lesssim \frac{\min\{N,L\}^{\frac{d-1}{2}-\varepsilon}} {\max\{N,L\}^{\frac{1}{2}-\varepsilon}\inner{N}^{\frac{1}{2}}\inner{L}^{\frac{1}{2}}}  \norm{\inner{x}F}_{L^2_{t,x}}\norm{\inner{x}G}_{L^2_{t,x}} . \label{eq weightedDD}
\end{align}
The same estimates hold with all the Duhamel integrals replaced by $\int_0^t$. 
\end{theorem}

Note that the frequency factor has precisely the same form as in Lemma~\ref{lem CK}:
\begin{align}
\frac{\min\{N,L\}^{\frac{d-1}{2}-\varepsilon}}  {\max\{N,L\}^{\frac{1}{2}-\varepsilon}}  =C(N,L)^{1-2\varepsilon}\min\{N,L\}^{(d-2)\varepsilon}.
\end{align}
Thus the loss is an arbitrarily small power of the frequency-separation
ratio. The additional factor $\inner{L}^{-\frac{1}{2}}$ in
\eqref{eq weightedLD}, and the factor $\inner{N}^{-\frac{1}{2}}\inner{L}^{-\frac{1}{2}}$ in
\eqref{eq weightedDD}, come from dual local smoothing.

\begin{proof}[Proof of Theorem~\ref{thm weightedCK}]
For every $M>0$, by duality of local smoothing, we have
\begin{align}
\norm{P_M\int_\R e^{-is\Delta}F(s)\,ds}_{L^2_x} \lesssim \inner{M}^{-1/2}\norm{\inner{x}F}_{L^2_{t,x}}. \label{eq dualweighted}
\end{align}

By the linear bilinear Strichartz estimate, we have the following unrestricted (in time) estimates
\begin{align}
\norm{e^{it\Delta}P_Nf\,P_L\int_\R e^{i(t-s)\Delta}F(s)\,ds}_{L^2_{t,x}} \lesssim C(N,L)\norm{P_Nf}_{L^2_x} \norm{P_L\int_\R e^{-is\Delta}F(s)\,ds}_{L^2_x}, \label{eq fullbilinear} \\
\norm{e^{it\Delta}P_Nf\,P_L\int_\R e^{i(t-s)\Delta}F(s)\,ds}_{L^\infty_{t,x}} \lesssim (NL)^{d/2}\norm{P_Nf}_{L^2_x} \norm{P_L\int_\R e^{-is\Delta}F(s)\,ds}_{L^2_x}. \label{eq fullbernstein}
\end{align}
Interpolating \eqref{eq fullbilinear} and \eqref{eq fullbernstein}, with weight $\varepsilon_0$ on \eqref{eq fullbernstein}, and then using \eqref{eq dualweighted}, we obtain
\begin{align}
\norm{e^{it\Delta}P_N f\, P_L \int_\R e^{i(t-s)\Delta} F(s)\,ds}_{L^{\frac{2}{1-\varepsilon_0}}_{t,x}} \lesssim C(N,L)^{1-\varepsilon_0}(NL)^{\frac{d\varepsilon_0}{2}}  \inner{L}^{-\frac{1}{2}}\norm{P_N }_{L^2_x}  \norm{\inner{x} F}_{L^2_{t,x}}. 
\end{align}
Then Christ-Kiselev lemma gives a restriced in time estimate
\begin{align}
\norm{e^{it\Delta}P_Nf\, P_L\int_{-\infty}^{t}e^{i(t-s)\Delta}F(s)\,ds}_{L^{\frac{2}{1-\varepsilon_0}}_{t,x}} \lesssim C(N,L)^{1-\varepsilon_0}(NL)^{\frac{d\varepsilon_0}{2}} \inner{L}^{-\frac{1}{2}}\norm{P_N f}_{L^2_x} \norm{\inner{x} F}_{L^2_{t,x}}. \label{eq retardedhigh}
\end{align}

To get another restricted in time estimate, we use Strichartz estimates again and write
\begin{align}
\norm{P_L\int_\R e^{i(t-s)\Delta}F(s)\,ds}_{L_{t,x}^{\frac{2(d+2)}{d}}} &\lesssim \norm{P_L\int_\R e^{-is\Delta}F(s)\,ds}_{L^2_x} \lesssim\inner{L}^{-\frac{1}{2}}\norm{\inner{x}F}_{L^2_{t,x}}.
\end{align}
then Christ-Kiselev gives
\begin{align}
\norm{P_L\int_{-\infty}^{t}e^{i(t-s)\Delta}F(s)\,ds}_{L_{t,x}^{\frac{2(d+2)}{d}}} \lesssim\inner{L}^{-\frac{1}{2}}\norm{\inner{x}F}_{L^2_{t,x}}. \label{eq linearretarded}
\end{align}

Now, H\"older and Strichartz allow the second resticted in time estimate
\begin{align}
\norm{e^{it\Delta}P_Nf\, P_L\int_{-\infty}^{t}e^{i(t-s)\Delta}F(s)\,ds}_{L^{\frac{d+2}{d}}_{t,x}} & \lesssim \norm{e^{it\Delta}P_Nf}_{L_{t,x}^{\frac{2(d+2)}{d}}} \norm{P_L\int_{-\infty}^{t}e^{i(t-s)\Delta}F(s)\,ds}_{L_{t,x}^{\frac{2(d+2)}{d}}}\\
& \lesssim\inner{L}^{-\frac{1}{2}}\norm{P_Nf}_{L^2_x} \norm{\inner{x}F}_{L^2_{t,x}}. \label{eq retardedlow}
\end{align}

Then interpolating \eqref{eq retardedhigh} and \eqref{eq retardedlow}, we obtain
\begin{align}
\norm{e^{it\Delta}P_Nf\, P_L\int_{-\infty}^{t}e^{i(t-s)\Delta}F(s)\,ds}_{L^2_{t,x}} & \lesssim_{\varepsilon_0}  \parenthese{C(N,L)^{1-\varepsilon_0}(NL)^{\frac{d\varepsilon_0}{2}}}^\theta \inner{L}^{-\frac{1}{2}}\norm{P_Nf}_{L^2_x}\norm{\inner{x}F}_{L^2_{t,x}} \\
& =  \frac{\min\{N,L\}^{\frac{d-1}{2}-\varepsilon} }{\max\{N,L\}^{\frac{1}{2}-\varepsilon} \inner{L}^{\frac{1}{2}}} \norm{P_Nf}_{L^2_x}\norm{\inner{x}F}_{L^2_{t,x}} 
\end{align}
where
\begin{align}
\varepsilon=\frac{d^2\varepsilon_0}{2\parenthese{d-2+(d+2)\varepsilon_0}} 
\end{align}
and the interpolation weight $\theta$ 
\begin{align}\label{eq thetadelta}
\theta=\frac{d-2}{d-2+(d+2)\varepsilon_0}
\end{align}
on \eqref{eq retardedhigh}. 
This now proves \eqref{eq weightedLD}.

The homogeneous bilinear estimate and Bernstein, interpolated exactly as above, give
\begin{align}
& \quad \norm{P_N\int_\R e^{i(t-s)\Delta}F(s)\,ds\, P_L\int_\R e^{i(t-s')\Delta}G(s')\,ds'}_{L^{\frac{2}{1-\varepsilon_0}}_{t,x}}\\
&\lesssim C(N,L)^{1-\varepsilon_0}(NL)^{\frac{d\varepsilon_0}{2}} \norm{P_N\int_\R e^{-is\Delta}F(s)\,ds}_{L^2_x} \norm{P_L\int_\R e^{-is'\Delta}G(s')\,ds'}_{L^2_x}\\
& \lesssim C(N,L)^{1-\varepsilon_0}(NL)^{\frac{d\varepsilon_0}{2}} \inner{N}^{-\frac{1}{2}}\inner{L}^{-\frac{1}{2}} \norm{\inner{x}F}_{L^2_{t,x}}\norm{\inner{x}G}_{L^2_{t,x}}. \label{eq DDfullhigh}
\end{align}
For fixed $G$, the left-hand side of \eqref{eq DDfullhigh} is the norm of a linear integral operator in $F$, then Christ-Kiselev truncates the $F$ integral and gives
\begin{align}
&\quad \norm{P_N\int_{-\infty}^{t}e^{i(t-s)\Delta}F(s)\,ds\, P_L\int_\R e^{i(t-s')\Delta}G(s')\,ds'}_{L^{\frac{2}{1-\varepsilon_0}}_{t,x}}\\
&\lesssim_{\varepsilon_0} C(N,L)^{1-\varepsilon_0}(NL)^{\frac{d\varepsilon_0}{2}} \inner{N}^{-\frac{1}{2}}\inner{L}^{-\frac{1}{2}} \norm{\inner{x}F}_{L^2_{t,x}}\norm{\inner{x}G}_{L^2_{t,x}}. 
\end{align}

Applying Christ-Kiselev gives
\begin{align}
&\quad \norm{P_N\int_{-\infty}^{t}e^{i(t-s)\Delta}F(s)\,ds\, P_L\int_{-\infty}^{t}e^{i(t-s')\Delta}G(s')\,ds'}_{L^{\frac{2}{1-\varepsilon_0}}_{t,x}}\\
&\lesssim C(N,L)^{1-\varepsilon_0}(NL)^{\frac{d\varepsilon_0}{2}} \inner{N}^{-\frac{1}{2}}\inner{L}^{-\frac{1}{2}} \norm{\inner{x}F}_{L^2_{t,x}}\norm{\inner{x}G}_{L^2_{t,x}}.\label{eq DDretardedhigh}
\end{align}

Applying \eqref{eq linearretarded} to both factors and using H\"older,
we have
\begin{align}
& \quad \norm{ P_N\int_{-\infty}^{t}e^{i(t-s)\Delta}F(s)\,ds\,P_L\int_{-\infty}^{t}e^{i(t-s')\Delta}G(s')\,ds'}_{L^{\frac{d+2}{d}}_{t,x}}\\
& \lesssim \inner{N}^{-\frac{1}{2}}\inner{L}^{-\frac{1}{2}} \norm{\inner{x}F}_{L^2_{t,x}}\norm{\inner{x}G}_{L^2_{t,x}}. \label{eq DDretardedlow}
\end{align}
Interpolating \eqref{eq DDretardedhigh} and \eqref{eq DDretardedlow}, with weight $\theta$ from \eqref{eq thetadelta} on the former, proves \eqref{eq weightedDD}. 

Now we finish the proof of Theorem~\ref{thm weightedCK}. 
\end{proof}


\section{Bilinear Frequency-Interaction Strichartz estimates}\label{sec PBilinear}

In this section, we generalize bilinear estimates for the linear Schr\"odinger equation with potentials in the setting in which the Littlewood-Paley projectors are applied after the Schr\"odinger flow operator. We consider the cases of when potential has a singularity and does not separately.

\PBilinear*

\begin{proof}[Proof of Theorem~\ref{thm PBilinear}]
First, we write the solution via the Duhamel principle, that is,
\begin{align}
u & = e^{it (\Delta -V)} u_0 = e^{it\Delta} u_0  -i \int_{0}^t e^{i(t-s)\Delta} V u (s) \, ds ,\\
v & = e^{it (\Delta -V)} v_0 = e^{it\Delta} v_0  -i \int_{0}^t e^{i(t-s)\Delta} V v (s) \, ds .
\end{align}

Next, we apply the projections to both sides of the equations above and expand the product into the following terms
\begin{align}
\norm{P_N u P_L v}_{L_{t}^pL_{x}^q} & \lesssim \norm{ e^{it\Delta} P_N u_0 e^{it\Delta} P_L v_0}_{L_{t}^pL_{x}^q} \\
& \quad + \norm{e^{it\Delta} P_N u_0 \parenthese{ \int_{0}^t e^{i(t-s)\Delta} P_L (V v)(s) \, ds }}_{L_{t}^pL_{x}^q} \\
& \quad + \norm{e^{it\Delta} P_L v_0 \parenthese{\int_{0}^t e^{i(t-s)\Delta} P_N (Vu)(s) \, ds }}_{L_{t}^pL_{x}^q}  \\
& \quad+ \norm{\parenthese{\int_{0}^t e^{i(t-s)\Delta} P_N (Vu)(s) \, ds} \parenthese{\int_{0}^t e^{i(t-s')\Delta} P_L (Vv)(s') \, ds'}}_{L_{t}^pL_{x}^q} \\
& = : I + II + III + IV . 
\end{align}

Next, we will treat the four terms separately. 
By assumption, we have
\begin{align}
I \leq C_{p, q}(N,L) \norm{u_0}_{L^2} \norm{v_0}_{L^2} .
\end{align}


For the remainder of the estimates, we will derive inhomogeneous estimates from homogeneous estimates and then apply interpolation.

First note that, using the dual to Corollary~\ref{cor BPSTZ resolvent}, we get
\begin{align}
\norm{\int_{\R} e^{-is \Delta} (Vu) (s) \, ds}_{L^2_x}=\norm{\int_{\R} e^{-is \Delta}|x|^{-1} |x|(Vu) (s) \, ds}_{L_x^{2}} \lesssim \norm{|x| Vu}_{L_t^2 L_x^{2}} .
\end{align}

Therefore, for $\frac{1}{p}+\frac{1}{q}\leq 1$,
    \begin{align}
        \norm{e^{it\Delta} P_N u_0 \parenthese{ \int_{\mathbb{R}} e^{i(t-s)\Delta} P_L (V v)(s) \, ds }}_{L_{t}^pL^q_x}
        & \leq C_{p, q}(N, L) \norm{ u_0}_{L_{x}^{2}} \norm{\int_{\R} e^{-is \Delta} P_L(Vv) (s) \, ds}_{L^{2}_{x}}\\ 
        & \leq C_{p, q}(N, L) \norm{ u_0}_{L_{x}^{2}} \norm{\int_{\R} e^{-is \Delta} V v(s) \, ds}_{L^{2}_{x}}\\
        & \lesssim C_{p, q}(N, L)\norm{u_0}_{L^2} \norm{| x | Vv}_{L^2_{t,x}} .
    \end{align}
Therefore, for $p>2$,  Christ-Kiselev Lemma (Lemma~\ref{lem:christkiselev}) with 
    \begin{align*}
        K(t, s)f(s):=  e^{it\Delta} P_N u_0 P_L\parenthese{  e^{i(t-s)\Delta} |x|^{-1} f(s)  } ,
    \end{align*}
    implies 
    \begin{align}
        \norm{e^{it\Delta} P_N u_0 \parenthese{ \int_{0}^t e^{i(t-s)\Delta} P_L (V v)(s) \, ds }}_{L_{t}^pL^q_x}  \lesssim C_{p, q}(N, L) \norm{u_0}_{L^2} \norm{| x | Vv}_{L^2_{t,x}} .
    \end{align}
Another application of Corollary~\ref{cor BPSTZ resolvent} yields 
    \begin{align}
        \norm{e^{it\Delta} P_N u_0 \parenthese{ \int_{0}^t e^{i(t-s)\Delta} P_L (V v)(s) \, ds }}_{L_{t}^pL^q_x}  \lesssim C_{p, q}(N, L) \norm{u_0}_{L^2} \norm{v_0}_{L^2}.
    \end{align}
Similarly,
    \begin{align*}
        III, IV \lesssim C_{p, q}(N, L) \norm{u_0}_{L^2} \norm{v_0}_{L^2}.
    \end{align*}

For $p=2$, we would like to interpolate the previous estimates with standard linear Strichartz estimates. Let $r=(d+2)/d$. Using linear Strichartz estimates,  we get
    \begin{align*}
        \norm{e^{it\Delta}\int_{\R} e^{-is \Delta} P_L(Vu) (s) \, ds}_{L^{2r}_{x, t}} &\lesssim \norm{\int_{\R} e^{-is \Delta} Vu (s) \, ds}_{L^2_x} \\
        &\lesssim \norm{| x | Vu}_{L_{t,x}^2} . 
    \end{align*}

    Now, since $2r>2$, Christ-Kiselev  implies
    \begin{align*}
        \norm{e^{it\Delta}\int_{0}^t e^{-is \Delta} P_L(Vv) (s) \, ds}_{L^{2r}_{t,x}} \lesssim  \norm{|x|Vv}_{L_{t,x}^2}. 
    \end{align*}
Since $V$ is a $(\kappa, \delta)$ potential, $|x|V(x) \lesssim |x|^{-1}$. Therefore, Corollary~\ref{cor BPSTZ resolvent} implies
    \begin{align*}
         \norm{| x | Vv}_{L_{t,x}^2} \lesssim_{\delta} \norm{v_0}_{L^2_x} . 
    \end{align*}
    Thus, H\"older's inequality now gives 
    \begin{align}
        \norm{e^{it\Delta} P_N u_0 \parenthese{ \int_{0}^t e^{i(t-s)\Delta} P_L (V v)(s) \, ds }}_{L_{t,x}^{r}} 
        & \leq\norm{e^{it\Delta} P_N u_0}_{L_{t,x}^{2r}} \norm{e^{it\Delta}\int_{0}^t e^{-is \Delta} P_L(Vv) (s) \, ds}_{L^{2r}_{t,x}}\\ 
        & \lesssim \norm{u_0}_{L^2} \norm{v_0}_{L^2} .
\end{align}

    Finally, interpolation yields
    \begin{align*}
        &\quad \norm{e^{it\Delta} P_N u_0 \parenthese{ \int_{0}^t e^{i(t-s)\Delta} P_L (V v)(s) \, ds }}_{L^{2}_{t}L^q_x} \\
        &\leq \parenthese{\norm{e^{it\Delta} P_N u_0 \parenthese{ \int_{0}^t e^{i(t-s)\Delta} P_L (V v)(s) \, ds }}_{L^{r}_{x, t}}}^{\varepsilon}\left(\norm{e^{it\Delta} P_N u_0 \parenthese{ \int_{0}^t e^{i(t-s)\Delta} P_L (V v)(s) \, ds }}_{L^{p_1}_{ t}L^{q_1}_x}\right)^{(1-\varepsilon)}\\
        &\lesssim C_{p_1, q_1}(N, L)^{1-\varepsilon} \norm{u_0}_{L^2} \norm{v_0}_{L^2} ,
    \end{align*}
    where
\begin{align*}
    \frac{1}{2}=\frac{\varepsilon}{r}+\frac{1-\varepsilon}{p_1},\quad
    \frac{1}{q}=\frac{\varepsilon}{r}+\frac{1-\varepsilon}{q_1}.
\end{align*}
    We note here $r= \frac{d+2}{d} < 2$ when $d \geq 3$. 
    This provides the desired bound for $II$. The estimates for $III$ and $IV$ are handled similarly.

    Now we complete the proof of Theorem~\ref{thm PBilinear}.
\end{proof}

For a potential without a singularity, we achieve a better estimate for the $p=2$ case. The proof of the following theorem follows the same structure as Theorem~\ref{thm PBilinear}. The difference is that we gain some regularity via smoothing estimates.

\PPBilinear*

\begin{proof}[Proof of Theorem~\ref{thm PBilinear1}]
First, we write the solution via the Duhamel principle, that is,
\begin{align}
u & = e^{it (\Delta -V)} u_0 = e^{it\Delta} u_0  -i \int_{0}^t e^{i(t-s)\Delta} V u (s) \, ds ,\\
v & = e^{it (\Delta -V)} v_0 = e^{it\Delta} v_0  -i \int_{0}^t e^{i(t-s)\Delta} V v (s) \, ds .
\end{align}

Next, we apply the projections to both sides of the equations above and expand the product into the following terms
\begin{align}
\norm{P_N u P_L v}_{L_{t,x}^2} & \lesssim \norm{ e^{it\Delta} P_N u_0 e^{it\Delta} P_L v_0}_{L_{t,x}^2} \\
& \quad + \norm{e^{it\Delta} P_N u_0 \parenthese{ \int_{0}^t e^{i(t-s)\Delta} P_L (V v)(s) \, ds }}_{L_{t,x}^2} \\
& \quad + \norm{e^{it\Delta} P_L v_0 \parenthese{\int_{0}^t e^{i(t-s)\Delta} P_N (Vu)(s) \, ds }}_{L_{t,x}^2}  \\
& \quad+ \norm{\parenthese{\int_{0}^t e^{i(t-s)\Delta} P_N (Vu)(s) \, ds} \parenthese{\int_{0}^t e^{i(t-s')\Delta} P_L (Vv)(s') \, ds'}}_{L_{t,x}^2} \\
& = : I + II + III + IV . 
\end{align}

Next, we will treat the four terms separately. 
By assumption, we have
\begin{align}
I \lesssim C(N,L) \norm{u_0}_{L^2} \norm{v_0}_{L^2}
\end{align}

 For $II$, we use Theorem~\ref{thm weightedCK} to get
    \begin{align*}
        &\quad \norm{e^{it\Delta} P_N u_0 \parenthese{ \int_{0}^t e^{i(t-s)\Delta} P_L (V v)(s) \, ds }}_{L_{t,x}^2} \\
        &\lesssim \inner{L}^{-\frac{1}{2}}C(N,L)^{1-2\varepsilon}\min\{N,L\}^{(d-2)\varepsilon} \norm{u_0}_{L^2} \norm{\inner{x} V v}_{L^2_{t,x}} , 
    \end{align*} 
    where $p_1,q_1$ are defined in Theorem~\ref{thm PBilinear}.
Now Corollary~\ref{cor BPSTZ resolvent} implies
    \begin{align*}
         \inner{L}^{-1/2}\norm{\langle x \rangle Vv}_{L_t^2 L^2_x} \lesssim \inner{L}^{-1/2}\norm{v_0}_{L^2_x}
    \end{align*}
    This provides the desired bound for $II$. 
    
    Similarly,
    \begin{align*}
        III \lesssim \inner{N}^{-\frac{1}{2}}C(N,L)^{1-2\varepsilon}\min\{N,L\}^{(d-2)\varepsilon} \norm{u_0}_{L^2} \norm{v_0}_{L^2} , 
    \end{align*}
    and 
    \begin{align*}
        IV \lesssim \inner{N}^{-\frac{1}{2}}\inner{L}^{-\frac{1}{2}}C(N,L)^{1-2\varepsilon}\min\{N,L\}^{(d-2)\varepsilon} \norm{u_0}_{L^2} \norm{v_0}_{L^2} .
    \end{align*}
Since $\inner{N}^{-\frac{1}{2}}\inner{L}^{-\frac{1}{2}} \leq \inner{N}^{-\frac{1}{2}}+\inner{L}^{-\frac{1}{2}}$, we finish the proof.
\end{proof}

\section{Bilinear Strichartz with Localized Initial Data (High Modes)}\label{sec BilinearPHigh}


Let $P_N$ be the Littlewood-Paley projector associated with the Euclidean Laplacian. Consider the following linear equation with $V$ a decaying potential:
\begin{align}\label{eqn linschro2}
\begin{cases}
i \partial_t u + \Delta u - Vu =0\\
u(0) = P_N f .
\end{cases}
\end{align}

The following representation of the solution to \eqref{eqn linschro2} will be very useful:
\begin{align}
u = e^{-itH_V} P_N f= e^{it(\Delta - V)} P_N f = e^{it\Delta} P_N f - i \int_0^t e^{i(t-s) \Delta} V (e^{is (\Delta -V)} P_N f) \, ds . 
\end{align}
We will use this identity to prove the following statement:

\bilocdatahighmodes*

\begin{proof}[Proof of Theorem~\ref{thm bilocdata high modes}]
We handle different terms of the product, $e^{-itH_V} P_N f \, e^{-itH_V} P_L g$ separately. 
\begin{align}
&\norm{e^{-itH_V} P_N f \, e^{-itH_V} P_L g}_{L_{t,x}^2\left(\mathbb{R} \times \mathbb{R}^d\right)} \\
& \leq \norm{e^{it\Delta} P_N f \, e^{it\Delta} P_L g}_{L_{t,x}^2\left(\mathbb{R} \times \mathbb{R}^d\right)} + \norm{e^{it\Delta} P_N f \, \int_0^t e^{i(t-s) \Delta} V (e^{-isH_V} P_L g) \, ds}_{L_{t,x}^2\left(\mathbb{R} \times \mathbb{R}^d\right)} \\
& \quad + \norm{e^{it\Delta} P_L g \, \int_0^t e^{i(t-s) \Delta} V (e^{-isH_V} P_N f) \, ds}_{L_{t,x}^2\left(\mathbb{R} \times \mathbb{R}^d\right)} \\
& \quad + \norm{\int_0^t e^{i(t-s) \Delta} V (e^{-isH_V} P_N f) \, ds \, \int_0^t e^{i(t-s) \Delta} V (e^{-isH_V} P_L g) \, ds}_{L_{t,x}^2\left(\mathbb{R} \times \mathbb{R}^d\right)} \\
& = : I + II + III + IV . 
\end{align}
Lemma~\ref{lem alldbourbi} allows us to estimate Term $I$:
\begin{align*}
    \norm{e^{it\Delta} P_N f \, e^{it\Delta} P_L g}_{L_{t,x}^2\left(\mathbb{R} \times \mathbb{R}^d\right)}\lesssim \frac{L^{\frac{d-1}{2}}}{N^{\frac{1}{2}}}\left\|P_Nf\right\|_2\left\|P_Lg\right\|_2.
\end{align*}

For the linear-nonlinear interaction, we have two terms, $II$ and $III$. For Term $II$, we use classical bilinear Strichartz estimates and Lemma~\ref{lem CK} to obtain
\begin{align}
II & = \norm{e^{it\Delta} P_N f \, \int_0^t e^{i(t-s) \Delta} V (e^{-isH_V} P_L g) \, ds}_{L_{t,x}^2\left(\mathbb{R} \times \mathbb{R}^d\right)} \\
& \lesssim \sum_{M} \norm{e^{it\Delta} P_N f \, P_M \int_0^t e^{i(t-s) \Delta} V (e^{-isH_V} P_L g) \, ds}_{L_{t,x}^2\left(\mathbb{R} \times \mathbb{R}^d\right)}  \\
& \lesssim \sum_{ M \leq L} \norm{e^{it\Delta} P_N f \, P_M \int_0^t e^{i(t-s) \Delta} V (e^{-isH_V} P_L g) \, ds}_{L_{t,x}^2\left(\mathbb{R} \times \mathbb{R}^d\right)}\\
&\hspace{1cm}+  \sum_{L<M\leq N} \norm{e^{it\Delta} P_N f \, P_M \int_0^t e^{i(t-s) \Delta} V (e^{-isH_V} P_L g) \, ds}_{L_{t,x}^2\left(\mathbb{R} \times \mathbb{R}^d\right)}  \\
&\hspace{2cm}+  \sum_{N<M} \norm{e^{it\Delta} P_N f \, P_M \int_0^t e^{i(t-s) \Delta} V (e^{-isH_V} P_L g) \, ds}_{L_{t,x}^2\left(\mathbb{R} \times \mathbb{R}^d\right)} \\
& = : II_a + II_b + II_c .
\end{align}


For Terms $II_a$ and $II_b$, first observe that Standard bilinear Strichartz estimates give us
    \begin{align}
        \norm{e^{it\Delta} P_N f \, P_M \int_{\mathbb{R}} e^{i(t-s) \Delta} V (e^{-isH_V} P_L g) \, ds}_{L_{t,x}^{2}\left(\mathbb{R} \times \mathbb{R}^d\right)} \lesssim \frac{M^{\frac{d-1}{2}}}{N^{\frac{1}{2}}}\norm{P_N f}_{L^2}  \norm{ P_M\int_{\mathbb{R}} e^{-is \Delta} V (e^{-isH_V} P_L g) \, ds }_{L^2_{ x} } ,\\
        \norm{e^{it\Delta} P_N f \, P_M \int_{\mathbb{R}} e^{i(t-s) \Delta} V (e^{-isH_V} P_L g) \, ds}_{L_{t,x}^{\infty}\left(\mathbb{R} \times \mathbb{R}^d\right)}  \lesssim (MN)^{\frac{d}{2}}\norm{P_N f}_{L^2}  \norm{ P_M\int_{\mathbb{R}} e^{-is \Delta} V (e^{-isH_V} P_L g) \, ds }_{L^2_{ x} } .
    \end{align}
Hence, interpolation between the bounds above (with weight $\varepsilon_0$ on the second bound) gives
    \begin{align}
        \norm{e^{it\Delta} P_N f \, P_M \int_{\mathbb{R}} e^{i(t-s) \Delta} V (e^{-isH_V} P_L g) \, ds}_{L_{t,x}^{\frac{2}{1-\varepsilon_0}}\left(\mathbb{R} \times \mathbb{R}^d\right)} \\
        \lesssim \frac{M^{\frac{d-1}{2}(1-\varepsilon_0)}}{N^{\frac{1}{2}(1-\varepsilon_0)}} (MN)^{\frac{d}{2}\varepsilon_0} \norm{P_N f}_{L^2}  \norm{ P_M\int_{\mathbb{R}} e^{-is \Delta} V (e^{-isH_V} P_L g) \, ds }_{L^2_{ x} } . \label{eq homMN}
    \end{align}
For $\alpha \in \mathbb{R}$, Corollary~\ref{cor:homsmooth} implies
\begin{align*}
        M^{\alpha}  \norm{ P_M\int_{\mathbb{R}} e^{-is \Delta} V (e^{-isH_V} P_L g) \, ds }_{L^2_{ x} } 
        &\lesssim \norm{|\nabla|^{
        \alpha}\int_{\mathbb{R}} e^{-is \Delta} V (e^{-isH_V} P_L g) \, ds}_{L^2_{ x} }\\
        &\lesssim \norm{|\nabla|^{\alpha-\frac{1}{2}}(\inner{x} V e^{-itH_V} P_L g)}_{L^2_{t, x} }.
    \end{align*}
Corollary~\ref{cor BDsmooth} implies 
    \begin{align*}
        \norm{|\nabla|^{\alpha-\frac{1}{2}}(\inner{x} V e^{-itH_V} P_L g)}_{L^2_{t, x} } \lesssim \norm{|\nabla|^{\alpha-1} P_L g}_{L^2_{x} }.
    \end{align*}
   
Hence
\begin{align}
    \norm{ P_M\int_{\mathbb{R}} e^{-is \Delta} V (e^{-isH_V} P_L g) \, ds }_{L^2_{ x} } 
    &\lesssim M^{-\alpha} \norm{|\nabla|^{\alpha-\frac{1}{2}}(\inner{x} V e^{-itH_V} P_L g)}_{L^2_{t, x} } , \label{eq11}\\
    \norm{ P_M\int_{\mathbb{R}} e^{-is \Delta} V (e^{-isH_V} P_L g) \, ds }_{L^2_{ x} } 
    &\lesssim M^{-\alpha} \norm{|\nabla|^{\alpha-1} P_L g}_{L^2_{x} } \sim  \frac{L^{\alpha -1} }{M^{\alpha}}\norm{P_L g}_{L^2_{x} } , \label{eq9}
\end{align}
which implies a gain of a full derivative when $\alpha \geq 1$. Therefore, when taking $\alpha =0$ and $\alpha = \frac{d-1}{2}+ \frac{1}{2} \varepsilon_0$, we have the unrestricted (in time) estimate
\begin{align}
    \eqref{eq homMN} \lesssim \frac{\min \{M,L\}^{\frac{d-1}{2}+ \frac{1}{2} \varepsilon_0} L^{-1}}{N^{\frac{1}{2}- \frac{d+1}{2} \varepsilon_0}} \norm{P_N f}_{L^2} \norm{P_L g}_{L^2} , 
\end{align}
which implies 
\begin{align}
     \norm{e^{it\Delta} P_N f \, P_M \int_{0}^t e^{i(t-s) \Delta} V (e^{-isH_V} P_L g) \, ds}_{L_{t,x}^{\frac{2}{1-\varepsilon_0}}\left(\mathbb{R} \times \mathbb{R}^d\right)} \lesssim \frac{\min \{M,L\}^{\frac{d-1}{2}+ \frac{1}{2} \varepsilon_0} L^{-1}}{N^{\frac{1}{2}- \frac{d+1}{2} \varepsilon_0}} \norm{P_N f}_{L^2} \norm{P_L g}_{L^2} . \label{eq12}
\end{align}
 We interpolate with the following estimate following from linear Strichartz, H\"older, Christ–Kiselev and \eqref{eq9} (with $\alpha = 0$):
    \begin{align}
        & \norm{e^{it\Delta} P_N f \, P_M \int_{0}^t e^{i(t-s) \Delta} V (e^{-isH_V} P_L g) \, ds}_{L_{t,x}^{\frac{d+2}{d}}\left(\mathbb{R} \times \mathbb{R}^d\right)} \\
        & \lesssim \norm{e^{it\Delta} P_N f}_{L_{t,x}^{\frac{2(d+2)}{d}}} \norm{P_M \int_{0}^t e^{i(t-s) \Delta} V (e^{-isH_V} P_L g) \, ds}_{L_{t,x}^{\frac{2(d+2)}{d}}} \\
        & \lesssim \norm{P_N f}_{L^2} \norm{|\nabla|^{-\frac{1}{2}}(\inner{x} V e^{-itH_V} P_L g)}_{L^2_{t, x} }\\
        & \lesssim \frac{1}{L}\norm{P_N f}_{L^2}  \norm{P_L g}_{L^2 } . \label{eq13}
    \end{align}
    Interpolating the two time restricted estimates \eqref{eq12} and \eqref{eq13}, we obtain
    \begin{align*}
        \norm{e^{it\Delta} P_N f \, P_M \int_{0}^t e^{i(t-s) \Delta} V (e^{-isH_V} P_L g) \, ds}_{L_{t,x}^{2}\left(\mathbb{R} \times \mathbb{R}^d\right)} & \lesssim \frac{\min \{M,L\}^{\frac{d-1}{2}-\varepsilon}L^{-1}}{N^{\frac{1}{2}-\varepsilon}}\norm{P_N f}_{L^2}\norm{P_L g}_{L^2} \\
        & = \frac{\min\{M,N,L\}^{\frac{d-1}{2}-\varepsilon}}{\max\{M,N\}^{\frac{1}{2}-\varepsilon} L}\norm{P_N f}_{L^2}\norm{P_L g}_{L^2} 
    \end{align*}
    with 
    \begin{align}
        \varepsilon = \frac{d^2 \varepsilon_0}{2(d-2 + (d+2)\varepsilon_0)}
    \end{align} 
    and the interpolation weight 
    \begin{align}
        \theta = \frac{d-2}{d-2 + (d+2) \varepsilon_0}
    \end{align}
    on \eqref{eq12}. 
    In fact, we also have an estimate for the indeterminate step (from Theorem~\ref{thm weightedCK})
    \begin{align}\label{eq10}
        & \norm{e^{it\Delta} P_N f \, P_M \int_{0}^t e^{i(t-s) \Delta} V (e^{-isH_V} P_L g) \, ds}_{L_{t,x}^{2}\left(\mathbb{R} \times \mathbb{R}^d\right)} \\
        & \lesssim \frac{\min\{M,N\}^{\frac{d-1}{2}-\varepsilon}}{\max \{M,N\}^{\frac{1}{2}-\varepsilon} M^{\frac{1}{2}}}  \norm{P_N f}_{L^2}   \norm{(\inner{x} V e^{-itH_V} P_L g)}_{L^2_{t, x} } . 
    \end{align}

Now, for Term $II_a$, where $M<L$, we have
\begin{align}
II_a & \lesssim \sum_{ M \leq L} \frac{M^{\frac{d-1}{2}-\varepsilon}}{N^{\frac{1}{2}-\varepsilon}L} \norm{P_N f}_{L^2}  \norm{P_L g}_{L^2 }  \lesssim  \frac{L^{\frac{d-1}{2}-1-\varepsilon}}{N^{\frac{1}{2}-\varepsilon}} \norm{P_N f}_{L^2}  \norm{P_L g}_{L^2 } .
\end{align}
    For Term $II_b$,
\begin{align}
II_b & \lesssim \sum_{ L<M\leq N} \frac{L^{\frac{d-1}{2}-1-\varepsilon}}{N^{\frac{1}{2}-\varepsilon}}\norm{P_N f}_{L^2}  \norm{P_L g}_{L^2 }  \lesssim \log(N/L) \frac{L^{\frac{d-1}{2}-1-\varepsilon}}{N^{\frac{1}{2}-\varepsilon}}\norm{P_N f}_{L^2}  \norm{P_L g }_{L^2_{x} } .
\end{align}

For high $M$, Term $II_c$, a similar calculation yields
\begin{align*}
    II_c & \lesssim\sum_{M > N} \frac{1}{M^{\frac{1}{2}-\varepsilon}} \norm{P_N f}_{L^2}  \norm{\nabla^{\frac{d-1}{2}-\frac{1}{2}-\varepsilon} [\inner{x} V e^{-itH_V} P_L g ]}_{L^2_{t, x} }\\
     & \lesssim\frac{L^{\frac{d-1}{2}-1-\varepsilon}}{N^{\frac{1}{2}-\varepsilon}} \norm{P_N f}_{L^2}  \norm{ P_L g }_{L^2 }.
\end{align*}

Combining terms, we have
\begin{align*}
   II & \lesssim II_a + II_b + II_c \lesssim \frac{L^{\frac{d-1}{2}-1 -\varepsilon}}{N^{\frac{1}{2}-\varepsilon}}\norm{P_N f}_{L^2}  \norm{ P_L g}_{L^2 } .
\end{align*}

For the other linear-nonlinear term, $III$, we use a similar decomposition and \eqref{eq10}
\begin{align}
III & = \norm{e^{it\Delta} P_L g \, \int_0^t e^{i(t-s) \Delta} V (e^{-isH_V} P_N f) \, ds}_{L_{t,x}^2\left(\mathbb{R} \times \mathbb{R}^d\right)} \\
& \lesssim \sum_{M} \norm{e^{it\Delta} P_L g \, P_M \int_0^t e^{i(t-s) \Delta} V (e^{-isH_V} P_N f) \, ds}_{L_{t,x}^2\left(\mathbb{R} \times \mathbb{R}^d\right)}  \\
& \lesssim \sum_{M \leq L} \frac{M^{\frac{d-1}{2}-\frac{1}{2}-\varepsilon}}{L^{\frac{1}{2}-\varepsilon}} \norm{P_L g}_{L^2}  \norm{\inner{x} V (e^{-itH_V} P_N f)}_{L^2_{t, x} } \\
&\hspace{.7cm}+  \sum_{M>L} \frac{L^{\frac{d-1}{2}-\varepsilon}}{M^{1-\varepsilon}}  \norm{P_L g}_{L^2} \norm{\inner{x} V (e^{-itH_V} P_N f)}_{L^2_{t, x} } .
\end{align}

In both cases, summing over $M$ gives 
    \begin{align*}
        III & = \norm{e^{it\Delta} P_L g \, \int_0^t e^{i(t-s) \Delta} V (e^{-isH_V} P_N f) \, ds}_{L_{t,x}^2\left(\mathbb{R} \times \mathbb{R}^d\right)} \\
        &\lesssim L^{\frac{d-1}{2}-1}\norm{P_L g}_{L^2}  \norm{\inner{x} V (e^{-itH_V} P_N f)}_{L^2_{t, x} }  . 
    \end{align*}
Observe that Corollary~\ref{cor BDsmooth} implies that 
    \begin{align*}
        \norm{\inner{x} V (e^{-itH_V} P_N f)}_{L^2_{t, x} } \lesssim \norm{ P_N f}_{H^{-1/2}_{ x} } \lesssim N^{-\frac{1}{2}} \norm{P_Nf}_{L^2_x} .
    \end{align*}
Therefore,
    \begin{align*}
       III & \lesssim \frac{L^{\frac{d-1}{2}-1}}{N^{\frac{1}{2}}}\norm{P_L g}_{L^2}  \norm{ P_N f}_{L^2 }. 
    \end{align*}

For the nonlinear-nonlinear term $IV$, we perform a decomposition for both terms simultaneously. Applying the same interpolation argument as in Cases $II$ and $III$, we have 
\begin{align}
    & \norm{P_K \int_0^t e^{i(t-s)\Delta} V(e^{-isH_V}P_N f) \, ds \, P_M \int_0^t e^{i(t-s)\Delta} V(e^{-isH_V}P_L g) \, ds}_{L^2_{t,x}\left(\mathbb{R}\times\mathbb{R}^d\right)}\\
    & \lesssim \frac{\min\{K,M\}^{\frac{d-1}{2}-\varepsilon}}{\max\{K,M\}^{\frac{1}{2}-\varepsilon}K^\alpha M^\beta} \norm{|\nabla|^{\alpha-\frac{1}{2}} (\inner{x}V e^{-itH_V}P_N f)}_{L^2_{t,x}} \norm{|\nabla|^{\beta-\frac{1}{2}} (\inner{x}V e^{-itH_V}P_L g)}_{L^2_{t,x}}, \label{eq14}
\end{align}
for $\alpha , \beta \geq 0$.

Similar to the previous decomposition
\begin{align}
    IV & = \norm{\int_{0}^t e^{i(t-s) \Delta} V (e^{-isH_V} P_N f) \, ds\, \int_{0}^t e^{i(t-s) \Delta} V (e^{-isH_V} P_L g) \, ds}_{L_{t,x}^2\left(\mathbb{R} \times \mathbb{R}^d\right)}\\
    & \lesssim \sum_{1\leq M\leq K}
    \norm{
    P_K \int_0^t e^{i(t-s)\Delta}
    V(e^{-isH_V}P_N f) \, ds \,
    P_M \int_0^t e^{i(t-s)\Delta}
    V(e^{-isH_V}P_L g) \, ds}_{L^2_{t,x}} \\
    & \quad +\sum_{1\leq K<M}
    \norm{
    P_K \int_0^t e^{i(t-s)\Delta}
    V(e^{-isH_V}P_N f) \, ds \,
    P_M \int_0^t e^{i(t-s)\Delta}
    V(e^{-isH_V}P_L g) \, ds}_{L^2_{t,x}} \\
    & \quad +\sum_{\min\{K,M\}<1}
    \norm{
    P_K \int_0^t e^{i(t-s)\Delta}
    V(e^{-isH_V}P_N f) \, ds \,
    P_M \int_0^t e^{i(t-s)\Delta}
    V(e^{-isH_V}P_L g) \, ds}_{L^2_{t,x}} \\
    & =: IV_a+IV_b+IV_c.
\end{align}

For the first sum, we absorb the $M^{\frac{d-1}{2}}$ term with Corollary~\ref{cor BDsmooth} on the $P_Lg$ term. We gain a finally factor of $N^{-1/2}$ from the same Corollary~\ref{cor BDsmooth} argument we used for $III$ (with $\alpha = \frac{1}{2}$ and $\beta = \frac{d-1}{2}-\varepsilon$ in \eqref{eq14})
\begin{align*}
    IV_a & \lesssim \sum_{1\leq M\leq K} \frac{M^{\frac{d-1}{2}-\varepsilon}}{K^{\frac{1}{2}-\varepsilon} K^{\frac{1}{2}}M^{\frac{d-1}{2}-\varepsilon}} \norm{\inner{x}V e^{-itH_V}P_N f}_{L^2_{t,x}} \norm{|\nabla|^{\frac{d-1}{2}-\frac{1}{2}-\varepsilon} (\inner{x}V e^{-itH_V}P_L g)}_{L^2_{t,x}}\\
    & \lesssim \parenthese{\sum_{K\geq1}\frac{\log K}{K^{1-\varepsilon}}} \norm{\inner{x}V e^{-itH_V}P_N f}_{L^2_{t,x}} \norm{|\nabla|^{\frac{d-1}{2}-\frac{1}{2}-\varepsilon} (\inner{x}V e^{-itH_V}P_L g)}_{L^2_{t,x}} \\
    & \lesssim \norm{\inner{x}V e^{-itH_V}P_N f}_{L^2_{t,x}} L^{\frac{d-1}{2}-1-\varepsilon} \norm{P_L g}_{L^2} \\
    & \lesssim \frac{L^{\frac{d-1}{2}-1-\varepsilon}} {N^{\frac{1}{2}}} \norm{P_N f}_{L^2} \norm{P_L g}_{L^2}.
\end{align*}

For the second sum, we perform a similar computation:
\begin{align*}
    IV_b & \lesssim \sum_{1\leq K<M} \frac{K^{\frac{d-1}{2}-\frac{1}{2}-\varepsilon}}{M^{\frac{d-1}{2}+\frac{1}{2}-2\varepsilon}} \norm{\inner{x}V e^{-itH_V}P_N f}_{L^2_{t,x}} \norm{|\nabla|^{\frac{d-1}{2}-\frac{1}{2}-\varepsilon} (\inner{x}V e^{-itH_V}P_L g)}_{L^2_{t,x}} \\
    & \lesssim \sum_{1\leq K<M} \frac{1}{M^{1-\varepsilon}} \norm{\inner{x}V e^{-itH_V}P_N f}_{L^2_{t,x}} \norm{|\nabla|^{\frac{d-1}{2}-\frac{1}{2}-\varepsilon} (\inner{x}V e^{-itH_V}P_L g)}_{L^2_{t,x}} \\
    & \lesssim \parenthese{\sum_{M\geq1}\frac{\log M}{M^{1-\varepsilon}}} \norm{\inner{x}V e^{-itH_V}P_N f}_{L^2_{t,x}} \norm{|\nabla|^{\frac{d-1}{2}-\frac{1}{2}-\varepsilon} (\inner{x}V e^{-itH_V}P_L g)}_{L^2_{t,x}} \\
    & \lesssim \norm{\inner{x}V e^{-itH_V}P_N f}_{L^2_{t,x}} L^{\frac{d-1}{2}-1-\varepsilon} \norm{P_L g}_{L^2} \\
    & \lesssim \frac{L^{\frac{d-1}{2}-1-\varepsilon}}{N^{\frac{1}{2}}} \norm{P_N f}_{L^2} \norm{P_L g}_{L^2}.
\end{align*}

Finally, for $IV_c$ we take $\alpha=\beta=0$. The corresponding frequency sum is finite, since
\begin{align*}
    & \sum_{\min\{K,M\}<1} \frac{\min\{K,M\}^{\frac{d-1}{2}-\varepsilon}} {\max\{K,M\}^{\frac{1}{2}-\varepsilon}}\\
    & \lesssim \sum_{K<1} \frac{1}{K^{\frac{1}{2}-\varepsilon}} \sum_{M\leq K}M^{\frac{d-1}{2}-\varepsilon} + \sum_{K\geq1} \frac{1}{K^{\frac{1}{2}-\varepsilon}} \sum_{M<1}M^{\frac{d-1}{2}-\varepsilon} \\
    & \lesssim \sum_{K<1}K^{\frac{d-2}{2}} + \sum_{K\geq1}\frac{1}{K^{\frac{1}{2}-\varepsilon}} \\
    & \lesssim 1.
\end{align*}
Therefore, Corollary~\ref{cor BDsmooth} with $s=0$ gives
\begin{align*}
    IV_c & \lesssim \norm{|\nabla|^{-1/2} (\inner{x}V e^{-itH_V}P_N f)}_{L^2_{t,x}} \norm{|\nabla|^{-1/2} (\inner{x}V e^{-itH_V}P_L g)}_{L^2_{t,x}} \\
    & \lesssim \frac{1}{N L} \norm{P_N f}_{L^2} \norm{P_L g}_{L^2}.
\end{align*}

Then
\begin{align}
    IV & \lesssim \frac{L^{\frac{d-1}{2}-1-\varepsilon}}{N^{\frac{1}{2}}} \norm{P_N f}_{L^2} \norm{P_L g}_{L^2}.
\end{align}
Combining all terms $I$, $II$, $III$, and $IV$, we conclude the desired estimate.

Now we complete the proof of Theorem~\ref{thm bilocdata high modes}.

\end{proof}


\section{Bilinear Strichartz with Localized Initial Data (Low Frequencies)}\label{sec BilinearPLow}

We would like to use a scaling argument to achieve bilinear estimates in the low frequency regime. As we saw in the previous section, high frequency expressions of local smoothing are beneficial, but low frequencies present issues. In order to handle low frequency estimates, we will need subcritical decay of the potential to achieve similar bilinear estimates.  Consequently, we obtain weaker bilinear estimates.

\bilocdatalowmodes*

    

\begin{proof}[Proof of Theorem~\ref{thm bilocdata lowmodes}]
 Let $0\leq \varepsilon \leq\min(\eta-2, 1)$. Just as in the proof of Theorem~\ref{thm bilocdata high modes}, we handle different terms of the product, $e^{-itH_V} P_N f \, e^{-itH_V} P_L g$, separately. 
\begin{align}
&\norm{e^{-itH_V} P_N f \, e^{-itH_V} P_L g}_{L_{t,x}^2\left(\mathbb{R} \times \mathbb{R}^d\right)} \\
& \leq \norm{e^{it\Delta} P_N f \, e^{it\Delta} P_L g}_{L_{t,x}^2\left(\mathbb{R} \times \mathbb{R}^d\right)} +  \norm{e^{it\Delta} P_L g \, \int_0^t e^{i(t-s) \Delta} V (e^{-isH_V} P_N f) \, ds}_{L_{t,x}^2\left(\mathbb{R} \times \mathbb{R}^d\right)} \\
& \quad + \norm{e^{it\Delta} P_N f \, \int_0^t e^{i(t-s) \Delta} V (e^{-isH_V} P_L g) \, ds}_{L_{t,x}^2\left(\mathbb{R} \times \mathbb{R}^d\right)}\\
& \quad + \norm{\int_0^t e^{i(t-s) \Delta} V (e^{-isH_V} P_N f) \, ds \, \int_0^t e^{i(t-s) \Delta} V (e^{-isH_V} P_L g) \, ds}_{L_{t,x}^2\left(\mathbb{R} \times \mathbb{R}^d\right)} \\
& = : I + II + III + IV . 
\end{align}
The hypothesis allows us to estimate Term $I$:
\begin{align*}
    \norm{e^{it\Delta} P_N f \, e^{it\Delta} P_L g}_{L_{t,x}^2\left(\mathbb{R} \times \mathbb{R}^d\right)}\lesssim \frac{L^{\frac{d-1}{2}}}{N^{\frac{1}{2}}}\left\|P_Nf\right\|_2\left\|P_Lg\right\|_2.
\end{align*}

For the linear-nonlinear interaction, we have two terms, $II$ and $III$. For Term $II$, we use a frequency decomposition argument as seen in the proof of Theorem~\ref{thm bilocdata high modes}
\begin{align}
II & = \norm{e^{it\Delta} P_L g \, \int_0^t e^{i(t-s) \Delta} V (e^{-isH_V} P_N f) \, ds}_{L^2_{t, x} } \\
& \lesssim \sum_{M} \norm{e^{it\Delta} P_L g \, P_M \int_0^t e^{i(t-s) \Delta} V (e^{-isH_V} P_N f) \, ds}_{L^2_{t, x} }  \\
& \lesssim   \sum_{M\leq L} \norm{e^{it\Delta} P_L g \, P_M \int_0^t e^{i(t-s) \Delta} V (e^{-isH_V} P_N f) \, ds}_{L^2_{t, x} }  \\
&\hspace{2cm}+  \sum_{L<M} \norm{e^{it\Delta} P_L g \, P_M \int_0^t e^{i(t-s) \Delta} V (e^{-isH_V} P_N f) \, ds}_{L^2_{t, x} } . 
\end{align}
By Lemma~\ref{lem CK}, we can take $\varepsilon_1<1$ to get 
\begin{align}
    II & \lesssim \sum_{M<L} \norm{e^{it\Delta} P_L g \, P_M \int_0^t e^{i(t-s) \Delta} V (e^{-isH_V} P_N f) \, ds}_{L^2_{t, x} }  \\
&\hspace{2cm}+  \sum_{L<M} \norm{e^{it\Delta} P_L g \, P_M \int_0^t e^{i(t-s) \Delta} V (e^{-isH_V} P_N f) \, ds}_{L^2_{t, x} }\\
&\lesssim \sum_{M \leq L} M^{\frac{d-2}{2}} \parenthese{\frac{M}{L}}^{(1-\varepsilon_1)/2}   \norm{g}_{L^2} \norm{V (e^{-itH_V} P_N f)}_{L_t^{2} L_x^{\frac{2d}{d+2}, 2}}\\
&\hspace{2cm}+\sum_{L< M} L^{\frac{d-2}{2}} \parenthese{\frac{L}{M}}^{(1-\varepsilon_1)/2}  \norm{g}_{L^2} \norm{V (e^{-itH_V} P_N f)}_{L_t^{2} L_x^{\frac{2d}{d+2}, 2}} . 
\end{align}

Now we can use Theorem~\ref{thm genmorawetz} (Corollary~\ref{cor flatgenmorawetz})  along with $\varepsilon < \min (2, d/2-1)$ to conclude that
\begin{align}
\norm{V (e^{-itH_V} P_N f)}_{L_t^2 L_x^{\frac{2d}{d+2},2} } &\lesssim \norm{|x|^{1+\varepsilon} V}_{L_t^{\infty} L_x^{d,\infty}} \norm{|x|^{-(1+\varepsilon)} (e^{-itH_V} P_N f)}_{L_{t}^2 L_x^{2,2}}  \\
&= \norm{|x|^{1+\varepsilon} V}_{L_t^{\infty} L_x^{d,\infty}} \norm{|x|^{-(1+\varepsilon)}(e^{-itH_V} P_N f)}_{L_{t,x}^2} \\
&\lesssim \norm{(-\Delta)^{\varepsilon/2}P_Nf}_{L^2}\\
&\lesssim N^{\varepsilon}\norm{P_Nf}_{L^2} . 
\end{align}
In both cases, summing over $M$ gives 
    \begin{align*}
        II & = \norm{e^{it\Delta} P_L g \, \int_0^t e^{i(t-s) \Delta} V (e^{-isH_V} P_N f) \, ds}_{L^2_{t, x}} \\
        &\lesssim L^{\frac{d-2}{2}}N^{\varepsilon}\norm{P_L g}_{L^2}  \norm{ P_N f }_{L^2 }  .
    \end{align*}
Similarly, for term III we get 
    \begin{align*}
        III \lesssim N^{\frac{d-2}{2}}L^{\varepsilon}\norm{P_L g}_{L^2}  \norm{ P_N f }_{L^2 } .
    \end{align*}

For the nonlinear-nonlinear term $IV$, 
\begin{align*}
    IV & = \norm{\int_0^t e^{i(t-s) \Delta} V (e^{-isH_V} P_N f) \, ds\, \int_0^t e^{i(t-s) \Delta} V (e^{-isH_V} P_L g) \, ds}_{L^2_{t, x}}\\
    & \lesssim \sum_{M<K} \norm{P_K\int_0^t e^{i(t-s) \Delta} V (e^{-isH_V} P_N f) \, ds \, P_M \int_0^t e^{i(t-s) \Delta} V (e^{-isH_V} P_L g) \, ds}_{L^2_{t, x} }  \\
&\hspace{2cm}+  \sum_{K<M} \norm{P_K\int_0^t e^{i(t-s) \Delta} V (e^{-isH_V} P_N f) \, ds \, P_M \int_0^t e^{i(t-s) \Delta} V (e^{-isH_V} P_L g) \, ds}_{L^2_{t, x} }
\end{align*}
By the same argument as for item $IV$ in the proof of Theorem~\ref{thm bilocdata high modes}, we use inequality \eqref{eq10} with some $\varepsilon_2<1/2$ to get
\begin{align*}
    IV  &\lesssim \sum_{ K<M \atop 1<M} \frac{K^{\frac{d-1}{2}-\frac{1}{2}-\varepsilon_2}}{M^{\frac{d-1}{2}+\frac{1}{2}}} \norm{\inner{x}V e^{-itH_V}P_N f}_{L^2_{t,x}} \norm{|\nabla|^{\frac{d-1}{2}-\frac{1}{2}+\varepsilon_2}(\inner{x}V e^{-itH_V}P_L g)}_{L^2_{t,x}}  \\
&\hspace{1cm}+   \sum_{M<K \atop 1<K} \frac{M^{\frac{d-1}{2}-\frac{1}{2}-\varepsilon_2}}{K^{\frac{d-1}{2}+\frac{1}{2}}} \norm{|\nabla|^{\frac{d-1}{2}-\frac{1}{2}+\varepsilon_2} (\inner{x}V e^{-itH_V}P_N f)}_{L^2_{t,x}} \norm{\inner{x}V e^{-itH_V}P_L g}_{L^2_{t,x}} \\
&\hspace{1cm}+\sum_{ K<M <1} \frac{K^{\frac{d-1}{2}-\frac{1}{2}-\varepsilon_2}}{M^{\frac{d-1}{2}+\frac{1}{2}-2\varepsilon_2}} \norm{\inner{x}V e^{-itH_V}P_N f}_{L^2_{t,x}} \norm{|\nabla|^{\frac{d-1}{2}-\frac{1}{2}-\varepsilon_2} (\inner{x}V e^{-itH_V}P_L g)}_{L^2_{t,x}} \\
&\hspace{1cm}+   \sum_{M<K<1} \frac{M^{\frac{d-1}{2}-\frac{1}{2}-\varepsilon_2}}{K^{\frac{d-1}{2}+\frac{1}{2}-2\varepsilon_2}} \norm{|\nabla|^{\frac{d-1}{2}-\frac{1}{2}-\varepsilon_2} (\inner{x}V e^{-itH_V}P_N f)}_{L^2_{t,x}} \norm{\inner{x}V e^{-itH_V}P_L g}_{L^2_{t,x}}
\end{align*}
Since $|\inner{x}V(x)| \leq |x|^{-(1+\varepsilon)}$, we can conclude that 
\begin{align*}
    IV  &\lesssim \sum_{1<M} M^{-\varepsilon_2} \norm{|x|^{-(1+\varepsilon)} e^{-itH_V}P_N f}_{L^2_{t,x}} \norm{|\nabla|^{\frac{d-1}{2}-\frac{1}{2}+\varepsilon_2} (\inner{x}Ve^{-itH_V}P_L g)}_{L^2_{t,x}}  \\
&\hspace{1cm}+   \sum_{1<K} K^{-\varepsilon_2} \norm{|\nabla|^{\frac{d-1}{2}-\frac{1}{2}+\varepsilon_2} (\inner{x}Ve^{-itH_V}P_N f)}_{L^2_{t,x}} \norm{|x|^{-(1+\varepsilon)} e^{-itH_V}P_L g}_{L^2_{t,x}}  \\
&\hspace{1cm}+\sum_{ M <1} M^{\varepsilon_2} \norm{|x|^{-(1+\varepsilon)}e^{-itH_V}P_N f}_{L^2_{t,x}} \norm{ |\nabla|^{\frac{d-1}{2}-\frac{1}{2}-\varepsilon_2} (\inner{x}Ve^{-itH_V}P_L g)}_{L^2_{t,x}}  \\
&\hspace{1cm}+   \sum_{K<1} K^{\varepsilon_2}\norm{|\nabla|^{\frac{d-1}{2}-\frac{1}{2}-\varepsilon_2} (\inner{x}Ve^{-itH_V}P_N f)}_{L^2_{t,x}} \norm{|x|^{-(1+\varepsilon)} e^{-itH_V}P_L g}_{L^2_{t,x}} 
\end{align*}
Summing everything we get 
\begin{align*}
    IV  &\lesssim  \norm{|x|^{-(1+\varepsilon)} e^{-itH_V}P_N f}_{L^2_{t,x}} \norm{|\nabla|^{\frac{d-1}{2}-\frac{1}{2}+\varepsilon_2} (\inner{x}Ve^{-itH_V}P_L g)}_{L^2_{t,x}}  \\
&\hspace{1cm}+    \norm{|\nabla|^{\frac{d-1}{2}-\frac{1}{2}+\varepsilon_2} (\inner{x}Ve^{-itH_V}P_N f)}_{L^2_{t,x}} \norm{|x|^{-(1+\varepsilon)} e^{-itH_V}P_L g}_{L^2_{t,x}}  \\
&\hspace{1cm}+ \norm{|x|^{-(1+\varepsilon)}e^{-itH_V}P_N f}_{L^2_{t,x}} \norm{ |\nabla|^{\frac{d-1}{2}-\frac{1}{2}-\varepsilon_2} (\inner{x}Ve^{-itH_V}P_L g)}_{L^2_{t,x}}  \\
&\hspace{1cm}+   \norm{|\nabla|^{\frac{d-1}{2}-\frac{1}{2}-\varepsilon_2} (\inner{x}Ve^{-itH_V}P_N f)}_{L^2_{t,x}} \norm{|x|^{-(1+\varepsilon)} e^{-itH_V}P_L g}_{L^2_{t,x}} 
\end{align*}
Now Corollary~\ref{cor flatgenmorawetz} and Corollary~\ref{cor BDsmooth} imply 
\begin{align*}
    IV  &\lesssim  L^{\varepsilon} \norm{P_L g}_{L^2_{x}} (\norm{P_N f}_{H^{(d-3)/2+\varepsilon_2}_{x}}+\norm{P_N f}_{H^{(d-3)/2-\varepsilon_2}_{x}} )\\
    &+ N^{\varepsilon} \norm{P_N f}_{L^2_{x}} (\norm{P_L g}_{H^{(d-3)/2+\varepsilon_2}_{x}}+\norm{P_L g}_{H^{(d-3)/2-\varepsilon_2}_{x}} ).
\end{align*}
Since $N \leq 1$, the inhomogeneity of the Sobolev norm implies 
\begin{align*}
    IV  &\lesssim   (L^{\varepsilon}+N^{\varepsilon}) \norm{P_L g}_{L^2_{x}} \norm{P_N f}_{L^2_{x}}.
\end{align*}
Combining all terms $I$, $II$, $III$, and $IV$, we conclude the desired estimate.

Now we complete the proof of Theorem~\ref{thm bilocdata lowmodes}.
\end{proof}

\section{Commutator Estimates}\label{sec Commbi}

In this section, we combine the two types of bilinear estimates we obtained earlier and show a commutator estimate.

\commbi*

    \begin{proof}[Proof of Theorem~\ref{thm commbi}]
        We will again use the identities
        \begin{align}
            e^{it(\Delta - V)} P_M f &= e^{it\Delta} P_M f - i \int_0^t e^{i(t-s) \Delta} V (e^{-isH_V} P_M f) \, ds ,\\
            e^{it(\Delta - V)} P_K g &= e^{it\Delta} P_K g - i \int_0^t e^{i(t-s) \Delta} V (e^{-isH_V} P_K g) \, ds . 
        \end{align}
        The standard Littlewood-Paley projection commutes with the evolution of the flat Schr\"odinger operator, therefore
        \begin{align*}
            P_N e^{-itH_V}P_Mf&= e^{it\Delta} P_NP_M f - i P_N\int_0^t e^{i(t-s) \Delta} V (e^{-isH_V} P_M f) \, ds\\
            P_Le^{-itH_V} P_K g &= e^{it\Delta} P_LP_K g - i P_L\int_0^t e^{i(t-s) \Delta} V (e^{-isH_V} P_K g) \, ds
        \end{align*}

        The linear-linear term can be handled easily with bilinear Strichartz:
        \begin{align*}
            \norm{ e^{it\Delta}P_N P_Mf \, e^{it\Delta}P_LP_Kg}_{L_{t,x}^2\left(\mathbb{R} \times \mathbb{R}^d\right)} \lesssim C(N,L)\norm{P_NP_Mf}_{L^2} \norm{P_LP_Kg}_{L^2} .
        \end{align*}
    For nonlinear-linear, we have, by Theorem~\ref{thm weightedCK}, 
        \begin{align*}
            &\norm{P_Ne^{it\Delta} P_M f \, P_L\int_0^t e^{i(t-s) \Delta} V e^{-is H_V} P_K g \, ds}_{L_{t,x}^2\left(\mathbb{R} \times \mathbb{R}^d\right)}\\
            &= \norm{P_Ne^{it\Delta} P_NP_M f \, \int_0^t e^{i(t-s) \Delta} P_LV e^{-is H_V} P_K g \, ds}_{L_{t,x}^2\left(\mathbb{R} \times \mathbb{R}^d\right)}\\
            &\lesssim C(N, L)^{1-2\varepsilon}\min\{N, L\}^{(d-2)\varepsilon} \inner{L}^{-1/2} \norm{P_NP_M f}_{L^2} \norm{\inner{x} V e^{-is H_V} P_K g}_{L_{t,x}^2} . 
        \end{align*}
        
By Corollary~\ref{cor BDsmooth},
        \begin{align*}
            \norm{\inner{x}^{-1}e^{-is H_V} P_K g}_{L_{t,x}^2} \lesssim \inner{K}^{-1/2}\norm{P_K g }_{L^2_x }.
        \end{align*}

Therefore,
     \begin{align*}
            &\norm{P_Me^{it\Delta} P_N f \, P_L\int_0^t e^{i(t-s) \Delta} V (e^{-isH_V} P_K g) \, ds}_{L_{t,x}^2\left(\mathbb{R} \times \mathbb{R}^d\right)}\\
            &\lesssim \left[C(N, L)^{1-2\varepsilon}\min\{N, L\}^{(d-2)\varepsilon}(\inner{K}\inner{L})^{-1/2}\norm{P_N P_M}_{L^2 \to L^2} \right] \norm{P_Mf}_{L^2}  \norm{P_K g}_{L^2 } .
        \end{align*}

 The other linear-nonlinear term can be handled similarly.

The nonlinear-nonlinear term is handled just as easily with the second estimate of Theorem~\ref{thm weightedCK}:
    \begin{align*}
        &\norm{P_N\int_0^t e^{i(t-s) \Delta} V (e^{-is H_V} P_M f) \, ds\, P_L\int_0^t e^{i(t-s) \Delta} V (e^{-is H_V} P_K g) \, ds}_{L_{t,x}^2\left(\mathbb{R} \times \mathbb{R}^d\right)}\\
       &\lesssim C(N, L)^{1-2\varepsilon}\min\{N, L\}^{(d-2)\varepsilon} \inner{L}^{-1/2} \inner{N}^{-1/2}\norm{\inner{x} V e^{-is H_V} P_M f}_{L_{t,x}^2} \norm{\inner{x} V e^{-is H_V} P_K g}_{L_{t,x}^2}\\
        &\lesssim  C(N, L)^{1-2\varepsilon}\min\{N, L\}^{(d-2)\varepsilon} \inner{L}^{-1/2} \inner{N}^{-1/2}\inner{M}^{-1/2} \inner{K}^{-1/2}\norm{ P_M f}_{ L_x^{2}} \norm{P_K g}_{ L_x^{2}} .
    \end{align*}

    Now we complete the proof of Theorem~\ref{thm commbi}.
\end{proof}

\appendix

\section{Proof of Lemma~\ref{lem KMVZZ smoothing}}\label{app KMVZZsmoothing}

Define the vector field $b : \mathbb{R}^d \to \mathbb{R}^d$ by 
    \begin{align*}
        b(x) := \frac{x}{R\langle x/R\rangle}.
    \end{align*}
The vector field, $b$, is designed so that 
    \begin{enumerate}[(a)]
        \item $b$ is bounded,
        \item The divergence of $b$ satisfies $\langle \nabla, b\rangle \geq 0$,
        \item $-\Delta (\langle \nabla, b\rangle )\geq 0$,
        \item $\langle x, b \rangle \leq \frac{|x|^2}{R}$,
        \item The symmetric $d\times d$ matrix, $\nabla b$, satisfies 
            \begin{align*}
                (\nabla b)_{jk} = \frac{\delta_{jk}}{R\inner{x/R}}- \frac{x_j x_k}{R^3\inner{x/R}^3} \geq \frac{\delta_{jk}}{R\inner{x/R}^3} .
            \end{align*}
    \end{enumerate}

Define the momentum functional $F$ by 
    \begin{align*}
        F(t):= \int_{\mathbb{R}^d} 2\langle\Im(\Bar{u}\nabla u), b\rangle.
    \end{align*}
Conservation of the energy, $\mathcal{E}(u):= \int |\nabla u|^2+ \int V|u|^2$, and Hardy's inequality implies $\|\nabla u\|_{L^2} \lesssim \|\nabla u_0\|_{L^2}$. In addition, $b$ is bounded, so Cauchy-Schwarz implies
    \begin{align*}
        |F(t)| \lesssim \|u_0\|_{L^2} \|\nabla u_0\|_{L^2} . 
    \end{align*}
In the typical fashion, we use the equation and integration by parts to get an estimate for $\partial_t F$:
    \begin{align*}
        \partial_t F(t) = \int_{\mathbb{R}^d} -\Delta(\langle \nabla, b \rangle) |u|^2 + 4 \langle \nabla\Bar{u}, (\nabla b )\nabla u \rangle -2\langle b , \nabla V\rangle |u|^2 .
    \end{align*}
By (c), 
    \begin{align*}
        \partial_t F(t) \geq    \int_{\mathbb{R}^d}4 \langle \nabla\Bar{u}, (\nabla b )\nabla u \rangle -2\langle b , \nabla V\rangle |u|^2 .
    \end{align*}
Now (e) implies 
    \begin{align*}
        \int_{\mathbb{R}^d}4 \langle \nabla\Bar{u}, (\nabla b )\nabla u \rangle -2\langle b , \nabla V\rangle |u|^2  \geq \int_{\mathbb{R}^d}4 \frac{|\nabla u|^2}{R \langle x/R\rangle^3} -2\langle b , \nabla V\rangle |u|^2.
    \end{align*}
For the second term in the sum, we note that $\inner{ x, \nabla V} = |x|\partial_{|x|}V$ and thus
    \begin{align*}
        -\inner{ x, \nabla V} = V-\partial_{|x|}(|x|V ) .
    \end{align*}
Since $V$ is a $(\kappa, \delta)$ potential, we can say that 
    \begin{align*}
        -\inner{ x, \nabla V} = V-\partial_{|x|}(|x|V ) \geq -\frac{2(\tfrac{d-2}{2})^2}{|x|^2}.
    \end{align*}
Therefore, 
    \begin{align*}
        -2\langle b , \nabla V\rangle= \frac{-2\inner{ x, \nabla V}}{R\inner{x/R}} \gtrsim -\frac{1}{R|x|^2}.
    \end{align*}
We now use the Fundamental Theorem of Calculus to conclude that 
    \begin{align*}
        \int_{\mathbb{R}}\int_{\mathbb{R}^d} \frac{|\nabla u|^2}{R \langle x/R\rangle^3} &\lesssim \lim_{T \to \infty}(F(T)-F(-T)) + \int_{\mathbb{R}} \int_{\mathbb{R}^d} \frac{|u|^2}{R|x|^2}\\
        &\lesssim \|u_0\|_{L^2} \|\nabla u_0\|_{L^2} + \frac{1}{R}\norm{|x|^{-1}u}^2_{L^2_{t, x}} . 
    \end{align*}
    By Corollary~\ref{cor BPSTZ resolvent} we can conclude that
    \begin{align*}
        \int_{\mathbb{R}}\int_{\mathbb{R}^d} \frac{|\nabla u|^2}{R \langle x/R\rangle^3} \lesssim \|u_0\|_{L^2} \|\nabla u_0\|_{L^2} + R^{-1}\|u_0\|_{L^2}^2. 
    \end{align*}
This also implies 
    \begin{align*}
        \int_{\mathbb{R}}\int_{B(0,R)} \frac{|\nabla u|^2}{R } \lesssim \|u_0\|_{L^2} \|\nabla u_0\|_{L^2} + R^{-1}\|u_0\|_{L^2}^2. 
    \end{align*}
Using $b_z(x):=b(x-z)$ in place of $b$ gives the second estimate of the lemma.
\bibliography{bibfile}

\begin{thebibliography}{10}

\bibitem{An08}
R.~Anton.
\newblock Cubic nonlinear {S}chr\"{o}dinger equation on three dimensional balls
  with radial data.
\newblock {\em Comm. Partial Differential Equations}, 33(10-12):1862--1889,
  2008.

\bibitem{BD07}
V.~Banica and T.~Duyckaerts.
\newblock Weighted {S}trichartz estimates for radial {S}chr\"{o}dinger equation
  on noncompact manifolds.
\newblock {\em Dyn. Partial Differ. Equ.}, 4(4):335--359, 2007.

\bibitem{Bour98}
J.~Bourgain.
\newblock Refinements of {S}trichartz' inequality and applications to
  {$2$}d-{NLS} with critical nonlinearity.
\newblock {\em Internat. Math. Res. Notices}, (5):253--283, 1998.

\bibitem{burq2004nonlinear}
N.~Burq, P.~G{\'e}rard, and N.~Tzvetkov.
\newblock On nonlinear schr{\"o}dinger equations in exterior domains.
\newblock {\em Annales de l'Institut Henri Poincar{\'e} C}, 21(3):295--318,
  2004.

\bibitem{BGT05}
N.~Burq, P.~G\'{e}rard, and N.~Tzvetkov.
\newblock Bilinear eigenfunction estimates and the nonlinear {S}chr\"{o}dinger
  equation on surfaces.
\newblock {\em Invent. Math.}, 159(1):187--223, 2005.

\bibitem{burq2004strichartz}
N.~Burq, F.~Planchon, J.~G. Stalker, and A.~S. Tahvildar-Zadeh.
\newblock Strichartz estimates for the wave and schr{\"o}dinger equations with
  potentials of critical decay.
\newblock {\em Indiana University mathematics journal}, pages 1665--1680, 2004.

\bibitem{Cam22}
N.~Camps.
\newblock Asymptotic stability of small ground states for {NLS} under random
  perturbations.
\newblock {\em Ann. Inst. H. Poincar\'e{} C Anal. Non Lin\'eaire}, 39(6):1--58,
  2022.

\bibitem{christ2001maximal}
M.~Christ and A.~Kiselev.
\newblock Maximal functions associated to filtrations.
\newblock {\em Journal of Functional Analysis}, 179(2):409--425, 2001.

\bibitem{CKSTT02}
J.~Colliander, M.~Keel, G.~Staffilani, H.~Takaoka, and T.~Tao.
\newblock Almost conservation laws and global rough solutions to a nonlinear
  {S}chr\"odinger equation.
\newblock {\em Math. Res. Lett.}, 9(5-6):659--682, 2002.

\bibitem{CKSTT08}
J.~Colliander, M.~Keel, G.~Staffilani, H.~Takaoka, and T.~Tao.
\newblock Global well-posedness and scattering for the energy-critical
  nonlinear {S}chr\"{o}dinger equation in {$\mathbb{R}^3$}.
\newblock {\em Ann. of Math. (2)}, 167(3):767--865, 2008.

\bibitem{CDKS01}
J.~E. Colliander, J.-M. Delort, C.~E. Kenig, and G.~Staffilani.
\newblock Bilinear estimates and applications to 2{D} {NLS}.
\newblock {\em Trans. Amer. Math. Soc.}, 353(8):3307--3325, 2001.

\bibitem{DPST07}
D.~De~Silva, N.~Pavlovi\'c, G.~Staffilani, and N.~Tzirakis.
\newblock Global well-posedness for a periodic nonlinear {S}chr\"odinger
  equation in 1{D} and 2{D}.
\newblock {\em Discrete Contin. Dyn. Syst.}, 19(1):37--65, 2007.

\bibitem{DSZ26}
M.~Deng, Y.~Song, and R.~Zhang.
\newblock Bilinear estimates for {S}chr\"odinger equation with repulsive
  inverse-square potential and its application.
\newblock {\em arXiv preprint arXiv:2609.15183}, 2026.

\bibitem{Deng23}
Y.~Deng.
\newblock On a bilinear restriction estimate for {S}chr\"odinger equations on
  2{D} waveguide.
\newblock {\em J. Differential Equations}, 416:1825--1836, 2025.

\bibitem{DDFZ24}
Y.~Deng, B.~Di, C.~Fan, and Z.~Zhao.
\newblock Bilinear estimate for {S}chr\"odinger equation on {$\Bbb{R}\times
  \Bbb{T}$}.
\newblock {\em J. Funct. Anal.}, 290(6):Paper No. 111312, 23, 2026.

\bibitem{DFYZZ24}
Y.~Deng, C.~Fan, K.~Yang, Z.~Zhao, and J.~Zheng.
\newblock On bilinear {S}trichartz estimates on waveguides with applications.
\newblock {\em J. Funct. Anal.}, 287(9):Paper No. 110595, 32, 2024.

\bibitem{Dod22}
B.~Dodson.
\newblock Global well-posedness for the defocusing, cubic nonlinear
  {S}chr\"{o}dinger equation with initial data in a critical space.
\newblock {\em Rev. Mat. Iberoam.}, 38(4):1087--1100, 2022.

\bibitem{engel2006short}
K.-J. Engel and R.~Nagel.
\newblock {\em A short course on operator semigroups}.
\newblock Springer, 2006.

\bibitem{FSWW}
C.~Fan, G.~Staffilani, H.~Wang, and B.~Wilson.
\newblock On a bilinear {S}trichartz estimate on irrational tori.
\newblock {\em Anal. PDE}, 11(4):919--944, 2018.

\bibitem{Gav26}
C.~Gavrus.
\newblock Bilinear {S}trichartz estimates on the hyperbolic space and the cubic
  {NLS}.
\newblock {\em arXiv preprint arXiv:2610.32358}, 2026.

\bibitem{GV85}
J.~Ginibre and G.~Velo.
\newblock The global {C}auchy problem for the nonlinear {S}chr\"odinger
  equation revisited.
\newblock {\em Ann. Inst. H. Poincar\'e{} Anal. Non Lin\'eaire}, 2(4):309--327,
  1985.

\bibitem{GLNY18}
Z.~Guo, J.~Li, K.~Nakanishi, and L.~Yan.
\newblock On the boundary {S}trichartz estimates for wave and {S}chr\"odinger
  equations.
\newblock {\em J. Differential Equations}, 265(11):5656--5675, 2018.

\bibitem{Hani12}
Z.~Hani.
\newblock A bilinear oscillatory integral estimate and bilinear refinements to
  {S}trichartz estimates on closed manifolds.
\newblock {\em Anal. PDE}, 5(2):339--363, 2012.

\bibitem{hernandez2018decomposition}
F.~Hernandez.
\newblock A decomposition for the schrodinger equation with applications to
  bilinear and multilinear estimates.
\newblock {\em Communications on Pure and Applied Analysis}, 17(2):627--646,
  2018.

\bibitem{Herr13}
S.~Herr.
\newblock The quintic nonlinear {S}chr\"odinger equation on three-dimensional
  {Z}oll manifolds.
\newblock {\em Amer. J. Math.}, 135(5):1271--1290, 2013.

\bibitem{HS15}
S.~Herr and N.~Strunk.
\newblock The energy-critical nonlinear {S}chr\"odinger equation on a product
  of spheres.
\newblock {\em Math. Res. Lett.}, 22(3):741--761, 2015.

\bibitem{HTT1}
S.~Herr, D.~Tataru, and N.~Tzvetkov.
\newblock Global well-posedness of the energy-critical nonlinear
  {S}chr\"{o}dinger equation with small initial data in {$H^1(\Bbb T^3)$}.
\newblock {\em Duke Math. J.}, 159(2):329--349, 2011.

\bibitem{HTT2}
S.~Herr, D.~Tataru, and N.~Tzvetkov.
\newblock Strichartz estimates for partially periodic solutions to
  {S}chr\"{o}dinger equations in {$4d$} and applications.
\newblock {\em J. Reine Angew. Math.}, 690:65--78, 2014.

\bibitem{IPT3}
A.~D. Ionescu and B.~Pausader.
\newblock The energy-critical defocusing {NLS} on {$\Bbb T^3$}.
\newblock {\em Duke Math. J.}, 161(8):1581--1612, 2012.

\bibitem{IPRT3}
A.~D. Ionescu and B.~Pausader.
\newblock Global well-posedness of the energy-critical defocusing {NLS} on
  {$\Bbb R\times \Bbb T^3$}.
\newblock {\em Comm. Math. Phys.}, 312(3):781--831, 2012.

\bibitem{Jao20}
C.~Jao.
\newblock Refined mass-critical {S}trichartz estimates for {S}chr\"odinger
  operators.
\newblock {\em Anal. PDE}, 13(7):1955--1994, 2020.

\bibitem{kato1989some}
T.~Kato and K.~Yajima.
\newblock Some examples of smooth operators and the associated smoothing
  effect.
\newblock {\em Reviews in Mathematical Physics}, 1(04):481--496, 1989.

\bibitem{KT98}
M.~Keel and T.~Tao.
\newblock Endpoint {S}trichartz estimates.
\newblock {\em Amer. J. Math.}, 120(5):955--980, 1998.

\bibitem{KPV96JAMS}
C.~E. Kenig, G.~Ponce, and L.~Vega.
\newblock A bilinear estimate with applications to the {K}d{V} equation.
\newblock {\em J. Amer. Math. Soc.}, 9(2):573--603, 1996.

\bibitem{KPV96TAMS}
C.~E. Kenig, G.~Ponce, and L.~Vega.
\newblock Quadratic forms for the {$1$}-{D} semilinear {S}chr\"odinger
  equation.
\newblock {\em Trans. Amer. Math. Soc.}, 348(8):3323--3353, 1996.

\bibitem{killip2017energy}
R.~Killip, C.~Miao, M.~Visan, J.~Zhang, and J.~Zheng.
\newblock The energy-critical nls with inverse-square potential.
\newblock {\em Discrete and Continuous Dynamical Systems}, 37(7):3831--3866,
  2017.

\bibitem{killip2018sobolev}
R.~Killip, C.~Miao, M.~Visan, J.~Zhang, and J.~Zheng.
\newblock Sobolev spaces adapted to the schr{\"o}dinger operator with
  inverse-square potential.
\newblock {\em Mathematische Zeitschrift}, 288(3):1273--1298, 2018.

\bibitem{KVClay}
R.~Killip and M.~Vi\c~san.
\newblock Nonlinear {S}chr\"odinger equations at critical regularity.
\newblock In {\em Evolution equations}, volume~17 of {\em Clay Math. Proc.},
  pages 325--437. Amer. Math. Soc., Providence, RI, 2013.

\bibitem{KV16}
R.~Killip and M.~Vi\c~san.
\newblock Scale invariant {S}trichartz estimates on tori and applications.
\newblock {\em Math. Res. Lett.}, 23(2):445--472, 2016.

\bibitem{LY24}
Z.~Lee and X.~Yu.
\newblock Global well-posedness and scattering for the defocusing septic
  one-dimensional {NLS} via new smoothing and almost {M}orawetz estimates.
\newblock {\em arXiv preprint arXiv:2407.07234}, 2024.

\bibitem{liskevich2003estimates}
V.~Liskevich and Z.~Sobol.
\newblock Estimates of integral kernels for semigroups associated with
  second-order elliptic operators with singular coefficients.
\newblock {\em Potential Analysis}, 18(4):359--390, 2003.

\bibitem{milman2004global}
P.~D. Milman and Y.~A. Semenov.
\newblock Global heat kernel bounds via desingularizing weights.
\newblock {\em Journal of Functional Analysis}, 212(2):373--398, 2004.

\bibitem{mizutani2020global}
H.~Mizutani and X.~Yao.
\newblock Global {K}ato smoothing and {S}trichartz estimates for
  {S}chr\"odinger type equations with rough decay potentials.
\newblock {\em Rev. Math. Phys.}, 37(7):Paper No. 2450050, 33, 2025.

\bibitem{muscalu2013classical}
C.~Muscalu and W.~Schlag.
\newblock {\em Classical and Multilinear Harmonic Analysis: Volume 1}, volume
  137.
\newblock Cambridge University Press, 2013.

\bibitem{ozawa2013sharp}
T.~Ozawa and K.~M. Rogers.
\newblock Sharp morawetz estimates.
\newblock {\em Journal d'Analyse Math{\'e}matique}, 121(1):163--175, 2013.

\bibitem{PTW14}
B.~Pausader, N.~Tzvetkov, and X.~Wang.
\newblock Global regularity for the energy-critical {NLS} on {$\mathbb{S}^3$}.
\newblock {\em Ann. Inst. H. Poincar\'{e} C Anal. Non Lin\'{e}aire},
  31(2):315--338, 2014.

\bibitem{PS24}
F.~Pusateri and A.~Soffer.
\newblock Bilinear estimates in the presence of a large potential and a
  critical {NLS} in 3{D}.
\newblock {\em Mem. Amer. Math. Soc.}, 299(1498):v+107, 2024.

\bibitem{Schippa24}
R.~Schippa.
\newblock Improved global well-posedness for mass-critical nonlinear
  {S}chr\"odinger equations on tori.
\newblock {\em J. Differential Equations}, 412:87--139, 2024.

\bibitem{ST25}
R.~Schippa and D.~Tataru.
\newblock Wave packet decompositions and sharp bilinear estimates for rough
  hamiltonian flows.
\newblock {\em arXiv preprint arXiv:2509.23551}, 2025.

\bibitem{Shao09}
S.~Shao.
\newblock Sharp linear and bilinear restriction estimates for paraboloids in
  the cylindrically symmetric case.
\newblock {\em Rev. Mat. Iberoam.}, 25(3):1127--1168, 2009.

\bibitem{Str77}
R.~S. Strichartz.
\newblock Restrictions of {F}ourier transforms to quadratic surfaces and decay
  of solutions of wave equations.
\newblock {\em Duke Math. J.}, 44(3):705--714, 1977.

\bibitem{Strunk15}
N.~Strunk.
\newblock Global well-posedness of the energy-critical defocusing {NLS} on
  rectangular tori in three dimensions.
\newblock {\em Differential Integral Equations}, 28(11-12):1069--1084, 2015.

\bibitem{SY20d}
M.~Sy and X.~Yu.
\newblock Global well-posedness and long-time behavior of the fractional {NLS}.
\newblock {\em Stoch. Partial Differ. Equ. Anal. Comput.}, 10(4):1261--1317,
  2022.

\bibitem{SY22}
M.~Sy and X.~Yu.
\newblock Global well-posedness for the cubic fractional {NLS} on the unit
  disk.
\newblock {\em Nonlinearity}, 35(4):2020--2072, 2022.

\bibitem{tao2010physical}
T.~Tao.
\newblock A physical space proof of the bilinear strichartz and local smoothing
  estimates for the {S}chr{\"o}dinger equation.

\bibitem{Tao03}
T.~Tao.
\newblock A sharp bilinear restrictions estimate for paraboloids.
\newblock {\em Geom. Funct. Anal.}, 13(6):1359--1384, 2003.

\bibitem{tao2006nonlinear}
T.~Tao.
\newblock {\em Nonlinear dispersive equations: local and global analysis}.
\newblock Number 106. American Mathematical Soc., 2006.

\bibitem{TiceInterpolation}
I.~Tice.
\newblock A crash course in interpolation theory.
\newblock Lecture notes, Department of Mathematical Sciences, Carnegie Mellon
  University, 2025.

\bibitem{Vil07}
M.~C. Vilela.
\newblock Inhomogeneous {S}trichartz estimates for the {S}chr\"{o}dinger
  equation.
\newblock {\em Trans. Amer. Math. Soc.}, 359(5):2123--2136, 2007.

\bibitem{WY23}
B.~Wilson and X.~Yu.
\newblock Global well-posedness and scattering for the defocusing mass-critical
  {S}chr{\"o}dinger equation in the three-dimensional hyperbolic space.
\newblock {\em Trans. Amer. Math. Soc.}, 2026.

\bibitem{Yaj87}
K.~Yajima.
\newblock Existence of solutions for {S}chr\"odinger evolution equations.
\newblock {\em Comm. Math. Phys.}, 110(3):415--426, 1987.

\bibitem{YY20}
X.~Yu and H.~Yue.
\newblock On the {G}lobal {W}ell-{P}osedness for the {P}eriodic {Q}uintic
  {N}onlinear {S}chr\"{o}dinger {E}quation.
\newblock {\em SIAM J. Math. Anal.}, 56(2):1851--1902, 2024.

\bibitem{Yue}
H.~Yue.
\newblock Global well-posedness for the energy-critical focusing nonlinear
  {S}chr\"{o}dinger equation on {$\Bbb T^4$}.
\newblock {\em J. Differential Equations}, 280:754--804, 2021.

\end{thebibliography}
\bibliographystyle{abbrv}
\end{document}